\documentclass[oneside,12pt,a4paper]{article}
\usepackage[T1]{fontenc}
\usepackage[utf8]{inputenc}
\usepackage[top=1cm,bottom=2cm,right=2cm,left=2cm,includehead]{geometry}
\usepackage{amsmath,amsfonts,amssymb,amscd,amsthm}
\usepackage[usenames,dvipsnames,svgnames]{xcolor}
\usepackage[numbers]{natbib}
\usepackage{doi}
\usepackage{graphicx}
\usepackage[font={small}]{caption} 
\usepackage[font={small}]{subcaption}
\usepackage{autonum} 
\usepackage{hyperref} 
\usepackage{enumitem}

\newcommand{\thetitle}{\uppercase{Distributed-memory parallelization of the aggregated unfitted  finite element method}}
\newcommand{\theauthors}{Francesc Verdugo\textsuperscript{a,}\footnote{Corresponding author. 
 \\Emails: \theemails},
 Alberto F. Mart\'in\textsuperscript{a,b}, and Santiago Badia\textsuperscript{a,b}}
\newcommand{\theauths}{F. Verdugo, A. F. Mart\'in, and S. Badia}
\newcommand{\theaffiliations}{
\textsuperscript{a}CIMNE – Centre Internacional de M\`etodes Num\`erics en
Enginyeria\\ Esteve Terradas 5, 08860
Castelldefels, Spain.\\[0.5em]

\textsuperscript{b}Department of Civil and Environmental Engineering, Universitat Polit\`ecnica de Catalunya\\ Jordi Girona 1-3, Edifici C1, 08034 Barcelona, Spain.
}
\newcommand{\theemails}{\texttt{fverdugo@cimne.upc.edu} (FV), \texttt{amartin@cimne.upc.edu} (AM), \texttt{sbadia@cimne.upc.edu} (SB)\\ \today }

\newcommand{\thethanks}{Financial support from the European Commission under the FET-HPC ExaQUte project (Grant agreement ID: 800898) within the Horizon 2020 Framework Programme is gratefully acknowledged. This work has been partially funded by the project MTM2014-60713-P from the ``Ministerio de Econom\'ia, industria y Competitividad'' of Spain. S. Badia gratefully acknowledges the support received from the Catalan Government through the ICREA Acad\`emia Research Program. F. Verdugo gratefully acknowledges the support received from the  \emph{Secretaria d'Universitats i Recerca} of the Catalan Government in the framework of the Beatriu Pinós Program (Grant Id.: 2016 BP 00145). The authors thankfully acknowledge the computer resources at Marenostrum-IV and the technical support provided by the Barcelona Supercomputing Center (RES-ActivityID: FI-2018-1-0014, FI-2018-2-0009, FI-2018-3-0029, FI-2019-1-0007).}

\hypersetup{
    bookmarks=true,
    breaklinks=true,
    bookmarksopen=true,
    pdftitle={\thetitle},    
    pdfauthor={\theauths},     
    colorlinks=true,       
    linkcolor=black,          
    citecolor=blue,        
    filecolor=black,      
    urlcolor=blue           
}

\usepackage{sectsty}
\sectionfont{\centering\sc\normalsize\selectfont}
\subsectionfont{\normalsize\selectfont}
\subsubsectionfont{\normalsize\selectfont}
\usepackage[compact]{titlesec}

\usepackage{fancyhdr}

\usepackage{newtxtext}
\usepackage{newtxmath}

\usepackage{acronym}
\usepackage[most]{tcolorbox}

\usepackage{booktabs}

\usepackage{tikz}
\definecolor{myellow}{RGB}{255,230,128}
\definecolor{gray20}{RGB}{204,204,204}
\definecolor{mygray}{RGB}{204,204,204}
\definecolor{mygreen}{RGB}{138,203,95}
\definecolor{myblue}{RGB}{77,151,214}

\definecolor{lstgrey}{rgb}{0.95,0.95,0.95}
\usepackage{listings}
\usepackage[textsize=small]{todonotes}

\usepackage[vlined]{algorithm2e}
\usepackage{xcolor}

\SetAlFnt{\small}
\SetArgSty{textnormal}

\SetAlCapSty{xAlCapSty}

\SetCommentSty{xCommentSty}

\SetNlSty{mynlfont}{}{} 

\LinesNumbered

\SetSideCommentRight

\DontPrintSemicolon

\RestyleAlgo{algoruled}

\SetInd{0.5em}{0.5em}

\acrodef{fe}[FE]{finite element}
\acrodef{dof}[DOF]{degree of freedom}
\acrodefplural{dof}[DOFs]{degrees of freedom}
\acrodef{vef}[VEF]{vertex, edge, and face}
\acrodefplural{vef}[VEFs]{vertices, edges, and faces}
\acrodef{dg}[DG]{discontinuous Galerkin}
\acrodef{vms}[VMS]{variational multiscale}
\acrodef{sps}[SPS]{symmetric projection stabilization}
\acrodef{agfem}[agFEM]{aggregated finite element method}
\acrodef{xfem}[XFEM]{extended finite element method}
\acrodef{agfe}[agFE]{aggregated finite element}
\acrodef{pde}[PDE]{partial differential equation}
\acrodef{amg}[AMG]{algebraic multigrid}
\acrodef{agg}[AgFEM]{aggregated unfitted finite element method}

\newcommand{\rev}{\color{black}}

\newcommand{\closure}[2][3]{%
{}\mkern#1mu\overline{\mkern-#1mu#2}}

\graphicspath{{figures/}{figures-rev1/}}

\def\FEMPAR{{\texttt{FEMPAR}}}
\def\dealii{{\texttt{deal.II}}}
\def\p4est{{\texttt{p4est}}}
\def\t8code{{\texttt{t8code}}}

\def\petsc{{\texttt{PETSc}}}
\def\trilinos{{\texttt{TRILINOS}}}

\begin{document}

\thispagestyle{empty}

%
%
%
%



\renewcommand*{\thefootnote}{\fnsymbol{footnote}}

\begin{center}
{ \bf {\thetitle}}

\vspace*{1em}

\theauthors

\vspace*{1em}

\theaffiliations

\end{center}

\setcounter{footnote}{0}
\renewcommand*{\thefootnote}{\arabic{footnote}}



\begin{center}

{\bf Abstract}

\vspace*{1em}

\begin{minipage}{0.9\textwidth}
\begin{small}

The aggregated unfitted finite element method (AgFEM) is a methodology recently introduced in order to address conditioning and stability problems associated with embedded, unfitted, or extended finite element methods. The method is based on removal of basis functions associated with badly cut cells by introducing carefully designed constraints, which results in well-posed systems of linear algebraic equations, while preserving the optimal approximation order of the underlying finite element spaces. The specific goal of this work is to present the implementation and performance of the method on distributed-memory platforms aiming at the efficient solution of large-scale problems.  In particular, we show that, by considering AgFEM, the resulting systems of linear algebraic equations can be effectively solved using standard algebraic multigrid preconditioners. This is in contrast with previous works that consider highly customized preconditioners in order to allow one the usage of iterative solvers in combination with unfitted techniques. Another novelty with respect to the methods available in the literature is the problem sizes that can be handled with the proposed approach. While most of previous references discussing linear solvers for unfitted methods are based on serial non-scalable algorithms, we propose a parallel distributed-memory method able to efficiently solve problems at large scales. This is demonstrated by means of a weak scaling test defined on complex 3D domains  up to 300M degrees of freedom and one billion cells on 16K CPU cores in the Marenostrum-IV platform. The parallel implementation of the AgFEM method is available in the large-scale finite element package \FEMPAR.

\end{small}

\end{minipage}
\end{center}

\vspace*{1em}






\section{Introduction}\label{sec:int}

Unfitted \ac{fe} methods (also known as embedded or immersed \ac{fe} methods) are useful techniques in order to perform numerical computations  associated with geometrically complex domains. They are specially useful in multi-phase and multi-physics applications with moving interfaces (e.g., fracture mechanics \cite{Sukumar2000}, fluid–structure interaction \cite{Massing2015}, free surface flows \cite{Sauerland2011}), in applications with varying domain topologies (e.g., shape or topology optimization \cite{Burman2018},  3D printing simulations \cite{chiumenti_numerical_2017}), in applications where the geometry is not described by CAD data (e.g., medical simulations based on CT-scan images \cite{Nguyen2017}), or in large-scale parallel computations, where generating and partitioning large unstructured meshes is particularly difficult. The main benefit of unfitted \ac{fe} methods is that they do not require the generation of \emph{body-fitted} meshes (i.e., meshes whose faces conform to the domain boundary). Instead, they embed the domain of interest in a geometrically simple background grid (usually a uniform or an adaptive Cartesian grid), which can be generated much more efficiently. The popularity of this approach is illustrated by the large number of different methodologies that have been proposed following  this rationale (see, e.g., the cutFEM method \cite{burman_cutfem:_2015}, the cutIGA method \cite{Elfverson2018}, the Finite Cell Method \cite{Schillinger2015}, the AgFEM method \cite{Badia2018}, the cg-FEM method \cite{Nadal2013}, the Immersed Boundary Method \cite{Mittal2005}, the i-spline method \cite{sanches_immersed_2011},  and some variants of the XFEM method \cite{sukumar_modeling_2001}).

Unfortunately, unfitted \ac{fe} methods have also well known drawbacks. One of the most notorious {\rev ones} is the so-called \emph{small cut cell problem}. The intersection of a background cell with the physical domain can be arbitrarily small and with arbitrarily high aspect ratios, which usually results in severely ill-conditioned systems of algebraic linear equations if no specific strategy is used to remedy it \cite{DePrenter2017}. This flaw makes the solution of the underlying linear systems much more challenging than for standard \ac{fe} methods based on body-fitted grids and it is still today one of the main limiting factors for the successful application of unfitted \ac{fe} techniques in realistic large-scale applications. Sparse direct solvers \cite{Davis2006} are usually considered when dealing with unfitted \ac{fe} techniques (see, e.g., \cite{Schillinger2015}) since they are robust enough to deal (up to a certain extent) with such ill-conditioned problems. However, its usage is prohibitive at large scales since their memory footprint and the algorithmic complexity scales supra-linearly with respect {\rev to} the problem size. In \ac{fe} analysis, the current way to effectively solve linear systems at large-scales is using iterative Krylov sub-space methods \cite{saad_iterative_2003} in combination with parallel and scalable preconditioners. Unfortunately, the well known scalable preconditioners based on multigrid \cite{Briggs2000,Chow2006} (either geometric \cite{Wesseling2001} or algebraic \cite{Ruge1987,Vanek2001}), or multi-level domain decomposition \cite{badia_multilevel_2016,Toselli2005} are mainly designed for body-fitted meshes and cannot readily deal with the ill-conditioning associated with the small cut cell problem.

Specific preconditioners for unfitted methods have been proposed in the literature in different contexts. These techniques are generally based on the idea of splitting the system matrix into blocks associated with    problematic \acp{dof}, i.e., the ones affected by small cuts, and unproblematic ones, i.e., the ones not affected (see, e.g., \cite{berger-vergiat_inexact_2012,DePrenter2017,hiriyur_quasi-algebraic_2012,Jomo2018,menk_robust_2011}). The methods differ in the particular techniques used to handle the block associated with non-problematic \acp{dof}, the block of problematic \acp{dof}, and the coupling between these two blocks. Even though these methodologies have been able to allow the usage of iterative solvers in combination with unfitted \ac{fe} methods in particular settings, they are mainly serial non-scalable algorithms. Considered problems are typically academic examples that go up to few hundreds of thousands of \acp{dof}. To the authors' best knowledge, the only exceptions are the works in \cite{badia_robust_2017,Jomo2018}. Reference \cite{Jomo2018} is based on a single level Additive Schwarz preconditioner and reports results for problems up to almost 4 million \acp{dof} solved with a hybrid OpenMP+MPI implementation with up to 16 MPI tasks and 64 OpenMP threads per MPI task. Even though this method is a step forward with respect to most previous works, it is well known that single-level Additive Schwarz preconditioners do not scale up algorithmically, and thus the method is not suitable for larger scales. On the other hand, the work in \cite{badia_robust_2017} considers a two-level balancing domain decomposition preconditioner \cite{Mandel2003} and reports results up to almost 6 million \acp{dof}. 

The goal of this work is to scale up to much larger problem sizes by considering an alternative approach. Instead of using a tailored preconditioner  able to deal with linear systems affected by the small cut cell problem,  we consider an enhanced \ac{fe} formulation so that the resulting method  leads to system matrices that are not affected by small cuts. This opens the door to consider well known scalable algorithms for conventional \ac{fe} analysis such as \ac{amg} \cite{Ruge1987,Vanek2001}, for which there are highly scalable parallel implementations in renowned scientific computing packages such as \trilinos~\cite{Heroux2005} or \petsc~\cite{petsc-web-page,petsc-user-ref}. This makes the method easy to use in parallel runs since the development of tailored linear solvers is not required.
To this end, we consider the recently introduced \ac{agg} \cite{Badia2018}. Other enhanced unfitted formulations could be also taken into account (e.g., the CutFEM method \cite{burman_cutfem:_2015}), but this task is not in the scope of the current work.
\ac{agg} is based on removal of basis functions associated with
badly cut cells by introducing carefully designed constraints and its formulation shares the good properties of body-fitted \ac{fe} methods such as stability, condition number bounds, optimal convergence, and continuity with respect to data. The detailed mathematical analysis of the method is included in \cite{Badia2018} for elliptic problems and in \cite{Badia2018a} for the Stokes equation. Even though the method is potentially applicable to parallel large-scale computations, previous works \cite{Badia2018,Badia2018a} focus mainly in the description of the methodology  in a serial context and its mathematical analysis. The main contribution of the current work is this two-fold. On the one hand, we extend \ac{agg} to large-scale computations by presenting an efficient distributed-memory implementation, which is based on standard functionality provided by distributed-memory \ac{fe} packages. In this case, we have implemented the method in the large-scale \ac{fe} software package \FEMPAR~ \cite{badia_fempar:_2017}.
The algorithmic strategies proposed for \ac{agg} in distributed-memory environments can also be used for the parallelization of other cell aggregation techniques for \ac{dg} methods \cite{Helzel2005,Johansson2013,Kummer2017}. It will also be essential in CutFEM \cite{burman_cutfem:_2015}, when using substructuring domain decomposition techniques to end up with well-posed subdomain problems.
On the other hand,
we show that, by considering the parallel \ac{agg}, the resulting systems of linear algebraic equations can be effectively solved at large scales using standard \ac{amg} preconditioners. Here, we have considered the well-known \petsc{}'s native \ac{amg} framework, referred to as \texttt{GAMG} \cite{GAMGweb}. We have customized as less as possible the parameters of this solver in order to demonstrate that users can take it as it is given ``out of the box'' by the library. These properties are demonstrated with weak scaling tests up to 16K cores and up to nearly 300M \acp{dof} in the Marenostrum-IV platform at the Barcelona Supercomputing Center \cite{bsc-home} using the Poisson equation on complex 3D domains as a model problem.

The outline of the paper is as follows. In Sect. \ref{sec:serial-aggFEM}, we briefly present the main components of the \ac{agg} method in a serial context as a basis for {\rev the subsequent} extension to distributed-memory platforms. In particular, we introduce the geometrical setup (Sects.~\ref{sec:imb_setup} and \ref{sec:cell-aggr}), the \ac{fe} spaces of the \ac{agg} method (Sect.~\ref{sec:fe-spaces}), the associated shape function basis (Sect. \ref{sec:shape-funs}), and the \ac{fe} assembly (Sect.~\ref{sec:serial-assembly}). In Sect.~\ref{sec:dm-impl}, we discuss the distributed-memory implementation. {\rev We detail the domain-decomposition setup (Sect.~\ref{sec:dd-setup}), the distributed \ac{fe} spaces (Sect.~\ref{sec:par-fe-spaces})}, several important algorithmic components (Sects.~\ref{sec:overview}, \ref{sec:par-cell-agg}, \ref{sec:data-import}, and \ref{sec:data-layout}), and the parallel \ac{fe} assembly (Sect.~\ref{sec:par-assembly}). In Sect.~\ref{sec:numericals}, we report the numerical examples. {\rev We detail the HPC environment (Sect.~\ref{sec:experimental-environment})}, the considered model problem (Sect.~\ref{sec:model-problem}), the considered linear solver (Sect.~\ref{sec:ls-setup}), the particular setup for the experiments (Sect.~\ref{sec:problem-setup}), and the results of the weak scaling test {\rev (Sects.~\ref{sec:weak-scaling-test} and \ref{sec:weak-scaling-test-full})}. Finally, we close the paper with the conclusions in Sect.~\ref{sec:conclussions}.

\section{The aggregated unfitted finite element method}  \label{sec:serial-aggFEM}

\subsection{Immersed boundary setup} \label{sec:imb_setup}

\begin{figure}[ht!]
  \centering
  \begin{subfigure}{0.24\textwidth}
    \centering
    \includegraphics[width=0.9\textwidth]{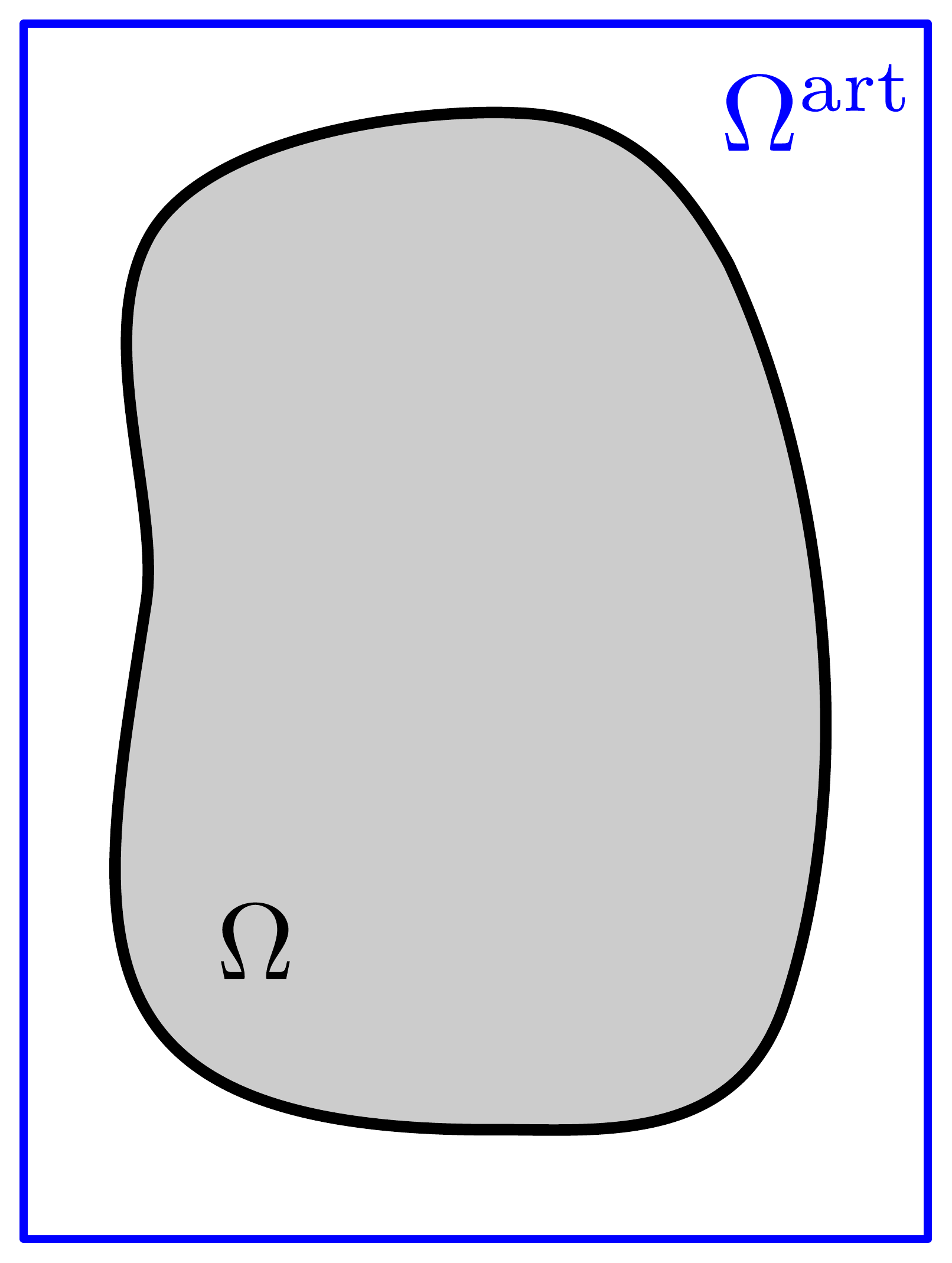}
    \caption{}
    \label{fig:immersed-setup-a}
  \end{subfigure}
  \begin{subfigure}{0.24\textwidth}
    \centering
    \includegraphics[width=0.9\textwidth]{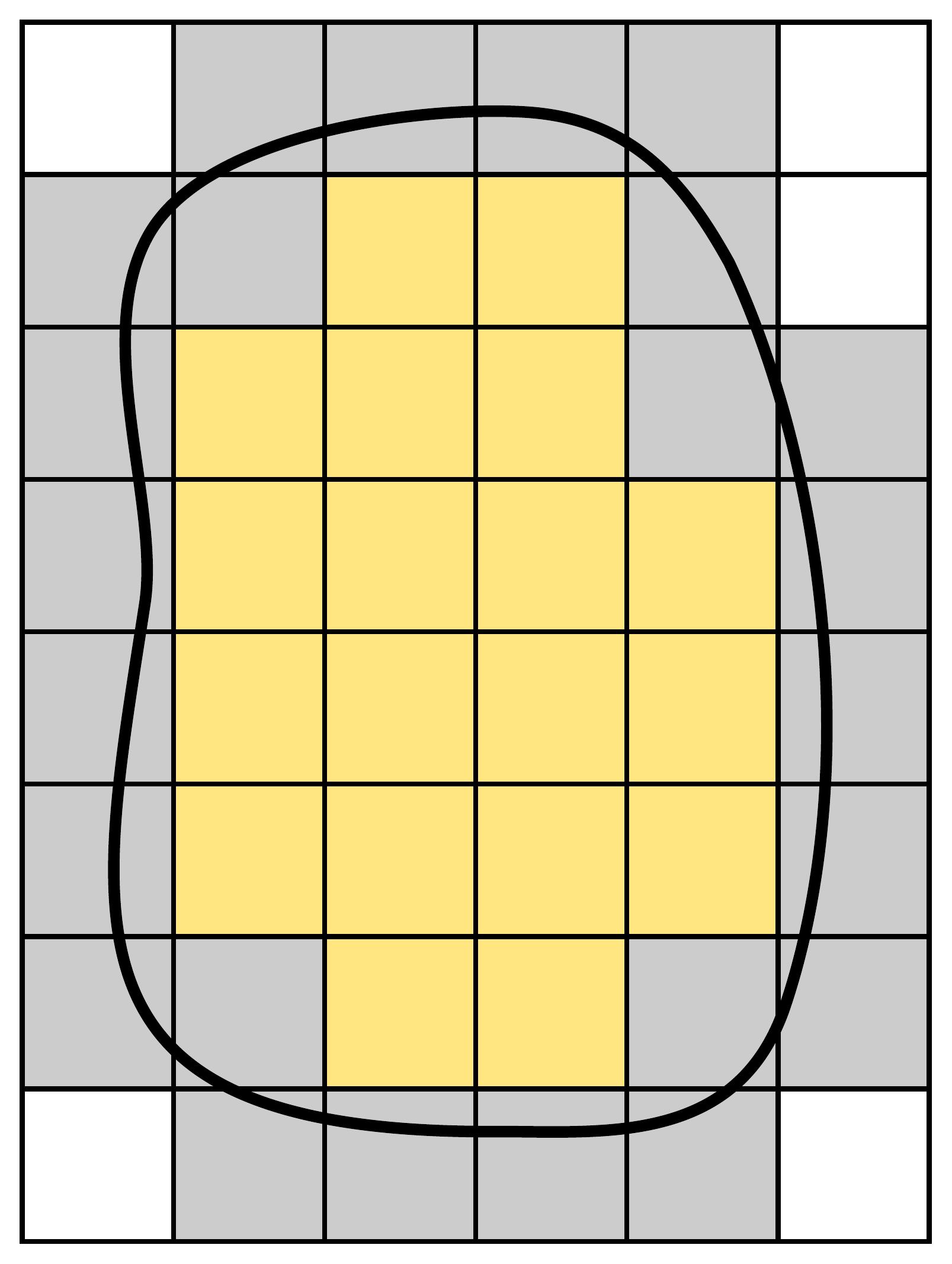}
    \caption{}
    \label{fig:immersed-setup-b}
  \end{subfigure}
    \begin{subfigure}{0.1\textwidth}
    \begin{tabular}{l}
    \tikz{\draw[fill=myellow]  (0,0) rectangle (1.2em,1.2em);} \small \emph{internal} cells
    \\
    \tikz{\draw[fill=gray20]  (0,0) rectangle (1.2em,1.2em);} \small \emph{cut} cells
    \\
    \tikz{\draw  (0,0) rectangle (1.2em,1.2em);} \small \emph{external} cells
    \end{tabular}
  \end{subfigure}
  \caption{Embedded boundary setup.}
  \label{fig:immersed-setup}
\end{figure}

Let $\Omega \subset \mathbb{R}^d$ be an open bounded polygonal domain with $d\in\{2,3\}$ the number of spatial dimensions.  Like in any other embedded boundary method, let us introduce an \emph{artificial} domain $\Omega^{\rm art}$ such that it has a simple geometry (e.g., a cuboid) and it includes the \emph{physical} domain $\Omega \subseteq\Omega^{\rm art}$  (see Fig. ~\ref{fig:immersed-setup-a}).  We denote by $\mathcal{T}^{\rm art}$ a partition of $\Omega^{\rm art}$ into \emph{cells} with characteristic cell size $h$.   We build $\mathcal{T}^{\rm art}$ as a Cartesian grid of $d$-dimensional parallelepipeds, i.e., hexahedra for $d=3$ and quadrilaterals for $d=2$ (see Fig. ~\ref{fig:immersed-setup-b}). We consider Cartesian grids since they can be efficiently generated and partitioned in parallel with the \p4est library \cite{burstedde_p4est_2011}. Note, however, that the following rationale holds also for unstructured background meshes and other cell types such as triangles and tetrahedra.

We introduce the usual classification of cells in  $\mathcal{T}^{\rm art}$ as \emph{internal, external or cut}. Specifically, $T \in \mathcal{T}^{\rm art}$ is an \emph{internal cell} if $T \subseteq \Omega$; $T$ is an \emph{external cell} if $T \cap \Omega = \emptyset$; otherwise, $T$ is a \emph{cut cell} (see Fig.~\ref{fig:immersed-setup-b}).
For the sake of simplicity and without any loss of generality, we represent the domain boundary $\partial\Omega$ as the zero level-set of a given scalar function $\psi^{\rm ls}$, namely $\partial\Omega\doteq\{x\in\mathbb{R}^d:\psi^{\rm ls}(x)=0\}$.  We define the physical domain as the set of points where the levelset function is negative, namely  $\Omega\doteq\{x\in\mathbb{R}^d:\psi^{\rm ls}(x)<0\}$. In this setup, classifying the cells into internal, external or cut is straightforward by sampling the levelset function at cell nodes. We note that the problem geometry could be described by other means, e.g., using 3D CAD data, by providing techniques to perform the cell classification and to perform numerical integration in cut cells (see, e.g., \cite{marco_exact_2015}). In any case, the way the geometry is handled does not affect the following exposition.
We represent the set of interior (resp., external and cut) cells with $\mathcal{T}^{\rm in}$ (resp., $\mathcal{T}^{\rm out}$ and $\mathcal{T}^{\rm cut}$). Furthermore, we define the set of \emph{active cells} as the union of the sets of interior and cut cells, namely $\mathcal{T}^{\rm act} \doteq \mathcal{T}^{\rm in} \cup \mathcal{T}^{\rm cut}$. The domain associated with active cells, namely $\Omega^{\rm act}\doteq\left(\bigcup_{T\in\mathcal{T}^{\rm act}}T\right) \cup ( \bigcup_{T,T'\in\mathcal{T}^{\rm act}} \closure[1]{T} \cap \closure[1]{T }')  $ is referred to as the \emph{active} domain, which satisfies $\Omega\subseteq\Omega^{\rm act}\subseteq\Omega^{\rm art}$. For simplicity, we refer to the set of active cells $\mathcal{T}^{\rm act}$ simply as $\mathcal{T}$ (i.e., we drop the label "act" in the notation).
In the following developments, it is convenient to assign a unique global identifier (referred to  as \emph{global id}) to each cell in $\mathcal{T}$. We denote by $T^k$ the cell with the global id $k$, for $k=1,\ldots,|\mathcal{T}|$, where $|\cdot|$ is the number of elements in a set. With this notation, we can write $\mathcal{T}=\{T^k\}_{k\in\mathcal{K}}$, where $\mathcal{K}\doteq\{1,\ldots,|\mathcal{T}|\}$ is the index set of all global cell ids. We also introduce $\mathcal{K}^{\rm in}$ (resp. $\mathcal{K}^{\rm cut}$) the subset of $\mathcal{K}$ that contains the ids of all interior (resp. cut) cells, namely $\mathcal{K}^{\rm in}\doteq\{k\in\mathcal{K}: T^k\in\mathcal{T}^{\rm in}\}$ and $\mathcal{K}^{\rm cut}\doteq\{k\in\mathcal{K}: T^k\in\mathcal{T}^{\rm cut}\}$.

\subsection{Cell aggregation}\label{sec:cell-aggr}

We briefly present here the (serial) cell aggregation strategy proposed in \cite{Badia2018}, which is the cornerstone of the \ac{agg} method. 
The extension of this cell aggregation method to  parallel distributed-memory frameworks is discussed later in Sect.~\ref{sec:par-cell-agg}.  The cell aggregation process consists in building a surjective map that for each cell in $\mathcal{T}$ returns an interior cell in $\mathcal{T}^{\rm in}$. For convenience, we represent this map as a transformation $R:\mathcal{K}\rightarrow\mathcal{K}^{\rm in}$ between cell ids. That is, $R$ takes a cell id $k\in\mathcal{K}$ and returns cell id $R(k)\in\mathcal{K}^{\rm in}$ corresponding to an interior cell. The transformation $R$ is called the \emph{root cell map} since, for a given cell $T^k\in\mathcal{T}$, the associated interior cell $T^{R(k)}\in\mathcal{T}^{\rm in}$  is referred to as the \emph{root} cell of $T^k$. The concept of root cell plays an important role in the construction of the so-called \ac{agfe} spaces, which will be discussed later in Sect~\ref{sec:fe-spaces}.

The root cell map $R$ induces a partition of the cells in $\mathcal{T}$ into so-called \emph{cell aggregates}. For each interior cell $T^k\in\mathcal{T}^{\rm in}$, we define a cell aggregate $\mathcal{A}_k$ as the set of cells that share $T^k$ as their root cell. Formally,  the aggregate associated with cell $T^k$ is defined as $\mathcal{A}_k\doteq \{T^{k'}\}_{k'\in R^{-1}(k)}$, where the inverse of $R$ is defined as  $R^{-1}(k)\doteq\{k'\in\mathcal{K}:\ k=R(k')\}$. Clearly, the set of aggregates $\{\mathcal{A}_k\}_{k\in\mathcal{K}^{\rm in}}$ is a partition of $\mathcal{T}$ since each cell in $\mathcal{T}$ belongs to one (and only one) aggregate.

As stated in \cite{Badia2018}, in order to construct \ac{agg} spaces, the map $R$ has to be such that the resulting cell aggregates $\mathcal{A}_k$ satisfy the following two properties for each $k\in\mathcal{K}^{\rm in}$: (a) the set $\mathcal{A}_k$ contains one and only one interior cell, namely $|\mathcal{A}_k\cap\mathcal{T}^{\rm in}|=1$; and (b) the (unique) interior cell $\tilde T_0\in \mathcal{A}_k\cap\mathcal{T}^{\rm in}$ can be reached from any cut cell $\tilde T_n\in \mathcal{A}_k\cap\mathcal{T}^{\rm cut}$  in the aggregate following a cell path  $\tilde T_0$, $\tilde T_1$, $\ldots$, $\tilde T_n$, where all cells are in the aggregate, i.e., $\tilde T_i\in \mathcal{A}_k$, $i=0,\ldots,n$, and two consecutive cells, $\tilde T_{i-1}$ and $\tilde T_i$, share a non-external face (for $d=3$) or a non-external edge (for $d=2$) for any $i=1,\ldots,n$. More precisely, we say that $\tilde T_{i-1}$ and $\tilde T_i$ share a non-external face ($d=3$) or a non-external edge ($d=2$) if $\tilde T_{i-1}\cap\tilde T_i$ is a manifold in $\mathbb{R}^d$ with co-dimension~1 and $\tilde T_{i-1}\cap\tilde T_i\cap\Omega\neq\emptyset$. Note that sharing only a corner (or an edge for $d=3$) is not enough. For simplicity, we will use the words "edge" and "face" as completely equivalent synonyms in the following for the  two-dimensional case ($d=2$).

\begin{figure}[ht!]
  \centering
  \begin{subfigure}{0.99\textwidth}
    \centering
    \begin{small}
      \begin{tabular}{llll}
         \tikz{\fill[fill=myellow]  (0,0) rectangle (1.4em,1.4em);} touched
         &
         \tikz{\fill[fill=gray20]    (0,0) rectangle (1.4em,1.4em);} untouched
         &
         \tikz{ \draw[line width=2pt] (0,0) -- (2em,0);} $\partial \Omega$
      \end{tabular}
    \end{small}
  \end{subfigure}
  \par
  \begin{subfigure}{0.24\textwidth}
    \centering
    \includegraphics[width=0.9\textwidth]{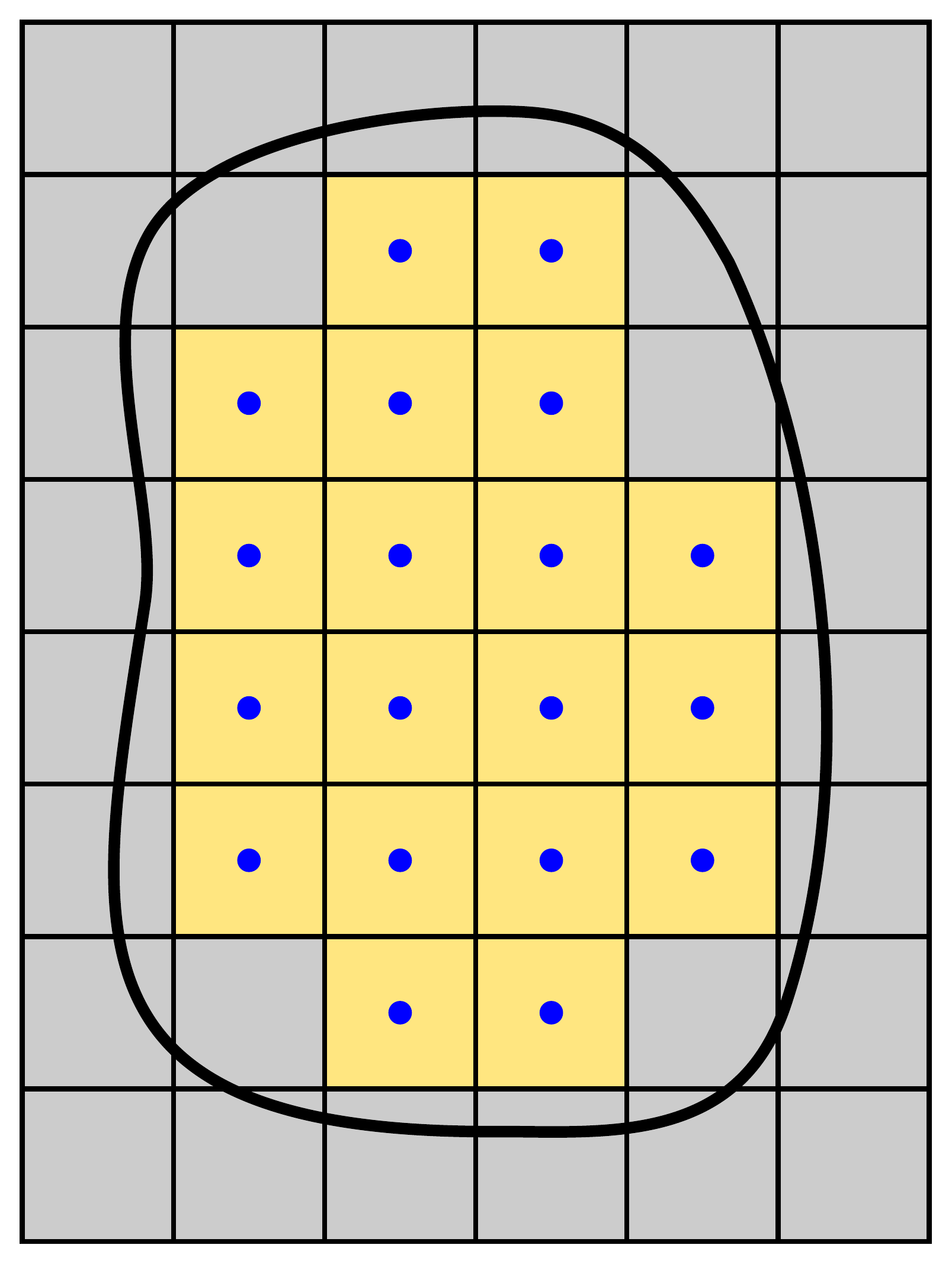}
    \caption{ Step 1.}
    \label{fig:aggr-steps-a}
  \end{subfigure}
  \begin{subfigure}{0.24\textwidth}
    \centering
    \includegraphics[width=0.9\textwidth]{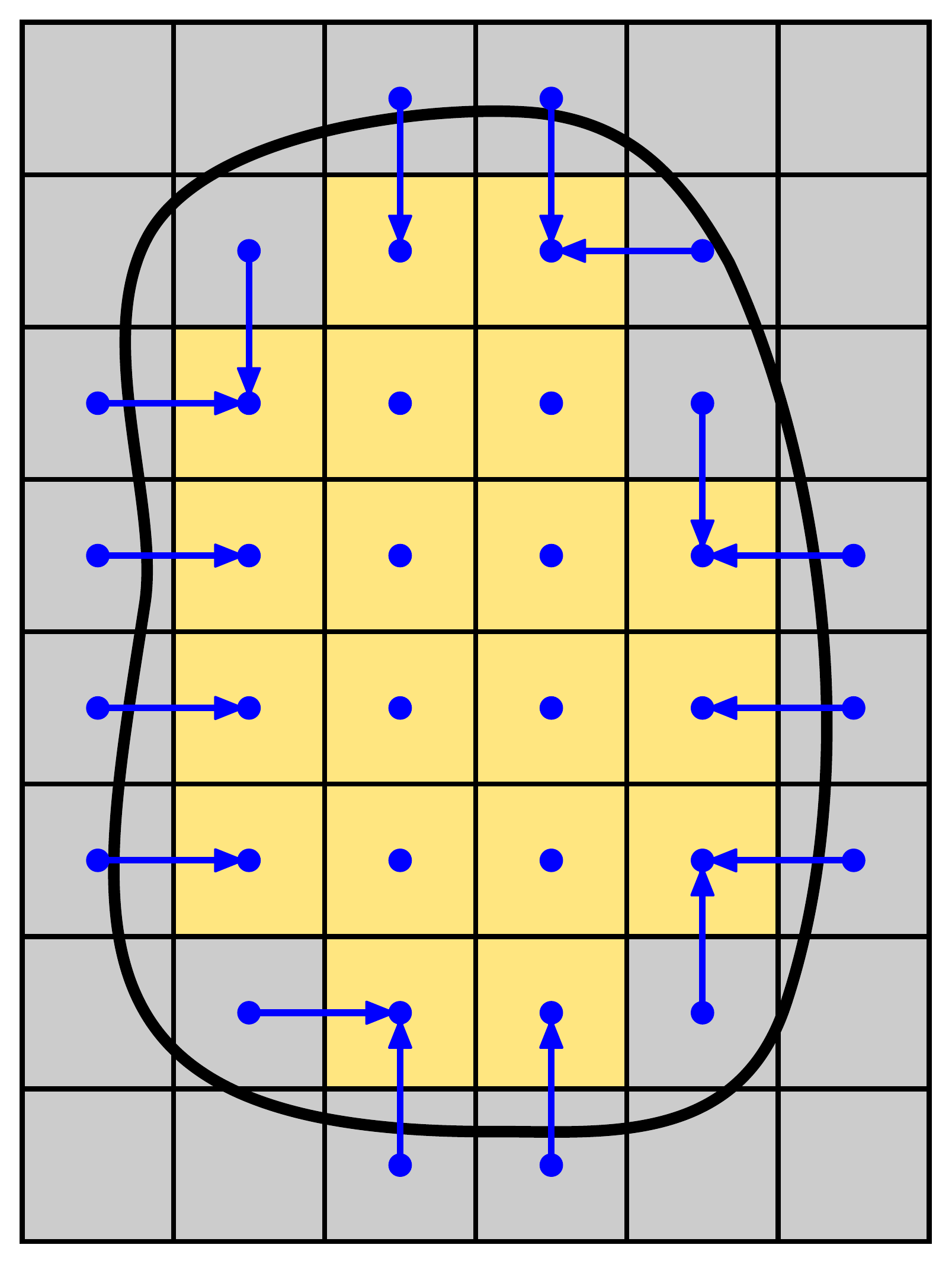}
    \caption{Step 2.}
    \label{fig:aggr-steps-b}
  \end{subfigure}
  \begin{subfigure}{0.24\textwidth}
    \centering
    \includegraphics[width=0.9\textwidth]{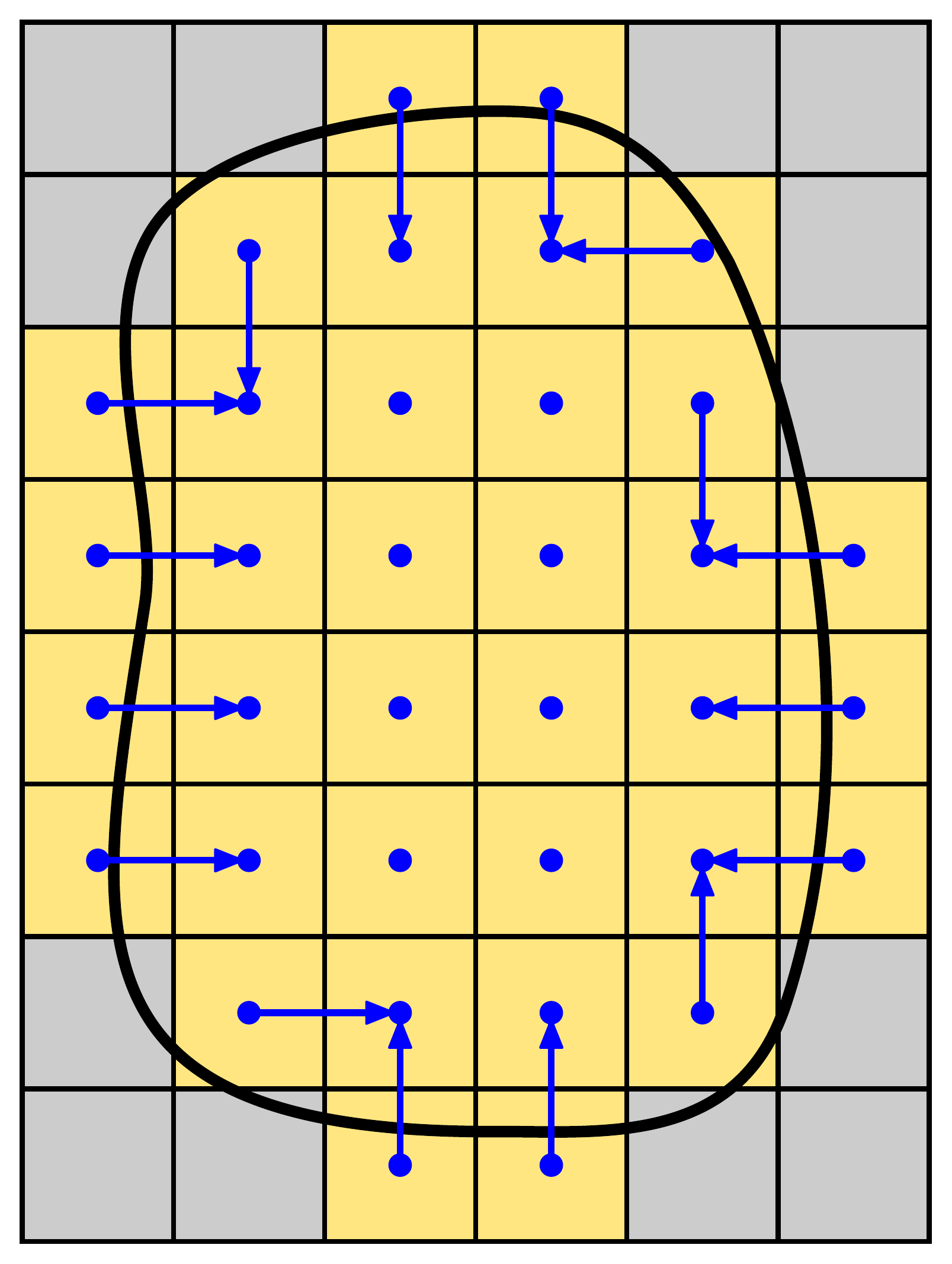}
    \caption{Step 3.}
    \label{fig:aggr-steps-c}
  \end{subfigure}
  \begin{subfigure}{0.24\textwidth}
    \centering
    \includegraphics[width=0.9\textwidth]{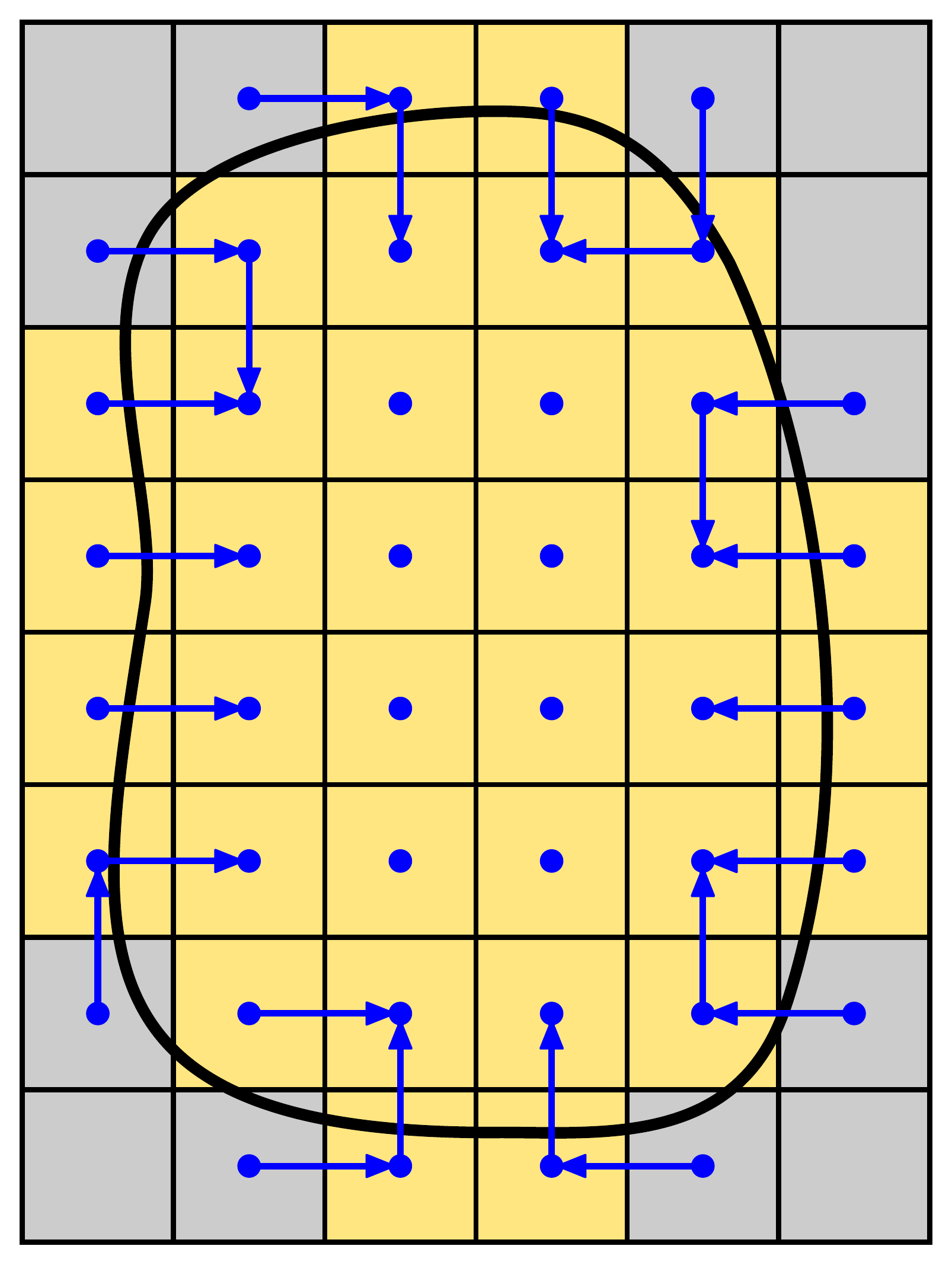}
    \caption{Step 2' (end).}
    \label{fig:aggr-steps-d}
  \end{subfigure}
  \caption{Illustration of the cell aggregation scheme defined in Algorithm \ref{alg:agg_sch}. A blue circle indicates cells for which the map $R$ is already defined with a valid value. Blue arrows indicate the face neighbor that is selected in line~\ref{ln:agg_sch_5} of the algorithm.}
  \label{fig:aggr-steps}
\end{figure}

In serial runs, a map $R$ fulfilling these properties is built with the strategy described in Alg.~\ref{alg:agg_sch}, and illustrated in Fig.~\ref{fig:aggr-steps} for a simple 2D case. The notation $\mathcal{K}^{\rm nei}(k)$ in line \ref{ln:agg_sch_2} of the algorithm denotes, for a given cell $T^k$ with id $k\in\mathcal{K}$, the set containing the ids  of the cells that share a face with cell $T^k$. In the following, we often denote the cells with ids in $\mathcal{K}^{\rm nei}(k)$ as the \emph{face neighbors} of cell $T^k$. Note that the algorithm consists in 3 basic steps (see Fig. ~\ref{fig:aggr-steps}). (a) We initialize the value of $R$ for interior cells by setting each interior cell as its own root cell, and we mark all interior cells as "touched". (b) We perform a loop in cut cells. For a given cut cell $T^k$, we define its root cell taking the root cell of one of the face neighbors of $T^k$ that have been already marked as "touched" (if such a {\rev neighbor exists}). (c) After the loop, we mark all cells that have been assigned with a root cell id as "touched". Then, steps (b) and (c) are repeated as many times as needed until all the cells are assigned with a root cell id.

In line~\ref{ln:agg_sch_5} of the algorithm, it is needed to choose one among the (possibly) several face neighbors that  have already been "touched". Several criteria can be considered for this purpose. In our implementation, we chose the neighbor having the closest root cell (in terms of the euclidean distance between cell barycenters). In case of a tie between two or more neighbors, we chose one arbitrarily (e.g., the neighbor with smaller cell id). During the algorithm, we also populate a map $N:\mathcal{K}\rightarrow\mathcal{K}$ that stores the ids of the chosen face neighbor. More precisely, for a given cut cell $T^k$  with id $k\in\mathcal{K}^{\rm cut}$,  $N(k)$ returns the id of face neighbor selected in line~\ref{ln:agg_sch_5}, whereas, for an interior cell $T^k$ with $k\in\mathcal{K}^{\rm in}$, $N(k)$ returns its own cell id, namely $N(k) = k$. {\rev The map $N$ plays an important role in the distributed-memory implementation}. Note that the aggregation scheme can be easily applied to arbitrary spatial dimensions. It only requires to find face neighbors and a classification of cells as internal, external, or cut. These procedures are widely available in embedded domain \ac{fe} codes.

\begin{algorithm}[ht!]
\caption{Serial cell aggregation algorithm}
\label{alg:agg_sch}
\SetKwIF{Unless}{ElseIf}{Else}{unless}{do}{else if}{else}{end}
\vspace*{0.5em}

Initialize $R(k)\leftarrow -1$ and $N(k)=-1$ for all $k\in\mathcal{K}$\nllabel{ln:agg_sch_0}\;
Set $R(k) \leftarrow k$ and $N(k)\leftarrow k$ for all $k\in\mathcal{K}^{\rm in}$\nllabel{ln:agg_sch_00}\;
Set as touched all cell ids $k\in\mathcal{K}$ such that $R(k)\neq -1$\nllabel{ln:agg_sch_1}\;
\For(\nllabel{ln:agg_sch_3}){$k\in\mathcal{K}^{\rm cut}$ such that $k$ has not been touched yet}{

  $\mathcal{K}^{\rm can}\leftarrow \emptyset$\;
  \For(\nllabel{ln:agg_sch_2}){ $k'\in\mathcal{K}^{\rm nei}(k)$ such that cell id $k'$ is already touched }{
    Let $F$  be the face shared by cells $T^k$ and $T^{k'}$\;
    \If{ $F\cap\Omega\neq\emptyset$ }{
    Add $k'$ to the set of candidates $\mathcal{K}^{\rm can}$\;
    }
  }

  \If{ $\mathcal{K}^{\rm can} \neq\emptyset$}{
   Choose an arbitrary cell id $k'\in\mathcal{K}^{\rm can}$ (e.g., the one with the closest root cell) \nllabel{ln:agg_sch_5}\;
       $R(k)\leftarrow R(k')$\;
       $N(k)\leftarrow k'$\nllabel{ln:agg_sch_4}\;
   }

}
\If{$R(k)= -1 $ for some cell id $k\in\mathcal{K}$}{
Go to line \ref{ln:agg_sch_1}\;
}

\end{algorithm}

\subsection{Standard and aggregated finite element spaces}\label{sec:fe-spaces}

Using the root cell map $R$, the authors in \cite{Badia2018} define \ac{agfe} spaces, which are the core of the \ac{agg} method. The aim of this section is to present such  \ac{agfe} spaces in a serial context, for its later extension to the distributed-memory case. To this end, we first define the usual conforming \ac{fe} space associated with the active part of the mesh, namely
\begin{equation}
\mathcal{V}^{\rm std} \doteq \{ v \in {C^0}(\Omega^{\rm act}) \, : \, v|_T \in
\mathcal{V}(T) \text{ for any } T \in \mathcal{T} \}.
\end{equation}
We denote by $\mathcal{V}(T)$ a vector space of functions defined on a generic cell $T\in\mathcal{T}$. {\rev In the numerical examples, we consider hexahedral cells and tri-linear shape functions. That is,  we define $\mathcal{V}(T) \doteq \mathcal{Q}_q(T)$ as  the space of polynomials that are of degree less or equal to $q$ ($q=1$ in the examples) with respect to each variable in $(x_1, \ldots, x_d)\in T$. Even though $d$-dimensional simplices are not considered in the numerical experiments, they can be also used in the presented method. In this case, we would define   $\mathcal{V}(T)\doteq \mathcal{P}_q(T)$ as the space of polynomials of order less or equal to $q$ in the variables $(x_1,\ldots,x_d)\in T$.}  In order to simplify notation, we define the elemental functional spaces $\mathcal{V}(T)$ in the physical cell $T\subset\Omega^{\rm act}$ (even though in the code we use reference parametric spaces as usual).
The space $\mathcal{V}^{\rm std}$ is the \emph{standard} functional space used in unfitted \ac{fe} methods (see, e.g., \cite{badia_robust_2017,DePrenter2017,Schillinger2015}).
 It is well known that it leads to arbitrary ill conditioned systems of linear algebraic equations (if no extra technique is used to remedy it).  This is the main motivation in order to introduce the \ac{agfe} spaces (see, e.g., \cite{Badia2018,Badia2018a}), since they lead to linear systems whose condition number is not affected by small cuts.

 In order to introduce the \ac{agfe} spaces we need some extra notation.
 For the sake of simplicity in the exposition and without loss of generality, we restrict ourselves to scalar-valued Lagrangian \acp{fe}. The extension to vector-valued or tensor-valued Lagrangian elements is straightforward by applying the same approach component by component. The extension to other types of \acp{fe} (e.g., Nédélec elements \cite{Olm2018}) is also possible. It only requires to work with the specific \acp{dof} of the element instead of the Lagrangian nodal values.

  Let $\{x^{a,k}\}_{a\in\Lambda(T^k)}$ be the set containing the (physical) coordinates of the Lagrangian nodes associated with cell $T^k\in\mathcal{T}$. In particular, $x^{a,k}\in\closure[1]{\Omega^{\rm act}}$ denotes the coordinate vector of node with id $a$ (in the local scope of cell $T^k$), and $\Lambda(T^k)\doteq\{1,\ldots,|\mathcal{V}(T^k)|\}$ is the set containing all local node ids in $T^k$.  Since (for scalar-valued Lagrangian spaces) each local node is associated with one \ac{dof} of the space $\mathcal{V}(T^k)$, we also refer to the local indices $a\in\Lambda(T^k)$ as the local \ac{dof} ids of cell $T^k$ and $\Lambda(T^k)$ as the set of local \ac{dof} ids in $T^k$. We also need to introduce the set $\{\phi^{a,k}\}_{a\in\Lambda(T^k)}$ of Lagrangian shape functions in the (physical) cell $T^k$. Following an analogous notation as for the cell nodal coordinates, the symbol $\phi^{a,k}$ refers to the shape function associated with the local \ac{dof} id $a\in\Lambda(T^k)$ in cell $T^k$. Clearly, $\phi^{a,k}(x^{b,k})=\delta_{ab}$ for any $a,b\in\Lambda(T^k)$, and $v^k=\sum_{a\in\Lambda(T^k)} v^{a,k} \phi^{a,k}$, with $v^{a,k}\doteq v^{k}(x^{a,k})$ being the cell nodal / \ac{dof} values for a given function $v^k\in\mathcal{V}(T^k)$.


On the other hand, we denote by  $\mathcal{X}\doteq\{x^i\}_{i\in\mathcal{I}}$ and  $\{\phi^i\}_{i\in\mathcal{I}}$  the set of \emph{global} nodal coordinates and the set of \emph{global} shape functions of space $\mathcal{V}^{\rm std}$, resp., where $x^i$ and $\phi^i$  are the coordinate vector and shape function of the mesh node with global id $i$ and $\mathcal{I}\doteq\{1,\ldots,|\mathcal{V}^{\rm std}|\}$ is the set of global node ids / \ac{dof} ids associated with the standard \ac{fe} space $\mathcal{V}^{\rm std}$. Clearly $\phi^i(x^j)=\delta_{ij}$, for any $i,j\in\mathcal{I}$, and $v=\sum_{i\in\mathcal{I}}v^i\phi^i$, with $v^i\doteq v(x^i)$ being the global nodal / \ac{dof} values for given function $v\in\mathcal{V}^{\rm std}$ .

We also need to introduce, for each cell $T^k\in\mathcal{T}$, the usual transformation between the cell local \ac{dof} ids and the corresponding global \ac{dof} ids. To this end, we denote by $g^{a,k}\in\mathcal{I}$ the global \ac{dof} id associated with the local \ac{dof} id $a$ in cell $T^k$.
Using the local-to-global map described by the values $g^{a,k}$, we define the set of global \ac{dof} ids touched by a cell $T^k$, i.e., $\mathcal{G}(k)\doteq\{g^{a,k}\}_{a\in\Lambda(T^k)}$. In addition, we introduce the notation $[\phi^i]_k$ representing  the local shape function in cell $T^k$ associated with the global one $\phi^i$, namely $[\phi^i]_k\doteq\phi^{a,k}$ with $i=g^{a,k}$. Note that $[\phi^i]_k$  and $\phi^i$ coincide inside the cell $T^k$, but they are different elsewhere. 

\begin{figure}[ht!]
\centering

  \begin{subfigure}{0.2\textwidth}
    \begin{tabular}{l}
      {\color{blue} $\bullet$} \small Free \acp{dof}
      \\
      {\color{red} $\times$} \small Constrained \acp{dof}
    \end{tabular}
  \end{subfigure}

\begin{subfigure}[b]{0.25\textwidth}
\includegraphics[width=0.95\textwidth]{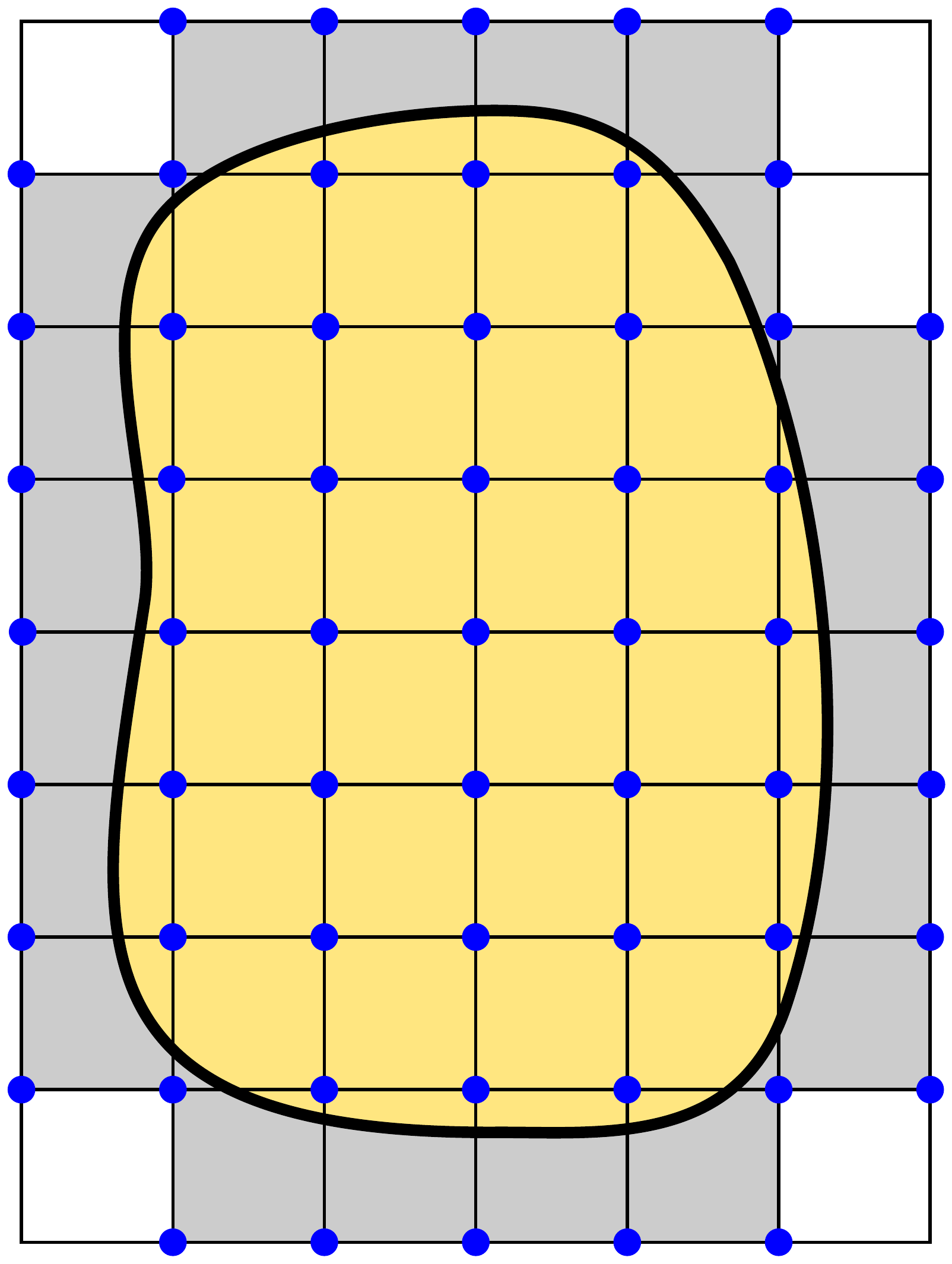}
\caption{$\mathcal{V}^{\rm std}$}
\label{fig:spaces-a}
\end{subfigure}
\begin{subfigure}[b]{0.25\textwidth}
\includegraphics[width=0.95\textwidth]{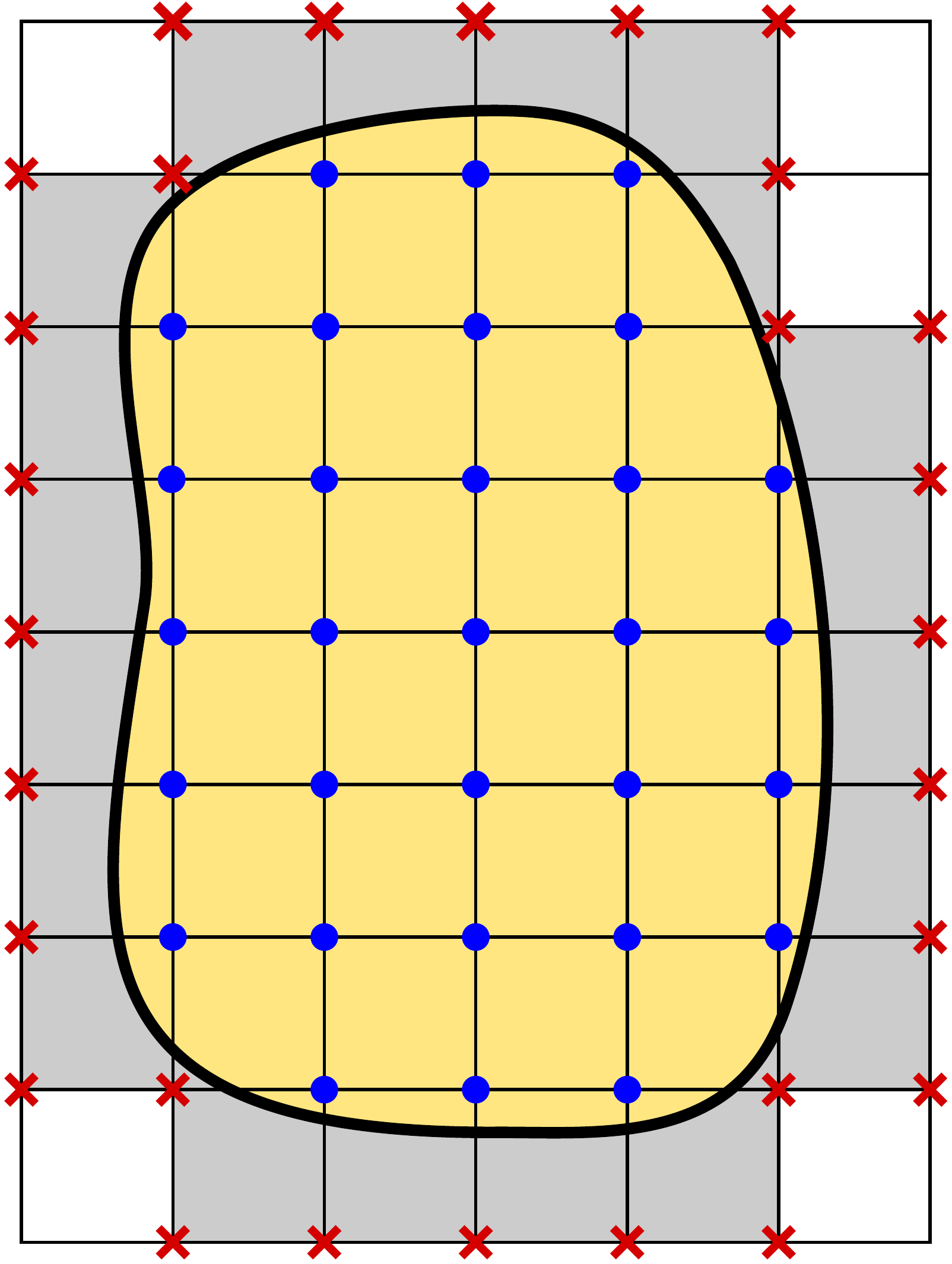}
\caption{$\mathcal{V}^{\rm agg}$}
\label{fig:spaces-b}
\end{subfigure}
\begin{subfigure}[b]{0.25\textwidth}
\includegraphics[width=0.943\textwidth]{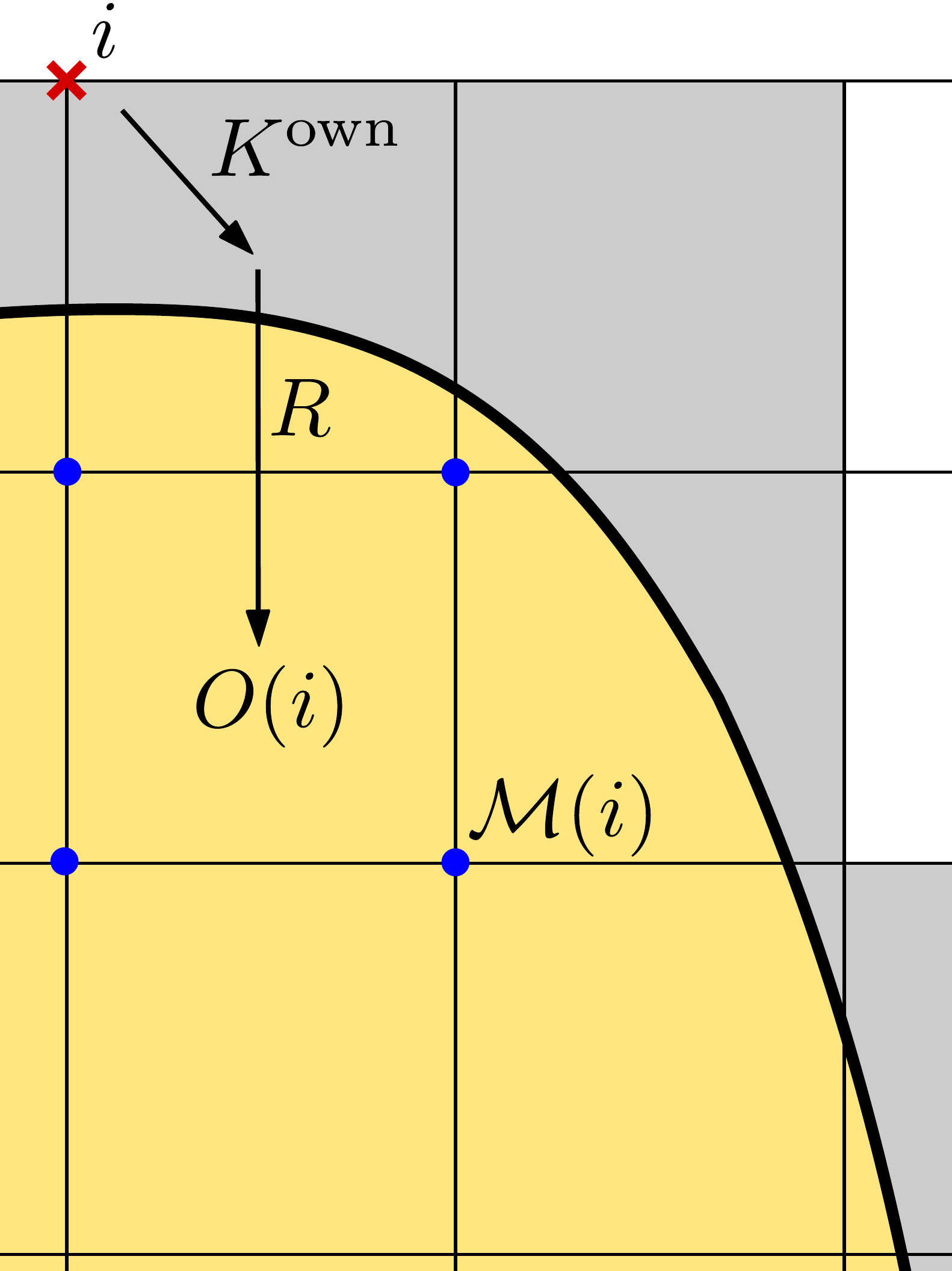}
\caption{$\mathcal{V}^{\rm agg}$ (detail)}
 \label{fig:spaces-c}
\end{subfigure}

 \caption{Illustration of the spaces $\mathcal{V}^{\rm std}$ and $\mathcal{V}^{\rm std}$ for first order Lagrangian elements.}
 \label{fig:spaces}
\end{figure}

Finally, we introduce the following classification of the global \ac{dof} ids into \emph{interior} and \emph{exterior} \ac{dof} ids. The set of interior \ac{dof} ids is defined as $\mathcal{I}^{\rm in}\doteq\{i\in\mathcal{I}:\ i\in\mathcal{G}(k)\text{ for some } k\in\mathcal{K}^{\rm in}\}$, i.e., {\rev all global \ac{dof} ids touched by the interior cells}. On the other hand, the set of exterior \ac{dof} ids is defined as $\mathcal{I}^{\rm out}=\mathcal{I}\setminus\mathcal{I}^{\rm in}$. Clearly, the ids in $\mathcal{I}^{\rm out}$ are only touched by cut cells. In Fig. \ref{fig:spaces-b}, interior and exterior \ac{dof} ids are represented with blue dots and red crosses, resp.
Let us also define a map $K^{\rm own}:\mathcal{I}^{\rm out}\rightarrow\mathcal{K}$ that takes a \ac{dof} id $i\in\mathcal{I}^{\rm out}$ and returns the cell id of one of the (possibly) several cells that contain the \ac{dof} id $i$ (see  Fig.~ \ref{fig:spaces-c}). In our implementation, we define $K^{\rm own}(i)$ as the smaller cell id $k\in\mathcal{K}$ such that $i\in\mathcal{G}(k)$.  By composing  the root cell map $R$ with $K^{\rm own}$, we obtain another map $O:\mathcal{I}^{\rm out}\rightarrow\mathcal{K}^{\rm in}$ that, for each exterior \ac{dof} id $i\in\mathcal{I}^{\rm out}$, returns an interior cell id. This latter map is defined as $O(i)\doteq R(K^{\rm own}(i))$. We denote cell $T^{O(i)}$ as the \emph{owner cell} of the \ac{dof} id $i$, or the \emph{root cell} associated with \ac{dof} id $i${\rev , } see  Fig.~ \ref{fig:spaces-c}.

With these notations, the \ac{agfe} space reads
\begin{equation}\label{eq:aggr-FES}
\mathcal{V}^{\rm agg} \doteq \{ v \in \mathcal{V}^{\rm std}  :   v^i = \sum_{j\in \mathcal{M}(i)} C_{ij}\ v^j \text{ for any } i\in \mathcal{I}^{\rm out}  \},
\end{equation}
where $C_{ij}\doteq[\phi^j]_{O(i)}(x^i)$ and $\mathcal{M}(i)\doteq\mathcal{G}(O(i))$.
Note that the space $\mathcal{V}^{\rm agg}$ is defined taking as starting point $\mathcal{V}^{\rm std}$, and adding some judiciously defined constraints to all exterior \acp{dof} (see Fig.~\ref{fig:spaces-a} and~\ref{fig:spaces-b}).  Only interior \acp{dof} are free in the space $\mathcal{V}^{\rm agg}$. Thus,   the dimension of $\mathcal{V}^{\rm agg}$ is equal to $|\mathcal{I}^{\rm in}|$. Moreover, $\mathcal{V}^{\rm agg}\subseteq \mathcal{V}^{\rm std}$ by construction. Each exterior \ac{dof} $i\in\mathcal{I}^{\rm out}$ is constrained as a linear combination of the \ac{dof} ids in the set $\mathcal{M}(i)$. Consequently, $\mathcal{M}(i)$ denotes the \emph{master \ac{dof} ids} of a given exterior \ac{dof} id $i\in\mathcal{I}^{\rm out}$. Note that $\mathcal{M}(i)$ contains all global \ac{dof} ids in the owner cell of $i$. Thus, constraints in \eqref{eq:aggr-FES} state that the value of $v$ at an outer node with id $i\in\mathcal{I}^{\rm out}$ is computed as the extrapolation of the values of $v$ at the nodes of the owner cell of $i$ (e.g., in Fig.~\ref{fig:spaces-c}, the value at the red cross is computed as an extrapolation of the values at the blue dots).

\subsection{Shape functions}\label{sec:shape-funs}

Before starting the discussion on the distributed-memory implementation, we detail how to assemble the system of linear algebraic equations associated with the \ac{agfe} space $\mathcal{V}^{\rm agg}$ in the serial case. To this end, we first need to introduce a shape function basis of $\mathcal{V}^{\rm agg}$, which is the purpose of this section. Let us define a restriction operator that, given a function $v\in\mathcal{V}^{\rm std}$, returns a function $\tilde v\in\mathcal{V}^{\rm agg}$ defined as:
\begin{equation}
\label{eq:aggr-restr}
\tilde v\doteq \sum_{j\in\mathcal{I}^{\rm in}} v(x^j)\ \phi^j +  \sum_{i\in\mathcal{I}^{\rm out}}  \sum_{j\in\mathcal{M}(i)} C_{ij} \ v(x^j)\ \phi^i.
\end{equation}
Note that function $\tilde v$ is computed in terms of the interior \ac{dof} ids only and it fulfills the constraints in \eqref{eq:aggr-FES} by construction. Moreover, $v$ and $\tilde v$ are different only in cut cells. The restriction does not affect the value of $v$ in interior cells.
A shape function basis $\{\tilde \phi^g\}_{g\in\mathcal{I}^{\rm in}}$  of $\mathcal{V}^{\rm agg}$ is obtained by applying restriction \eqref{eq:aggr-restr} to the shape functions of $\mathcal{V}^{\rm std}$ associated with interior \acp{dof}. More precisely, the shape function of $\mathcal{V}^{\rm agg}$ associated with the \ac{dof} id $g\in\mathcal{I}^{\rm in}$ is defined as
\begin{equation}
\tilde\phi^g \doteq \sum_{j\in\mathcal{I}^{\rm in}} \phi^g(x^j)\ \phi^j +  \sum_{i\in\mathcal{I}^{\rm out}}  \sum_{j\in\mathcal{M}(i)} C_{ij} \ \phi^g(x^j)\ \phi^i.
\end{equation}
Since the global shape functions of $\mathcal{V}^{\rm std}$ are such that $\phi^i(x^j)=\delta_{i j}$ for any $i, j\in\mathcal{I}$, previous formula simplifies to
\begin{equation}
\label{eq:aggr-sfuns}
\tilde\phi^g\doteq \phi^g + \sum_{i\in\mathcal{Z}(g)}C_{ig}\  \phi^i,
\end{equation}
where we have introduced the set $\mathcal{Z}(g)\doteq\{i\in\mathcal{I}^{\rm out}:\ g\in\mathcal{M}(i)\}$ referred to as the \emph{ slave \ac{dof} ids} of an interior \ac{dof} id $g\in\mathcal{I}^{\rm in}$.  Note that the functions $\tilde\phi^i$ satisfy $\tilde\phi^i(x^j)=\delta_{ij}$ for any $i,j\in\mathcal{I}^{\rm in}$. Thus, $\{\tilde\phi^i\}_{i\in\mathcal{I}^{\rm in}}$ is the shape function basis of $\mathcal{V}^{\rm agg}$ associated with the set of nodes $\{x^j\}_{j\in\mathcal{I}^{\rm in}}$. Note that the shape functions of the \ac{agfe} are associated with the interior nodes of the mesh.

\subsection{Finite element assembly}\label{sec:serial-assembly}

We present now the assembly operations needed to build the matrix $\mathbf{A}$ and vector $\mathbf{b}$ representing the discrete counterparts of a given bi-linear form $\mathrm{a}(\cdot,\cdot)$ and a linear form $\mathrm{b}(\cdot)$ with respect to the \ac{agfe} space $\mathcal{V}^{\rm agg}$. {\rev The strategy we have followed is related with the way linear constraints are imposed in the assembly process in reference \cite{Shephard1984}.}

 Using the previously defined shape functions of $\mathcal{V}^{\rm agg}$, the entries of the matrix $\mathbf{A}$ and the vector $\mathbf{b}$ are defined as $\left[\mathbf{A}\right]_{ij}\doteq \mathrm{a}(\tilde\phi^i,\tilde\phi^j)$ and  $\left[ \mathbf{b} \right]_i\doteq \mathrm{b}(\tilde\phi^i)$ for each $i,j\in\mathcal{I}^{\rm in}$. As usual, the global matrix and vector are computed as the assembly of elemental matrices and vectors
$[\mathbf{A}_k]_{ab}\doteq\mathrm{a}_k(\phi^{a,k},\phi^{b,k})$ and $[\mathbf{b}_k]_{a}\doteq\mathrm{b}_k(\phi^{a,k})$,  $a,b\in\Lambda
(T^k)$, associated with a generic cell $T^k\in\mathcal{T}$. Here, $\mathrm{a}_k(\cdot,\cdot)$ (resp. $\mathrm{b}_k(\cdot)$) is the restriction of $\mathrm{a}(\cdot,\cdot)$ (resp. $\mathrm{b}(\cdot)$) to cell $T^k\in\mathcal{T}$.

For local \ac{dof} ids $a\in\Lambda(T^k)$ such that its corresponding global \ac{dof} id $g^{a,k}\in\mathcal{I}^{\rm in}$ is interior, the assembly is standard, i.e., the value $[\mathbf{b}_k]_{a}$ is assembled into the global entry $[\mathbf{b}]_{i}$ such that $i=g^{a,k}$.  However, for local ids such that $g^{a,k}\notin\mathcal{I}^{\rm in}$, the assembly is slightly different. To clarify this, we use the definition of $\tilde\phi^g$ (see Eq.~\eqref{eq:aggr-sfuns})  in order to get a more detailed expression for $\left[\mathbf{A}\right]_{ij}$ and $\left[ \mathbf{b} \right]_i$, namely
\begin{equation}
\label{eq:bi_detail}
\left[ \mathbf{b} \right]_i=\mathrm{b}(\phi^i)+\sum_{n\in\mathcal{Z}(i)}C_{ni}\mathrm{b}(\phi^n),\qquad \text{for }i\in\mathcal{I}^{\rm in}, \text{ and}
\end{equation}
\begin{equation}
\begin{aligned}
\left[\mathbf{A}\right]_{ij} =\mathrm{a}(\phi^i,\phi^j) + \sum_{n\in\mathcal{Z}(i)}C_{ni}\ \mathrm{a}(\phi^n,\phi^j) +  \sum_{m\in\mathcal{Z}(j)}C_{mj}\ \mathrm{a}(\phi^i,\phi^m)
+  \sum_{n\in\mathcal{Z}(i)} \sum_{m\in\mathcal{Z}(j)} C_{ni} C_{mj}\ \mathrm{a}(\phi^n,\phi^m),
\end{aligned}
\end{equation}
for $i,j\in\mathcal{I}^{\rm in}$. In addition, we rewrite formula  \eqref{eq:bi_detail} in terms of cell-local quantities:
\begin{equation}
\label{eq:bi_detail_k}
[\mathbf{b}]_i = \sum_{k\in\mathcal{B}(i)} [\mathbf{b}_k]_a + \sum_{n\in\mathcal{Z}(i)}  \sum_{k\in\mathcal{B}(n)} C_{ni}  [\mathbf{b}_k]_b\quad\text{ for }i\in\mathcal{I}^{\rm in},
\end{equation}
with $i=g^{a,k}$, $n=g^{b,k}$, and $\mathcal{B}(i)\doteq\{k\in\mathcal{K}: i=g^{c,k} \text{ for some } c\in\Lambda(T^k)\}$ being the set containing the ids of cells sharing \ac{dof} id $i$.
It is clear by observing formula \eqref{eq:bi_detail_k} that the assembly of the entry $[\mathbf{b}_k]_a$  associated with a local \ac{dof} id $a\in\Lambda(T^k)$ such that $g^{a,k}\in\mathcal{I}^{\rm in}$ is done as usual, i.e., $[\mathbf{b}_k]_a$ is added to the corresponding global entry $[\mathbf{b}]_i$ with $i=g^{a,k}$. However, the  assembly is different for a local \ac{dof} id $a\in\Lambda(T^k)$ such that $g^{a,k}\notin\mathcal{I}^{\rm in}$. By observing the last term in Eq. \eqref{eq:bi_detail_k},  the local entry $[\mathbf{b}_k]_b$ contributes to the global entry $\left[ \mathbf{b} \right]_i$ for  all global \ac{dof} ids $i\in\mathcal{I}^{\rm in}$ that satisfy $g^{b,k}\in\mathcal{Z}(i)$, which is equivalent to say that $i\in\mathcal{M}(g^{b,k})$. Thus, the value $[\mathbf{b}_k]_b$ scaled by a factor $C_{ni}$, with $n=g^{b,k}$, is added to all global entries  $\left[ \mathbf{b} \right]_i$ that satisfy $i\in\mathcal{M}(g^{b,k})$. That is, the local entries $[\mathbf{b}_k]_a$ associated with constrained \ac{dof} ids $g^{a,k}\notin\mathcal{I}^{\rm in}$ contribute to the global entries corresponding to the master \ac{dof} ids in $\mathcal{M}(g^{a,k})$. This is a convenient way of imposing constraints, since it does not lead to saddle point problems as it would be the case if we would use the method of Lagrange multipliers.

The same rationale applies for the entries of the system matrix. In this case, the assembly operation to be performed depends on the two local \ac{dof} ids associated with   $\left[\mathbf{A}_k\right]_{ab}$, leading to four different cases. A summary of all different assembly operations is given in Table~\ref{tab:assembly-cases} for all possible cases. The symbol += stands for an in-place sum, i.e., the value of the right is added into the value on the left. Note that the procedures needed to assemble in $\mathcal{V}^{\rm agg}$ are closely related to the ones required in octree-based adaptive mesh refinement codes with non-conforming meshes. The imposition of continuity constraints at so called \emph{hanging nodes} (see, e.g., \cite{Badia2019b}) leads to assembly operations equivalent to the ones presented in Table \ref{tab:assembly-cases}. Thus, \ac{agfe} spaces can be easily incorporated in existing octree-based adaptive codes.

\begin{table}[ht!]

\newlength{\lrow}
\setlength{\lrow}{0.6em}

\begin{tabular}{p{0.25\textwidth}p{0.7\textwidth}}
\toprule
Case & Assembly operation\\
\midrule
Right hand side vector\\
$g^{a,k}\in\mathcal{I}^{\rm in}$ & $[\mathbf{b}]_i \text{ += }   [\mathbf{b}_k]_a$, $i=g^{a,k}$ \quad (standard case)\\[\lrow]
$g^{a,k}\notin\mathcal{I}^{\rm in}$ & $[\mathbf{b}]_i \text{ += }   C_{ni}\ [\mathbf{b}_k]_a$,  $n=g^{a,k}$,  $\forall i\in\mathcal{M}(n)$\\[2\lrow]
System matrix\\
$g^{a,k}\in\mathcal{I}^{\rm in}$ and  $g^{b,k}\in\mathcal{I}^{\rm in}$ & $[\mathbf{A}]_{ij}\text{ += }  [\mathbf{A}_k]_{ab}$, $i=g^{a,k}$, $j=g^{b,k}$ \quad (standard case)\\[\lrow]
$g^{a,k}\notin\mathcal{I}^{\rm in}$ and  $g^{b,k}\in\mathcal{I}^{\rm in}$ & $[\mathbf{A}]_{ij}\text{ += }  C_{ni}\ [\mathbf{A}_k]_{ab}$, $n=g^{a,k}$, $j=g^{b,k}$, $\forall i\in\mathcal{M}(n)$\\[\lrow]
$g^{a,k}\in\mathcal{I}^{\rm in}$ and  $g^{b,k}\notin\mathcal{I}^{\rm in}$ & $[\mathbf{A}]_{ij}\text{ += }  C_{mj}\ [\mathbf{A}_k]_{ab}$, $i=g^{a,k}$   , $m=g^{b,k}$, $\forall j\in\mathcal{M}(m)$\\[\lrow]
$g^{a,k}\notin\mathcal{I}^{\rm in}$ and  $g^{b,k}\notin\mathcal{I}^{\rm in}$ & $[\mathbf{A}]_{ij}\text{ += }  C_{ni} C_{mj}\ [\mathbf{A}_k]_{ab}$, $n=g^{a,k}$   , $m=g^{b,k}$, $\forall i\in\mathcal{M}(n)$, $\forall j\in\mathcal{M}(m)$\\
\bottomrule
\end{tabular}
\vspace{0.5em}

\caption{Summary of the assembly operations related with the \ac{agfe} space $\mathcal{V}^{\rm agg}$. {\rev These assembly operations are related with the methodology introduced in \cite{Shephard1984}.}  }
\label{tab:assembly-cases}
\end{table}

\section{Distributed-memory implementation}\label{sec:dm-impl}

\subsection{Domain decomposition setup}\label{sec:dd-setup}

\begin{figure}[ht!]
  \centering
  \begin{subfigure}[b]{0.24\textwidth}
    \centering

     \begin{small}

      \hfill
      {\color[RGB]{255,230,128} \rule{10pt}{8pt}} $D_1$\hfill
      {\color[RGB]{204,204,204} \rule{10pt}{8pt}} $D_2$\hfill

     \end{small}

    \includegraphics[width=0.9\textwidth]{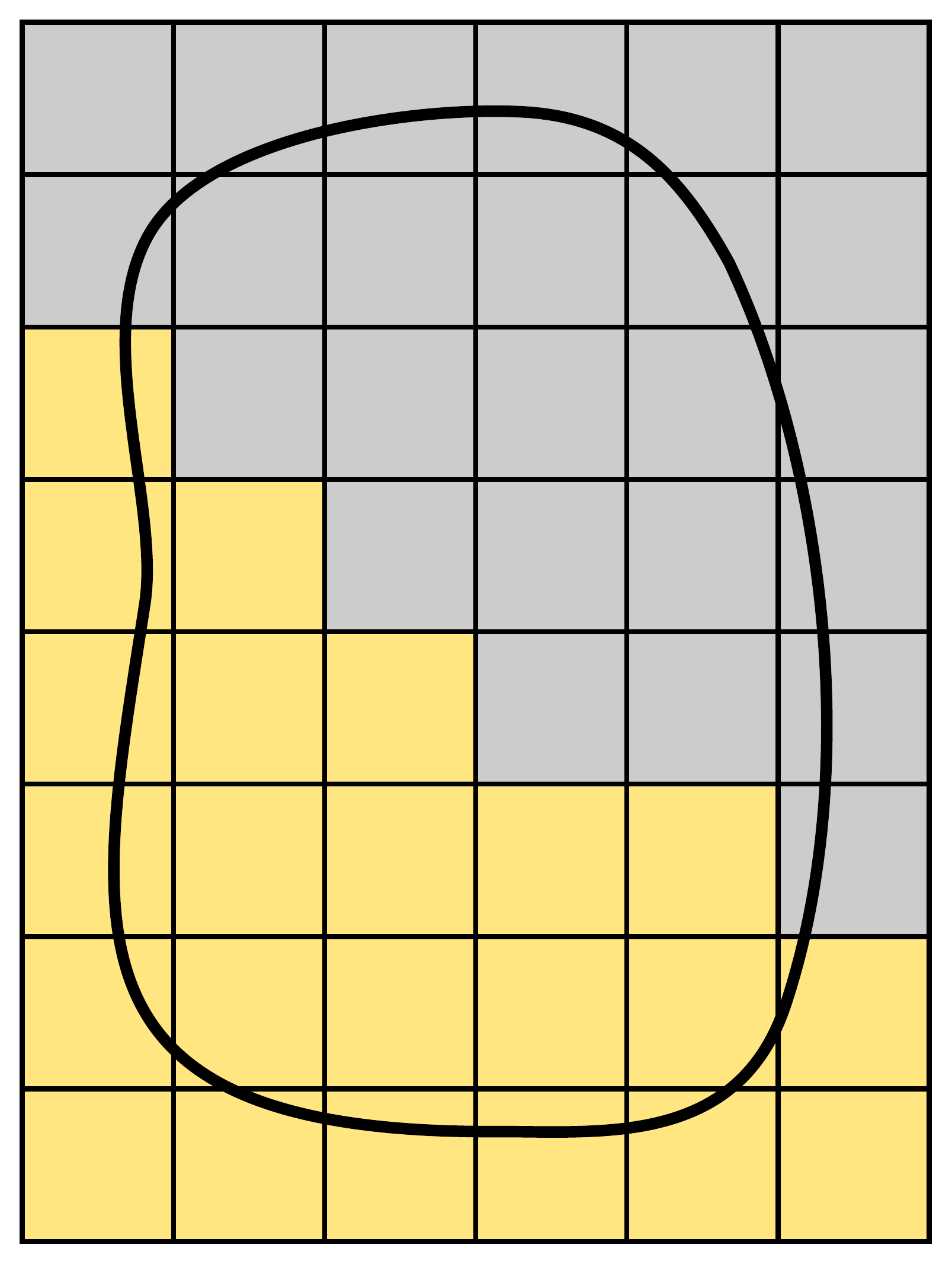}
    \caption{$\mathcal{D}=\{D_1,D_2\}$.}
    \label{fig:partition-setup-a}
  \end{subfigure}
  \begin{subfigure}[b]{0.24\textwidth}
    \centering
    \includegraphics[width=0.9\textwidth]{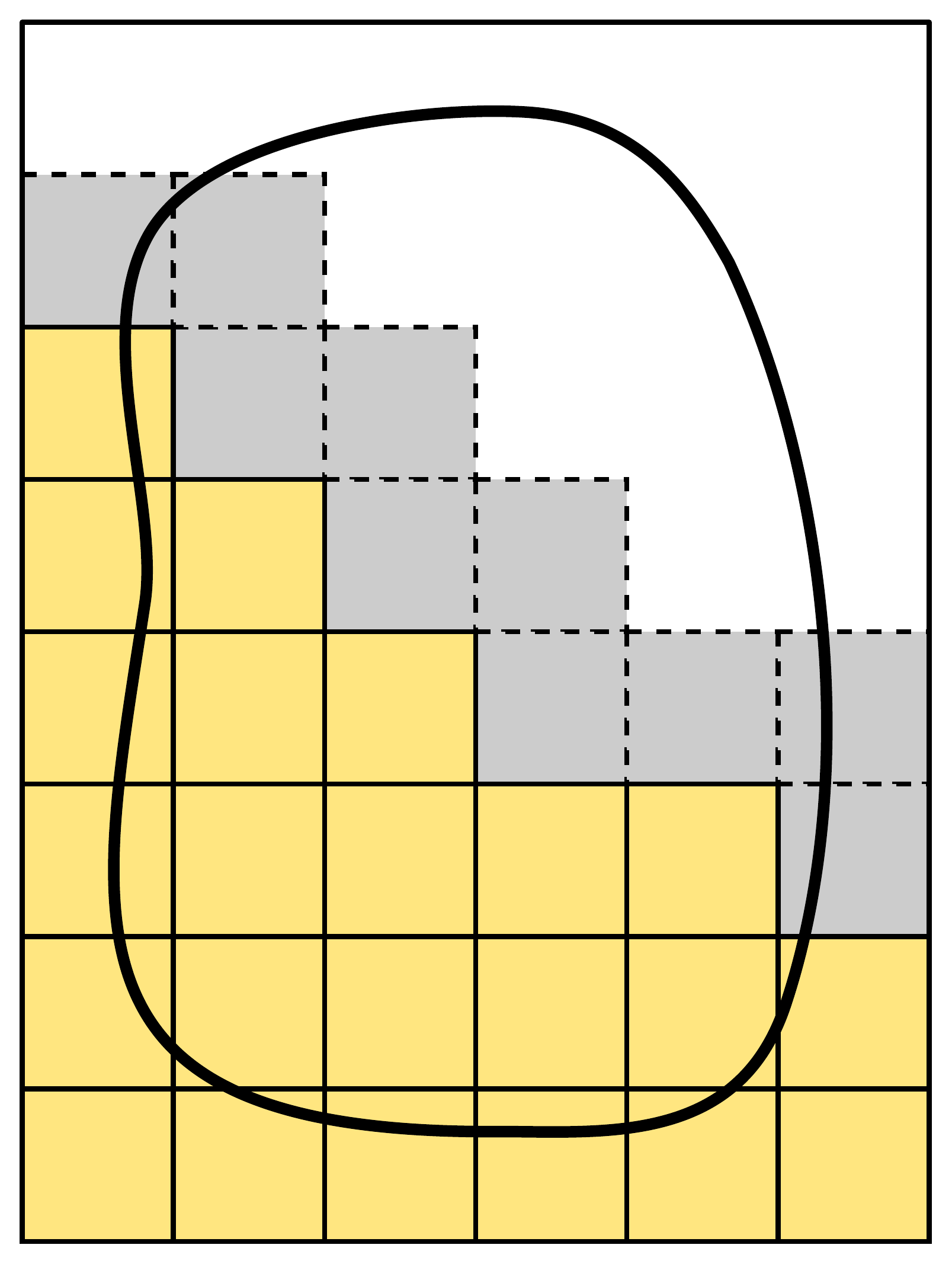}
    \caption{$\mathcal{T}^{\rm L}_1$ and $\mathcal{T}^{\rm G}_1$.}
    \label{fig:partition-setup-b}
  \end{subfigure}
    \begin{subfigure}[b]{0.24\textwidth}
    \centering
    \includegraphics[width=0.9\textwidth]{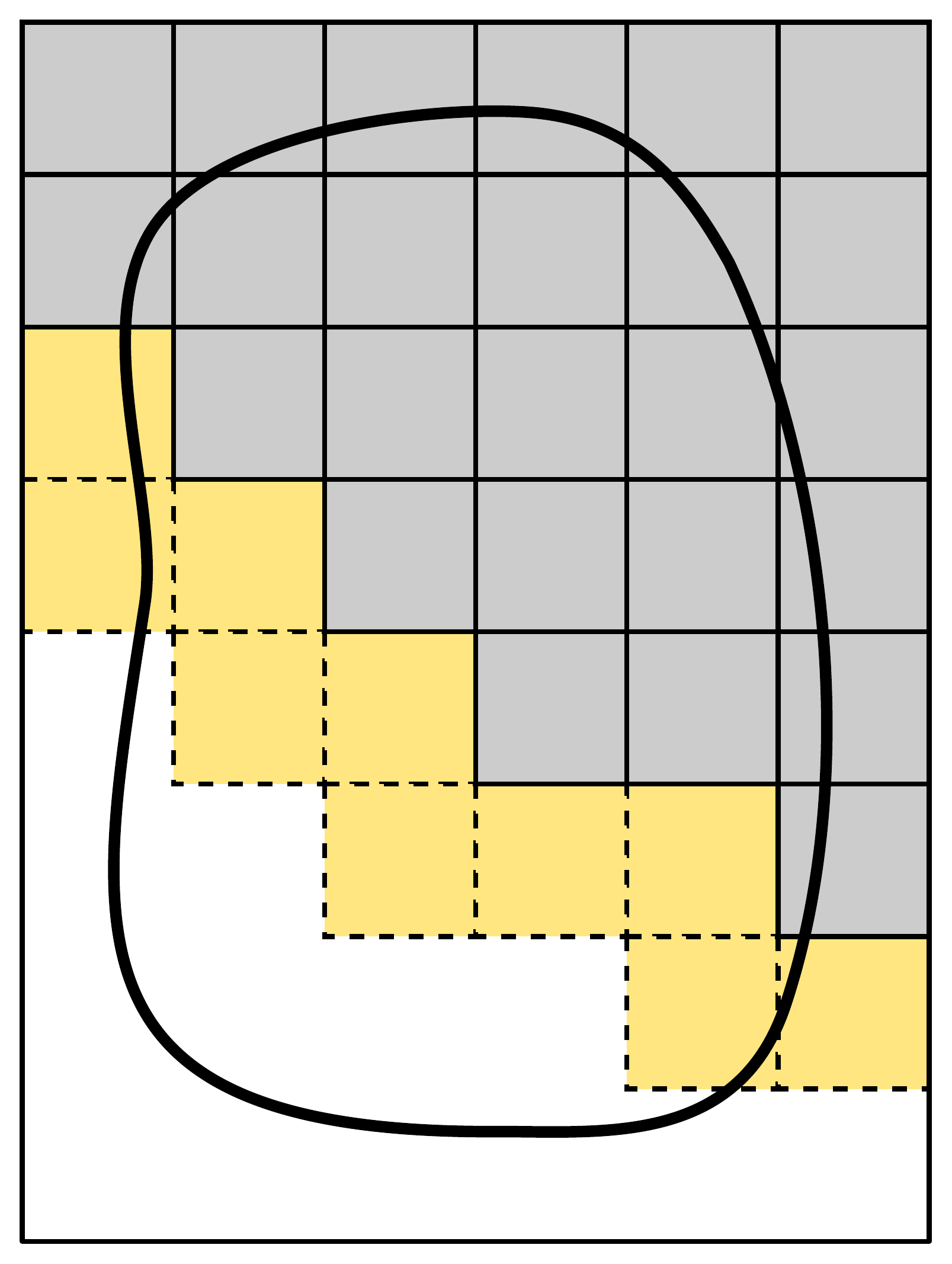}
    \caption{$\mathcal{T}^{\rm L}_2$ and $\mathcal{T}^{\rm G}_2$.}
    \label{fig:partition-setup-c}
  \end{subfigure}
  \caption{Domain decomposition setup illustrated for two sub-domains.
   ({\sc a}) Partition into sub-domains.  ({\sc b--c}) Local and ghost cells for each sub-domain. Local cells are represented with solid lines, whereas ghost cells are represented with dashed ones.}
  \label{fig:partition-setup}
\end{figure}

After the presentation of the \ac{agg} in a serial context, we discuss its implementation in a distributed-memory framework. We start by introducing some notation associated with the partition of the mesh into sub-domains, which is the basis of the distributed-memory implementation.  Let $\mathcal{D}$ be a partition of $\Omega^{\rm art}$ into sub-domains obtained by the union of cells in the background mesh $\mathcal{T}^{\rm art}$, i.e., for each cell $T\in\mathcal{T}^{\rm art}$, there is a sub-domain $D\in\mathcal{D}$ such that $T\subset D$ (see Fig.~\ref{fig:partition-setup-a}).
The generation of the sub-domain partition $\mathcal{D}$ is specially simple in the context of embedded boundary methods since Cartesian grids are allowed.
Here, we consider the load-balancing mechanism available in \p4est, which is based on space-filling curves \citep{Bader2012} and has shown excellent scalability up to hundreds of thousands of cores. {\rev Note that our approach does not require graph-based partitioning techniques such as \texttt{metis} or \texttt{parmetis} \cite{Padua2011}}. 

For convenience, we consider a numbering of the sub-domains in $\mathcal{D}$, where $D_s$ is the sub-domain with id $s$ and $\mathcal{S}\doteq\{1,\ldots,|\mathcal{D}|\}$ is the index set containing all sub-domain ids. With this notation the set of sub-domains $\mathcal{D}$ is expressed as $\mathcal{D}=\{D_s\}_{s\in\mathcal{S}}$.
For a given $s\in\mathcal{S}$, we define the set of \emph{local cells} $\mathcal{T}^{\rm L}_s$ as the cells included in $D_s$, i.e., $\mathcal{T}^{\rm L}_s\doteq\{T\in\mathcal{T}^{\rm art} :\ T\subset D_s\}$. On the other hand, we define the set of \emph{ghost cells} $\mathcal{T}^{\rm G}_s$  of $D_s$ as the off-subdomain cells that share at least one vertex, edge or face with $D_s$, namely $\mathcal{T}^{\rm G}_s\doteq\{T\in (\mathcal{T}^{\rm art} \setminus \mathcal{T}^{\rm L}_s  ):\ \closure[1]{T} \cap \closure[1]{D_s} \neq \emptyset\}$ (see Fig.~\ref{fig:partition-setup-b}, and \ref{fig:partition-setup-c}). The union of local and ghost cells, $\mathcal{T}_s\doteq\mathcal{T}^{\rm L}_s\cup\mathcal{T}^{\rm G}_s$, is called the set of \emph{locally relevant} cells in sub-domain $D_s$.

For the following developments, we need a \emph{local} numbering for the locally relevant cells $\mathcal{T}_s$. That is, we write $\mathcal{T}_s=\{T^\ell_s\}_{\ell\in\mathcal{L}_s}$, where $T^\ell_s$ denotes the cell with local id $\ell$ in the scope of sub-domain $D_s$ and $\mathcal{L}_s\doteq\{1,\ldots,|\mathcal{T}_s|\}$ is the set of all local cell ids associated with $D_s$.  {\rev In addition, we define $\mathcal{L}^{\rm L}_s\doteq\{\ell\in\mathcal{L}_s:\ T^\ell_s\in\mathcal{T}^{\rm L}_s\}$ and $\mathcal{L}^{\rm G}_s\doteq\{\ell\in\mathcal{L}_s:\ T^\ell_s\in\mathcal{T}^{\rm G}_s\}$, which contain the local cell ids of local and ghost cells resp. We also need to identify which is the sub-domain of each cell in terms of local cell ids. To this end, we denote as $S_s(\ell)$ the id of the sub-domain, where the cell with id $\ell\in\mathcal{L}_s$ is a local cell. Clearly, $S_s(\ell)=s$ for all local cell ids $\ell\in\mathcal{L}^{\rm L}_s$ and $S_s(\ell)\neq s$ for ghost cell  ids $\ell\in\mathcal{L}^{\rm G}_s$.}

 For simplicity and without any loss of generality, we will assume (as we have done in the serial case) that exterior cells have been removed since they do not play any role. In particular, it implies that $\mathcal{T}^{\rm L}_s$, $\mathcal{T}^{\rm G}_s$, and $\mathcal{T}_s$ contain only interior or cut cells. As already stated in Sect.~ \ref{sec:imb_setup}, each cell in the active mesh $\mathcal{T}$ receives a global id $k\in\mathcal{K}$, which is unique across all cells in  $\mathcal{T}$.  Since $\mathcal{T}_s$ is a subset of $\mathcal{T}$, all cells in $\mathcal{T}_s$ have a global cell id.
{\rev We denote by $K_s(\ell)$ the global cell id associated with a local cell id $\ell\in\mathcal{L}_s$. Conversely, we denote as $L_s(k)$ the local cell id associated with the global cell id $k\in\mathcal{K}_s$. Here, $\mathcal{K}_s\doteq\{k\in\mathcal{K}:\ T^k\in\mathcal{T}_s\}$  is the} set of all global ids assigned to cells in $\mathcal{T}_s$. With this notation, we can write $\mathcal{T}_s=\{T^k\}_{k\in\mathcal{K}_s}$. {\rev Finally, we need to introduce $\mathcal{S}^{\rm nei}_s\doteq\{s'\in\mathcal{S}\setminus\{s\}:\ \closure[1]{D_{s'}}\cap \closure[1]{D_s}\neq \emptyset\}$, i.e., the set containing the sub-domain ids of the \emph{nearest neighbors} of $D_s$.}

\subsection{Distributed finite element spaces}\label{sec:par-fe-spaces}

{\rev The goal of this section is to express the \ac{agfe} space  $\mathcal{V}^{\rm agg}$ in terms of sub-domain local quantities in the framework of the domain decomposition setup}. By doing so, most of the notation needed in the subsequent text will be introduced in a natural way. For the sake of simplicity, we assume from now on (if not otherwise indicated) that we are located in the local scope of a generic sub-domain $D_s$ with id $s\in\mathcal{S}$. Also for simplicity (and without any loss of generality), we assume that we deal with scalar-valued Lagrangian \ac{fe} interpolations.

 First, we define $\mathcal{V}^{\rm std}_s\doteq\mathcal{V}^{\rm std}|_{D_s}$ the restriction of the standard \ac{fe} space $\mathcal{V}^{\rm std}$ to the sub-domain $D_s$.  Given a \ac{fe} function $v\in\mathcal{V}^{\rm std}$, we also introduce its sub-domain restriction: $v_s \doteq v|_{D_s}$. Note that $v_s$ is characterized by a subset of the \acp{dof} of $v$, which are the ones associated with the nodes of the mesh $\mathcal{T}^{\rm L}_s$, i.e., the nodes in the set $\mathcal{X}_s\doteq\{x\in\mathcal{X}:\ x\in \closure[1]{T} \text{ for some } T\in \mathcal{T}^{\rm L}_s\}$. As we have done previously for cells, we introduce a sub-domain local numbering for  nodes in $\mathcal{X}_s$, i.e., $\mathcal{X}_s=\{x^j_s\}_{j\in\mathcal{J}_s}$, where $x^j_s$ denotes the nodal coordinate of the mesh node with  local id $j$ in sub-domain $D_s$ and $\mathcal{J}_s\doteq\{1,\ldots,|\mathcal{X}_s|\}$ is the set of all sub-domain local node ids (also referred as set of local \ac{dof} ids since there is a one-to-one correspondence between nodes and \acp{dof}). Thus, each node in $\mathcal{X}_s$ has two different ids, a local one in the local mesh $\mathcal{T}_s$ and a {\rev global one in the entire mesh} $\mathcal{T}$. Let $\mathcal{I}_s\doteq\{i\in\mathcal{I}:\ x^i\in\mathcal{X}_s\}$ be the set of the global node ids (or global \ac{dof} ids) in sub-domain $D_s$ and let $I_s:\mathcal{J}_s\rightarrow\mathcal{I}_s$ be the map that, for a given local node id $j$, returns the corresponding global id $i=I_s(j)$ such that $x^i=x^j_s$.
With these notations, we can write a generic function $v_s\in\mathcal{V}^{\rm std}_s$ as a linear combination of sub-domain local quantities, i.e., $v_s(x)=\sum_{j\in\mathcal{J}_s} v^j_s\ \phi^j_s(x)$, where $v^j_s\doteq v_s(x^j_s) =v^i=v(x^i)$ and $\phi^j_s\doteq\phi^i|_{D_s}$ with $i=I_s(j)$ for all $j\in\mathcal{J}_s$. Clearly, $v^j_s$ and $\phi^j_s$ are the \ac{dof} and the shape function associated with the sub-domain local id $j$.

On the other hand, we introduce $\mathcal{V}^{\rm agg}_s\doteq\mathcal{V}^{\rm agg}|_{D_s}$ the restriction of the \ac{agfe} space $\mathcal{V}^{\rm agg}$ to sub-domain $D_s$. In order to characterize $\mathcal{V}^{\rm agg}_s$ in terms of sub-domain local quantities, we need to introduce some further notation. First, let us define $\phi^a_{s,\ell}$ the shape function associated with a generic \ac{dof} id $a\in\Lambda(T^\ell_s)$ in cell $T^\ell_s$, namely $\phi^a_{s,\ell}\doteq\phi^{a,k}$ with $k=K_s(\ell)$. Analogously, we introduce $x^a_{s,\ell}\doteq x^{a,k}$ and  $g^a_{s,\ell} \doteq g^{a,k}$, with $k=K_s(\ell)$, which represent the nodal coordinates and the global \ac{dof} id associated with a generic \ac{dof} id $a\in\Lambda(T^\ell_s)$ in cell $T^\ell_s$. Note that the symbols $g^a_{s,\ell}$ and $g^{a,k}$  express the same quantities, but the former expresses them in terms of local cell ids $\ell\in\mathcal{L}_s$, and the latter in terms of global cell ids $k\in\mathcal{K}_s$. This is true also for  $x^a_{s,\ell}$ and $x^{a,k}$, and $\phi^a_{s,\ell}$ and $\phi^{a,k}$, resp. We also define $j^a_{s,\ell}$ as the sub-domain local \ac{dof} id in $\mathcal{J}_s$ associated with the cell local \ac{dof} id $a\in\Lambda(T^\ell_s)$ in cell $T^\ell_s$. Formally, $j^a_{s,\ell}$ is expressed as $j^a_{s,\ell}\doteq I^{-1}_s(g^a_{s,\ell})$ {in terms of the inverse application $I^{-1}_s$ (which is not needed to be generated in practice). The process that generates $g^a_{s,\ell}$ and $j^a_{s,\ell}$ for each $a$, $\ell$, and $s$ is standard in \ac{fe} packages.
The sub-domain local cell-wise \ac{dof} ids $j^a_{s,\ell}$ are simply generated by applying independently in each sub-domain the standard procedure for generating \ac{dof} ids in the serial case (see, e.g., \citep{badia_fempar:_2017} for details). On the other hand, the generation of the global \ac{dof} ids $g^a_{s,\ell}$ is also done by means of standard tools (see, e.g., \cite{bangerth_algorithms_2012,Badia2019b}) as it is detailed later at the end of Sect.~\ref{sec:data-layout}.


We also need to introduce a classification of the locally relevant cells into interior, exterior, and cut. To this end, we define the index sets $\mathcal{L}_{s}^{\rm in}\doteq\{\ell\in\mathcal{L}_s:\ T^\ell_s\in\mathcal{T}^{\rm in}\}$ (resp. for $\mathcal{L}^{\rm out}_s$ and $\mathcal{L}^{\rm cut}_s$) containing the subset of $\mathcal{L}_s$ corresponding to interior cells (resp. exterior and cut).  
In our case, the set $\mathcal{L}^{\rm out}_s$ is empty since we have assumed that exterior cells have been removed from the mesh.
The classification of cells as interior or cut induces the following partition of sub-domain local \ac{dof} ids in $\mathcal{J}_s$.  Let  $\mathcal{J}^{\rm in}_s\doteq\{j\in\mathcal{J}_s:\ j=j^a_{s,\ell} \text{ for some } a\in\Lambda(T^\ell_s) \text{ and } \ell\in\mathcal{L}^{\rm in}_s \}$ be the set of \emph{interior} sub-domain local \ac{dof} ids, and let $\mathcal{J}^{\rm out}_s\doteq\mathcal{J}_s\setminus\mathcal{J}^{\rm in}_s$ be the set of \emph{exterior} sub-domain local \ac{dof} ids. Note that the \ac{dof} ids in $\mathcal{J}^{\rm in}_s$ are the ones touched by interior cells with ids in $\mathcal{L}^{\rm in}_s$, whereas the ids in $\mathcal{J}^{\rm out}_s$ are only touched by cut cells.

The last ingredient to characterize $\mathcal{V}^{\rm agg}_s$ in terms of sub-domain quantities is a distributed version of the root {\rev cell map $R$ and the map $K^\mathrm{\rm own}$}. On the one hand, the distributed root cell map  $R_s:\mathcal{L}_s\rightarrow\mathcal{K}^{\rm in}$ is a transformation that takes a local cell id $\ell\in\mathcal{L}_s$ and returns the global cell id $R_s(\ell)\in\mathcal{K}^{\rm in}$ of the root cell of $T^\ell_s$, namely $R_s(\ell)\doteq R(K_s(\ell))$. We encode the input value of $R_s(\cdot)$ with a local cell id in $\mathcal{L}_s$ since we only need to define $R_s$ for locally relevant cells in $\mathcal{T}_s$. However, the returned value of $R_s$ has to be encoded with a global cell id in $\mathcal{K}^{\rm in}$, and not with a local cell id in $\mathcal{L}_s$, since the root cell of $T^\ell_s$ \emph{can be outside the set of locally relevant cells} $\mathcal{T}_s$ (which are the only ones that have a local id in sub-domain $D_s$). Note also that, even though we have formally defined $R_s$ using $R$, the map $R$ is never built in practice. Constructing and storing $R$ so that each processor can retrieve the value $R(k)$ for any global cell id $k\in\mathcal{K}$ in the entire mesh would kill parallel (memory) scalability. Instead, we only build the distributed map $R_s$ as shown later in Sect.~\ref{sec:par-cell-agg}. In addition to $R_s$, we also construct another map $S^{\rm root}_s:\mathcal{L}_s\rightarrow\mathcal{S}$ that for each local cell id $\ell\in\mathcal{L}_s$ returns the id $S^{\rm root}_s(\ell)\in\mathcal{S}$ of the sub-domain in which the corresponding root cell is a local cell, namely  $s'=S^{\rm root}_s(\ell)$ is such that  $T^{\ell'}_{s'}\in\mathcal{T}^{\rm L}_{s'}$  with $\ell'=L_{s'}(R_s(\ell))$. Finally, we need to build a transformation $L^{\rm own}_s:\mathcal{J}^{\rm out}_s\rightarrow\mathcal{L}_s$, which is a distributed version of $K^\mathrm{\rm own}$.  For a local \ac{dof} id $j\in\mathcal{J}^{\rm out}_s$, the value $L^{\rm own}_s(j)$ is the local cell id $\ell\in\mathcal{L}_s$ of one of the (possibly several) cells  $T^\ell_s$ such that $j=j^a_{s,\ell}$ for some $a\in\Lambda(T^\ell_s)$. That is, $L^{\rm own}_s(j)$ is the local id of one of the (possibly) several cells that contain the sub-domain \ac{dof} id $j$. When more than one cell can be selected, we {\rev choose} one arbitrarily, e.g., the cell that has smaller global id $k\in\mathcal{K}_s$. Note that one can choose one cell arbitrarily, but the chosen criterion has to be consistent across all sub-domains. E.g., choosing the cell with smaller global cell id is consistent, but choosing the cell with smaller local id is not.

{\rev With} these notations, we can finally write $\mathcal{V}^{\rm agg}_s$ in terms of sub-domain quantities:
\begin{equation}\label{eq:agg_space_s}
\mathcal{V}^{\rm agg}_s=\left\lbrace v_s\in\mathcal{V}^{\rm std}_s :\ v^j_s = \sum_{a\in \Lambda(T^{\ell'}_{s'})}  C^{ja}_s\ v^a_{s',\ell'} \text{ for all }j\in\mathcal{J}^{\rm out}_s\right\rbrace,
\end{equation}
with $C^{ja}_s\doteq\phi^{a}_{s',\ell'}(x^j_s)$, $\ell' = L_{s'}(R_s(L^{\rm own}_s(j)))$  and $s'= S^{\rm root}_s(L^{\rm own}_s(j))$. By observing Eq.~\eqref{eq:agg_space_s}, the main difficulty of the implementation of the method becomes apparent. {\rev Note} that in the scope of the current sub-domain $D_s$, not all quantities are available in order to define the constraints in Eq.~\eqref{eq:agg_space_s}. In order to evaluate the coefficient $C^{ja}_s$, one needs access to the shape function in physical space $\phi^a_{s',\ell'}$,  which is one of the local shape basis functions of the cell with local id $\ell' = L_{s'}(R_s(L^{\rm own}_s(j)))$ in sub-domain $s'= S^{\rm root}_s(L^{\rm own}_s(j))$, which can be interpreted as the root cell associated with the \ac{dof} id $j$. Building up  $\phi^a_{s',\ell'}$ requires the physical coordinates of the nodes of the root cell, namely $x^b_{s',\ell'}$ for $b\in\Lambda
(T^{\ell'}_{s'})$. Since, in general, the sub-domain where the root cell is located is different from the current sub-domain, i.e., $s'\neq s$, the quantities $x^b_{s',\ell'}$ are not available in the scope of $D_s$ (see Fig.~\ref{fig:remote-root}). Moreover, one cannot directly use the communication pattern used to exchange cell-based data between nearest neighbor sub-domains
 because id $s'$ is not even necessarily the id of a nearest neighbor sub-domain, namely $s'\notin\mathcal{S}^{\rm nei}_s$ in general. Thus, an ad-hoc communication strategy has to be built in order to import quantities like $x^b_{s',\ell'}$ into the scope of {\rev the current processor}.


\begin{figure}[ht!]
\centering
\begin{subfigure}{0.5\textwidth}
\includegraphics[width=0.9\textwidth]{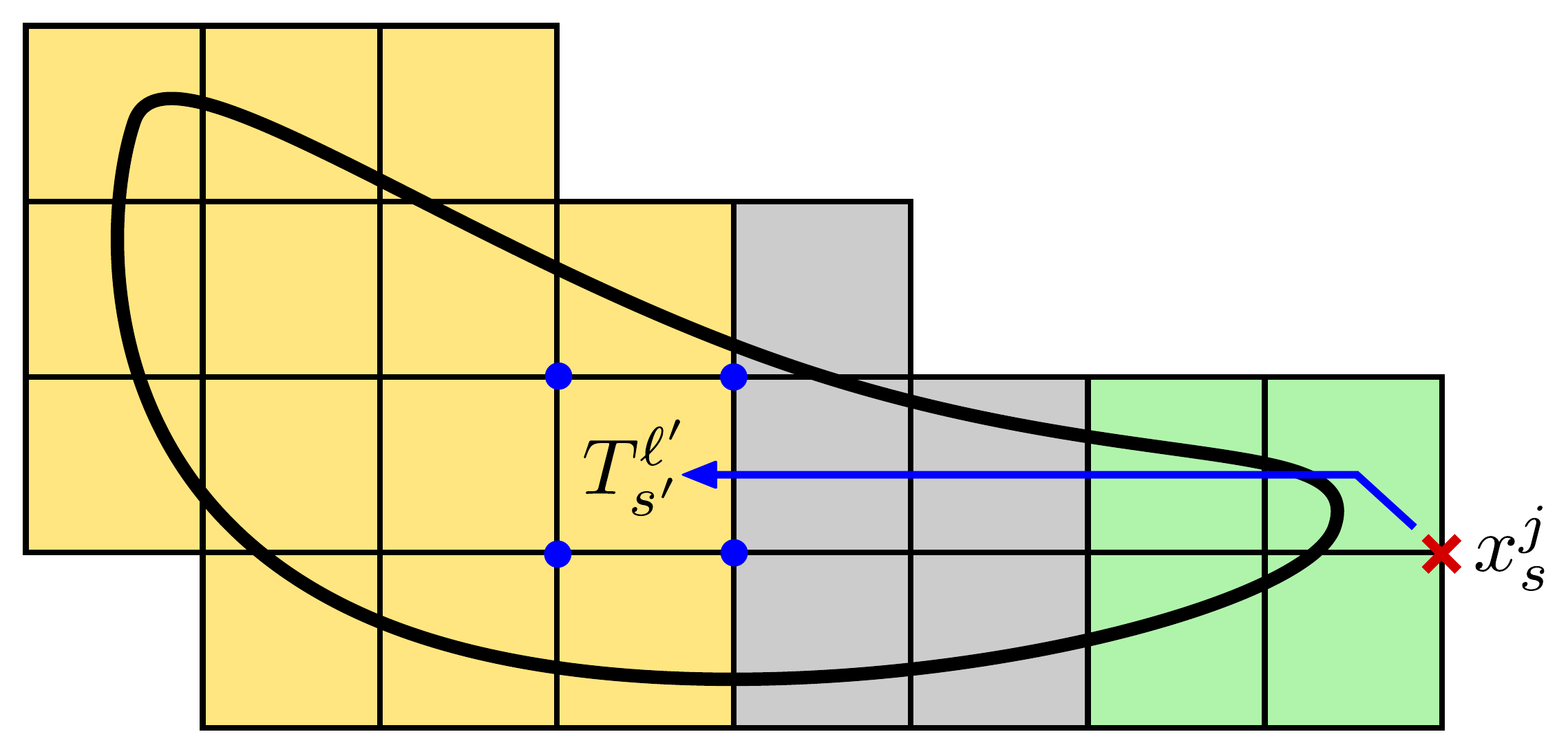}
\end{subfigure}

\begin{subfigure}{0.5\textwidth}
     \begin{small}
      \hfill
      {\color[RGB]{255,230,128} \rule{10pt}{8pt}} $D_{s'}$\hfill
            {\color[RGB]{204,204,204} \rule{10pt}{8pt}} $D_{s''}$\hfill
      {\color[RGB]{176,244,171} \rule{10pt}{8pt}} $D_s$\hfill
     \end{small}
\end{subfigure}

\caption{Illustration of a case where an exterior \ac{dof} is constrained by a root cell in a remote sub-domain. The value at $x^j_s$ is constrained as the extrapolation of the nodal values of cell $T^{\ell'}_{s'}$. In order to set up such constraint it is necessary to send data from sub-domain $D_{s'}$ to sub-domain $D_s$, which are not nearest neighbor sub-domains.}
\label{fig:remote-root}
\end{figure}

\subsection{Implementation overview}\label{sec:overview}

In the following sections (Sects.~\ref{sec:par-cell-agg}--\ref{sec:par-assembly}), we discuss the implementation of the main algorithmic phases of the distributed-memory \ac{agg} method. We start by presenting the distributed version of the cell aggregation method presented in Alg.~\ref{alg:agg_sch} (see Sect.~\ref{sec:par-cell-agg}), where we show the actual computation of the distributed root cell map $R_s$ and the auxiliary maps $K^{\rm own}_s$ and $S^{\rm root}_s$. Then, in Sect.~\ref{sec:data-import}, we detail how to build the communication strategy in order to import data from root cells that are located in remote sub-domains. 
Finally,  we detail the data distribution model used for representing the linear system of equations in the distributed-memory implementation (Sect.~\ref{sec:data-layout}) and the parallel \ac{fe} assembly associated with the \ac{agfe} space (Sect.~\ref{sec:par-assembly}).

\subsection{Parallel cell aggregation}\label{sec:par-cell-agg}

The goal of this section is to parallelize the cell aggregation method presented in  Alg.~\ref{alg:agg_sch}. To this end, the most naive approach would be to run Alg.~\ref{alg:agg_sch} on the local cells $\mathcal{T}^{\rm L}_s$ in each processor independently, trying to greedily associate each cut cell with a root cell in the same sub-domain. However, this approach does not work when a sub-domain has cut cells that cannot be aggregated to any local cell (i.e., the cut cell does not share a face with any local internal cell), and therefore, it is not general enough. Moreover, even when the method is algorithmically applicable  (e.g., for domain decomposition setups where each sub-domain has at least one interior cell), it does not necessarily lead to the same aggregates as the serial method in Alg.~\ref{alg:agg_sch}. 
{\rev The mathematical} analysis of the \ac{agg} method in \citep{Badia2018} states that the approximation properties of the \ac{agfe} spaces are negatively affected as the aggregate size increases. Thus, this naive parallel strategy that can potentially lead to large aggregates has to be avoided.

\begin{algorithm}[t!]
\caption{Parallel cell aggregation algorithm}
\label{alg:par_agg_sch}
\SetKwIF{Unless}{ElseIf}{Else}{unless}{do}{else if}{else}{end}
\vspace*{0.5em}

Initialize $R_s(\ell)\leftarrow -1 $, $S^{\rm root}_s(\ell)\leftarrow -1 $  and $N_s(\ell)\leftarrow -1 $ for all local cells $\ell\in\mathcal{L}^{\rm L}_s$\nllabel{lin:par_agg_sch_0}\;
Set $R_s(\ell)\leftarrow K_s(\ell) $, $S^{\rm root}_s(\ell)\leftarrow s $ and $N_s(\ell)\leftarrow K_s(\ell)$ for all interior local cells $\ell\in\mathcal{L}^{\rm L}_s\cap\mathcal{L}^{\rm in}_s$\;
Recover the values of $R_s(\ell)$ and $S^{\rm root}_s(\ell)$ in ghost cells $\ell\in\mathcal{L}^{\rm G}_s$ with a nearest neighbor exchange \nllabel{lin:goto-par-agg}\nllabel{lin:par_agg_sch_1}\;
Set as touched all cell ids $\ell\in\mathcal{L}_s$ such that $R_s(\ell)\neq -1$\nllabel{lin:par_agg_sch_10} \;
\For(\nllabel{lin:par_agg_sch_11}){$\ell\in\mathcal{L}^{\rm L}_s\cap\mathcal{L}^{\rm cut}_s$ such that $\ell$ has not been touched yet}{

  $\mathcal{L}^{\rm can}_s\leftarrow \emptyset$\;
  \For(\nllabel{lin:par_agg_sch_00}){ $\ell'\in\mathcal{L}^{\rm nei}_{s}(\ell)$ such that cell id $\ell'$ is already touched }{
    Let $F$  be the face shared by cells $T^\ell_s$ and $T^{\ell'}_s$\;
    \If{ $F\cap\Omega\neq\emptyset$}{
    Add $\ell'$ to the set of candidates $\mathcal{L}^{\rm can}_s$\;
    }
  }

  \If{ $\mathcal{L}^{\rm can}_s \neq\emptyset$}{
   Choose an arbitrary cell id $\ell'\in\mathcal{L}^{\rm can}_s$ (e.g., the one with the closest root cell) \nllabel{lin:par_agg_sch_4}\;
       $R_s(\ell)\leftarrow R_s(\ell')$\;
    $S^{\rm root}_s(\ell)\leftarrow S^{\rm root}_s(\ell')$\;
    $N_s(\ell)\leftarrow K_s(\ell')$\nllabel{lin:par_agg_sch_12}\;
   }

}
\eIf{$R_s(\ell)\neq -1 $ for all cell ids $\ell\in\mathcal{L}_s$}{
$b_s\leftarrow$ true\;}
{
$b_s\leftarrow$ false\;
}

\If (\nllabel{lin:par_agg_sch_2}) { $b_{s'}$ is still false in some sub-domain $s'\in\mathcal{S}$}{
Go to line \ref{lin:goto-par-agg} in all processors\;
}

\end{algorithm}

Instead, we consider a more sophisticated parallel algorithm, which leads to the same aggregates as the ones obtained with the serial aggregation method. This enhanced distributed-memory cell aggregation method is detailed in Alg.~\ref{alg:par_agg_sch} and illustrated in Fig.~\ref{fig:par-aggr-steps} for a simple two dimensional case. The notation $\mathcal{L}^{\rm nei}_s(\ell)$, in line \ref{lin:par_agg_sch_00} of the algorithm represents the set of local cell ids of the face neighbors of cell $T^\ell_s$. That is, the notations  $\mathcal{K}^{\rm neig}(k)$ (see Sect.~\ref{sec:cell-aggr}) and $\mathcal{L}^{\rm nei}_s(\ell)$ both represent the face neighbors of a given cell, but the former in terms of global cell ids and the latter in terms of local cell ids. The parallel cell aggregation method in Alg.~\ref{alg:par_agg_sch} follows essentially the same steps as the serial counterpart plus some communications between neighbor sub-domains, which allow one to aggregate a cell into a root belonging to a different sub-domain. The method consists in three main steps (see Fig.~\ref{fig:par-aggr-steps}). (a) We initialize the value of $R_s$ for all interior local cells and perform a nearest neighbor communication in order to determine the value of $R_s$ at  ghost cells. At this point, we mark all interior (local and ghost) cells as "touched". (b) As in the serial case, we perform a loop in (local) cut cells. For a given cut cell, we define its root cell taking the root cell of one of its face neighbors that have been already marked as "touched" (if such a {\rev neighbor exists}). (c) After the loop, we perform a nearest neighbor communication in order to update the value of $R_s$ at ghost cells. At this point, we mark as "touched" all (local and ghost) cells for with $R_s$ is already set with a valid root cell id. Finally, we repeat steps (b) and (c) in all processors as many times as needed until all the cells are assigned with a root cell id. {\rev In the numerical experiments,  3 repetitions of steps (b) and (c) were needed to aggregate all cut cells for all cases.} 

 Note that, during the execution of Alg.~\ref{alg:par_agg_sch}, we also build the map $S^{\rm root}_s:\mathcal{L}_s\rightarrow\mathcal{S}$, which was formally presented in Sect.~\ref{sec:par-fe-spaces}. In addition, we build an auxiliary map $N_s:\mathcal{L}_s\rightarrow\mathcal{K}$, which is the distributed version of the map $N:\mathcal{K}\rightarrow\mathcal{K}$ presented in Sect.~\ref{sec:cell-aggr}. Precisely, $N_s(\ell)$ returns the global id of the cell that was selected in line \ref{lin:par_agg_sch_4} of the algorithm, i.e., the next cell in the path towards the root cell. This information will be required later in the next section (Sect.~\ref{sec:data-import}) in order to build the communication strategy for importing data from remote root cells.


%
\begin{figure}[t!]
  \centering
  \begin{subfigure}{0.99\textwidth}
    \centering
    \begin{small}
      \begin{tabular}{llll}
         \tikz{\fill[fill=myellow]  (0,0) rectangle (1.4em,1.4em);} touched
         &
         \tikz{\fill[fill=gray20]    (0,0) rectangle (1.4em,1.4em);} untouched
      \end{tabular}
    \end{small}
  \end{subfigure}
  \par
  \begin{subfigure}{0.24\textwidth}
    \centering
    \includegraphics[width=0.9\textwidth]{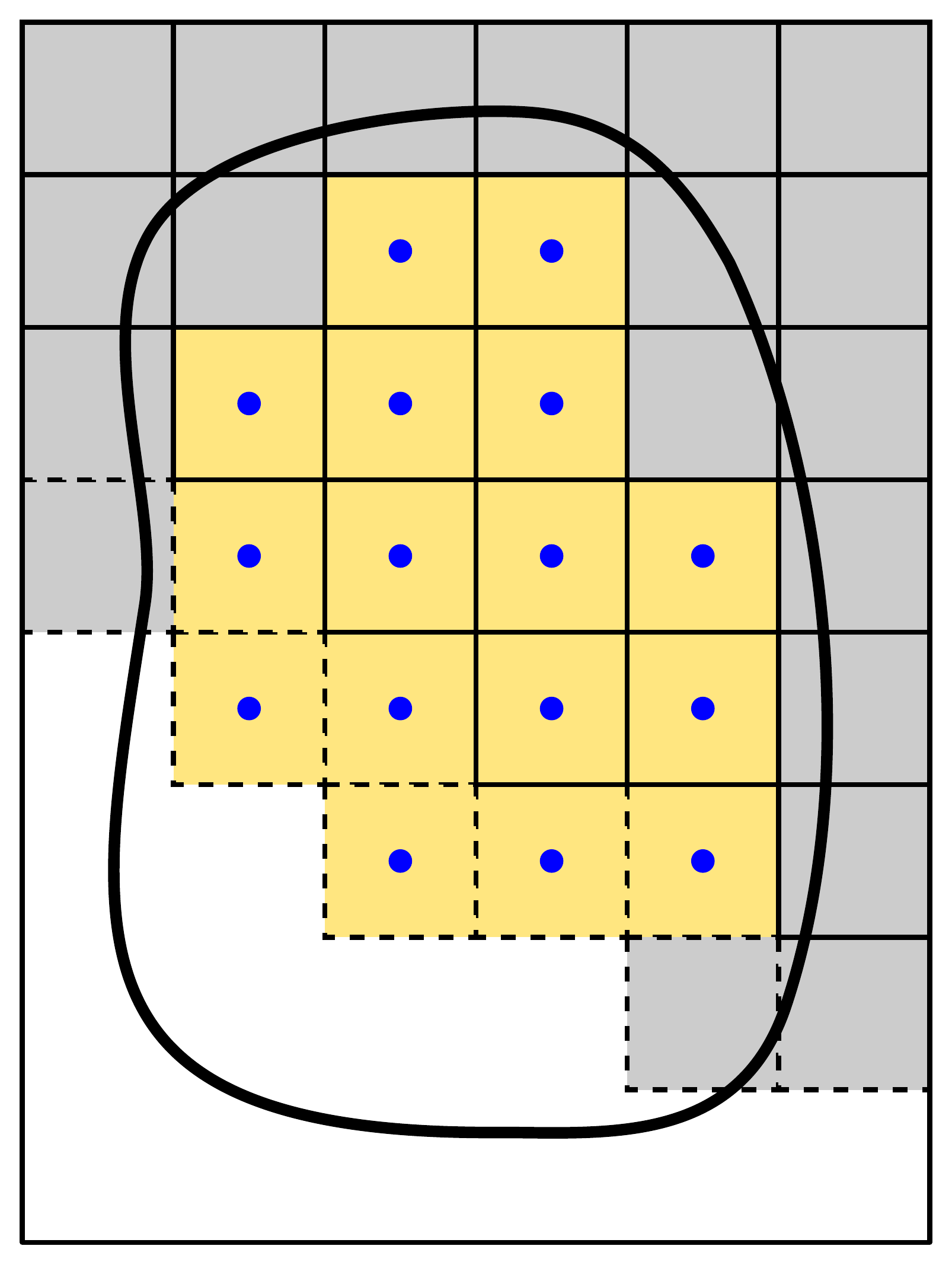}

    \includegraphics[width=0.9\textwidth]{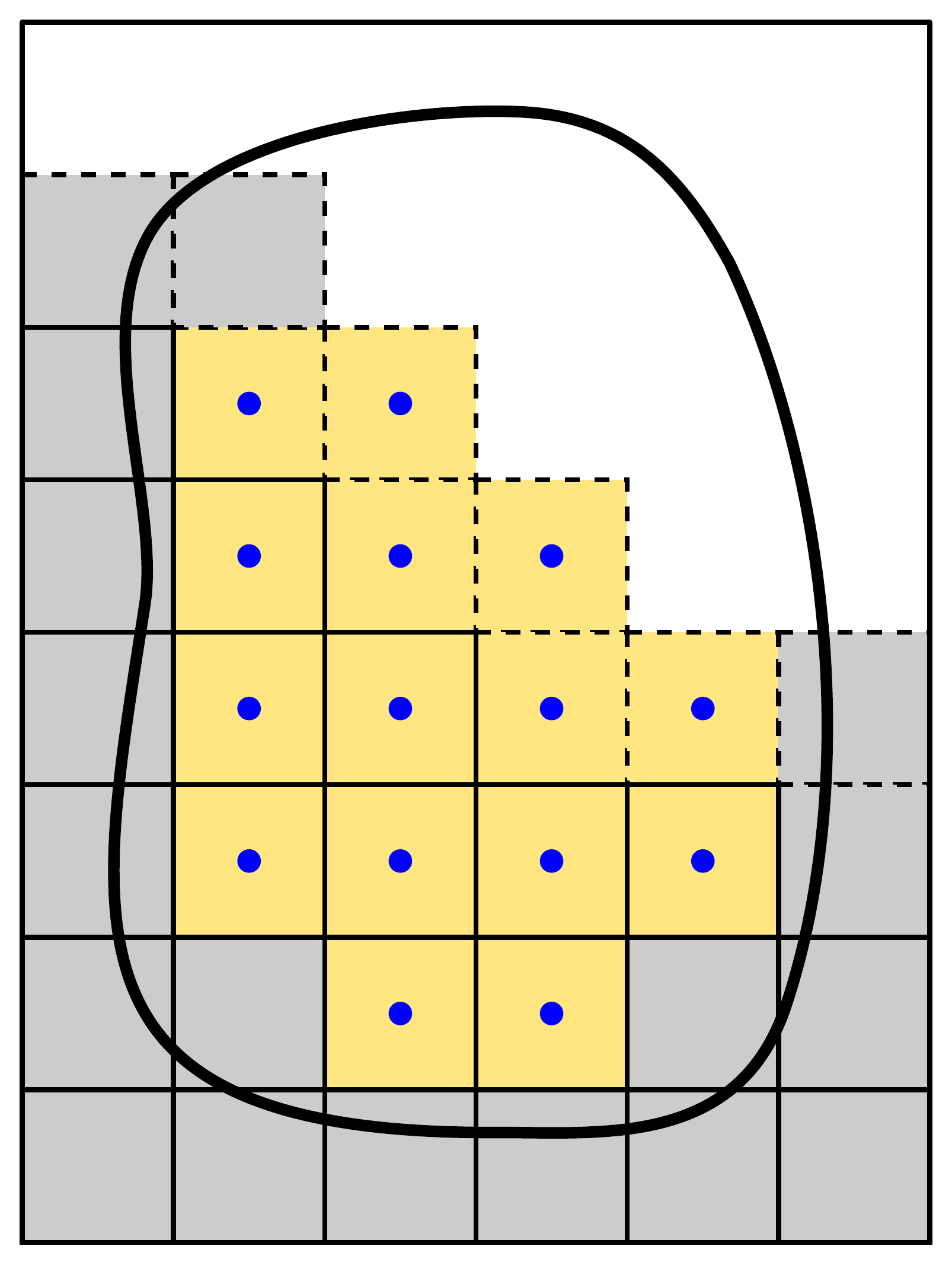}
    \caption{Step 1.}
    \label{fig:par-aggr-steps-a}
  \end{subfigure}
  \begin{subfigure}{0.24\textwidth}
    \centering
    \includegraphics[width=0.9\textwidth]{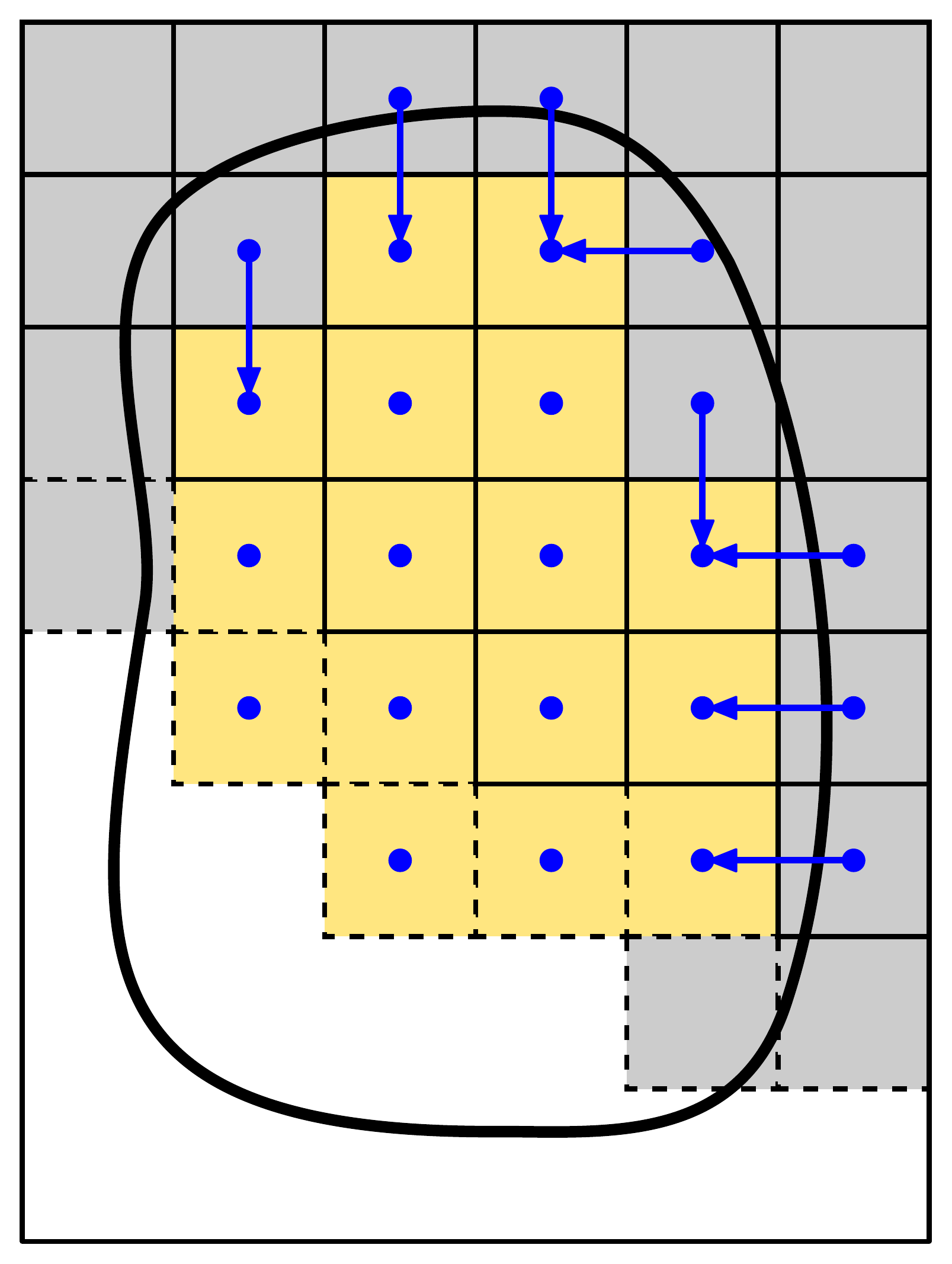}

    \includegraphics[width=0.9\textwidth]{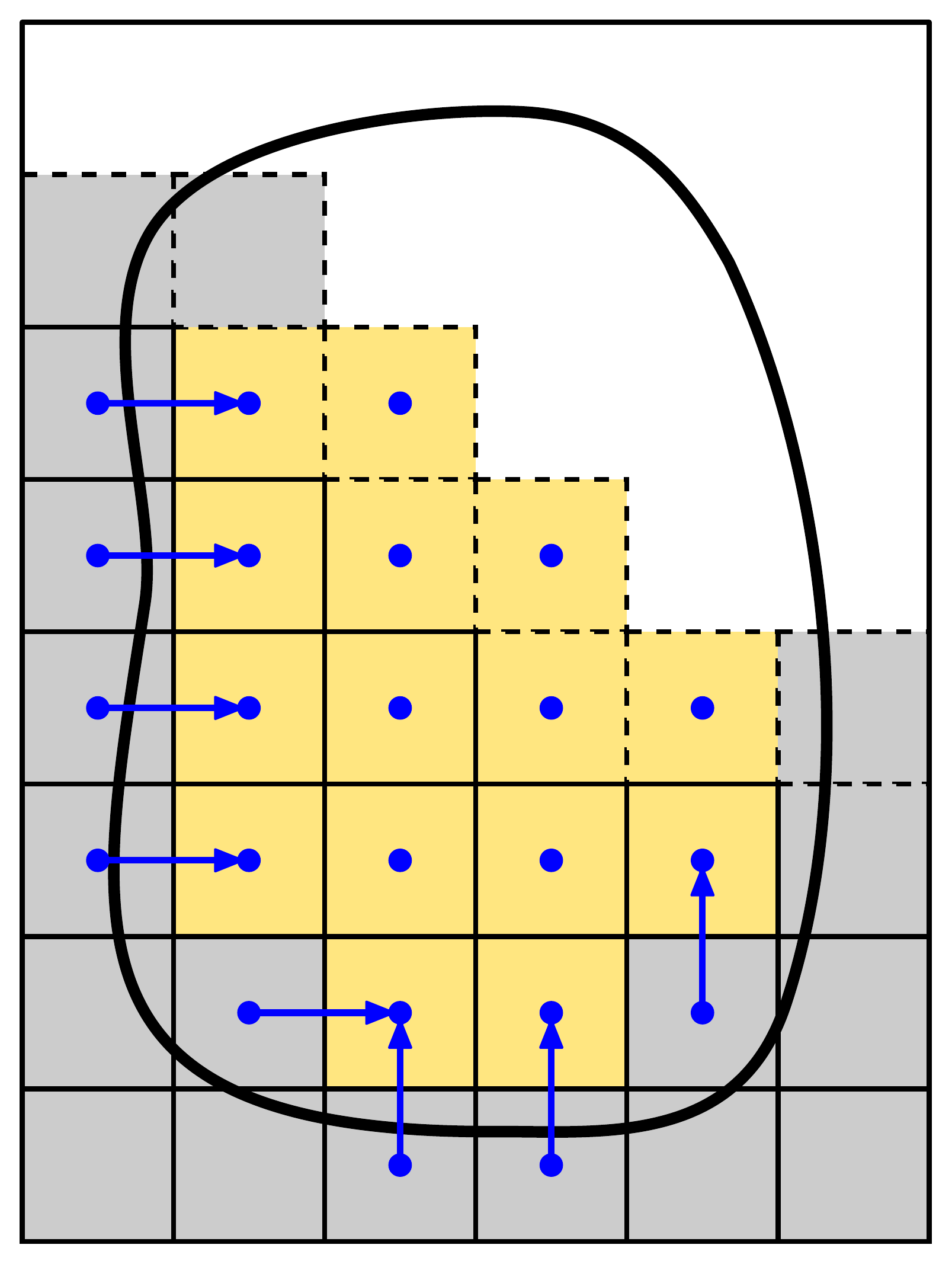}
    \caption{Step 2.} 
    \label{fig:par-aggr-steps-b}
  \end{subfigure}
  \begin{subfigure}{0.24\textwidth}
    \centering
    \includegraphics[width=0.9\textwidth]{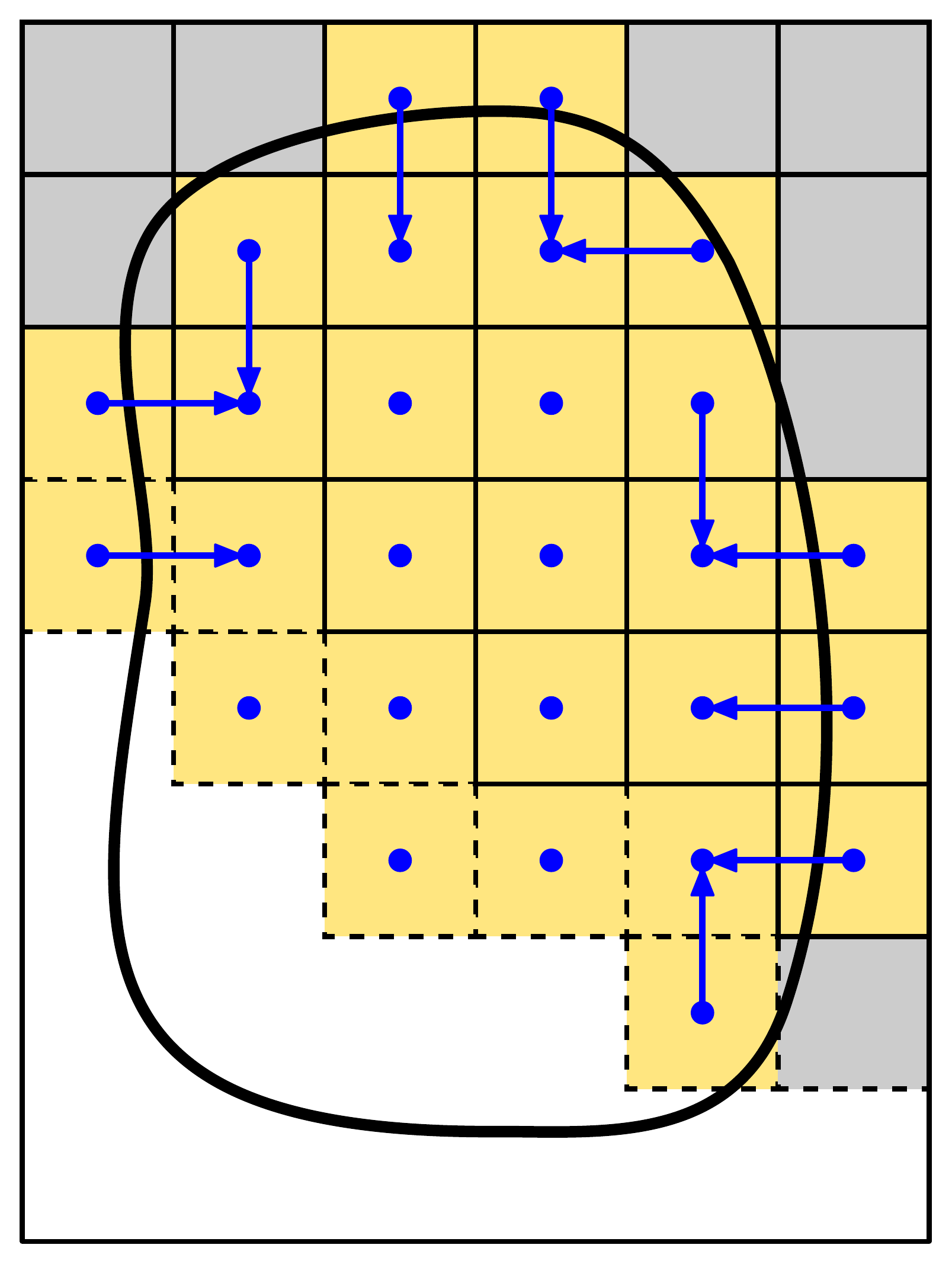}

    \includegraphics[width=0.9\textwidth]{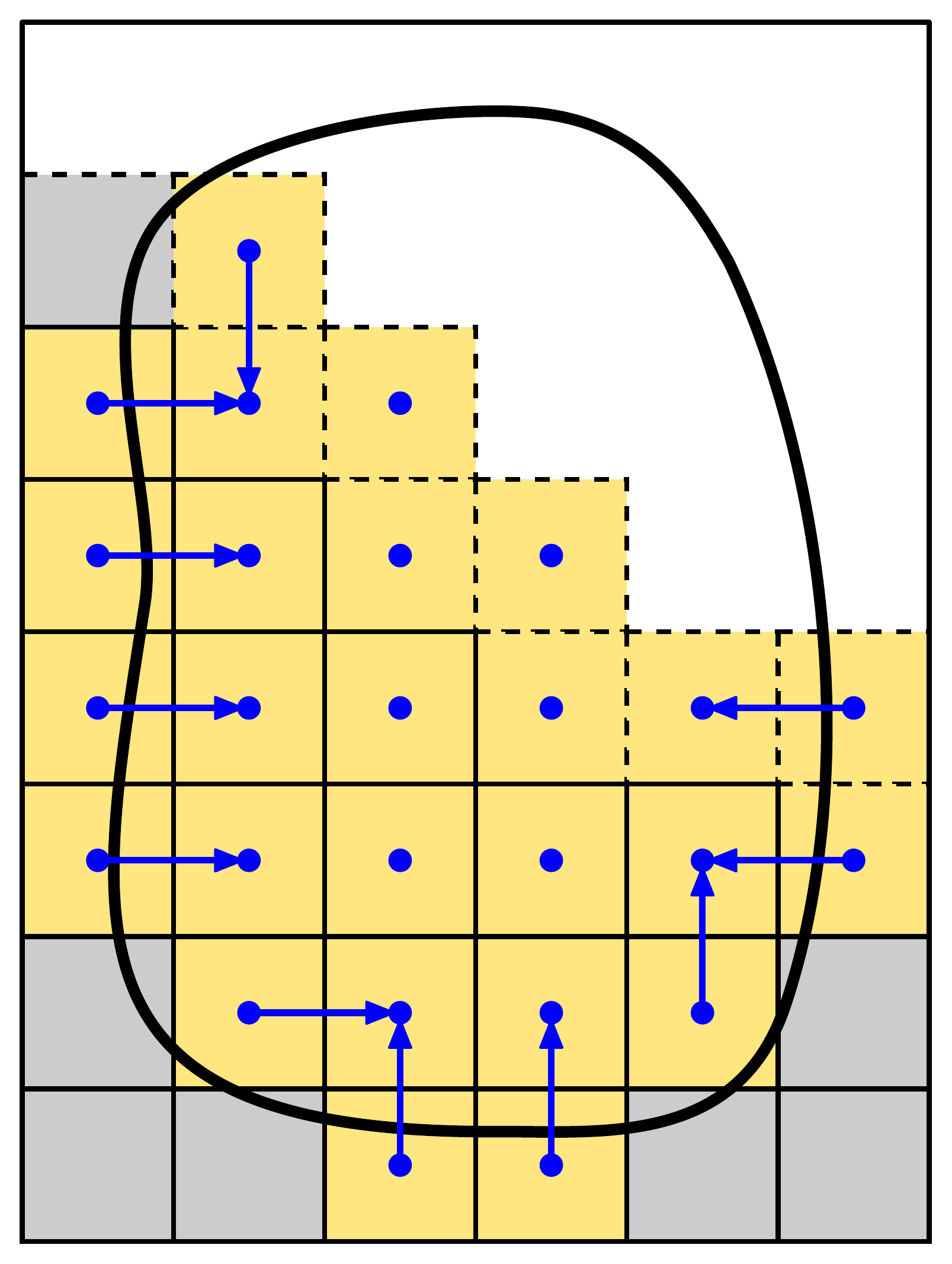}
    \caption{Step 3.} 
    \label{fig:par-aggr-steps-c}
  \end{subfigure}
  \begin{subfigure}{0.24\textwidth}
    \centering
    \includegraphics[width=0.9\textwidth]{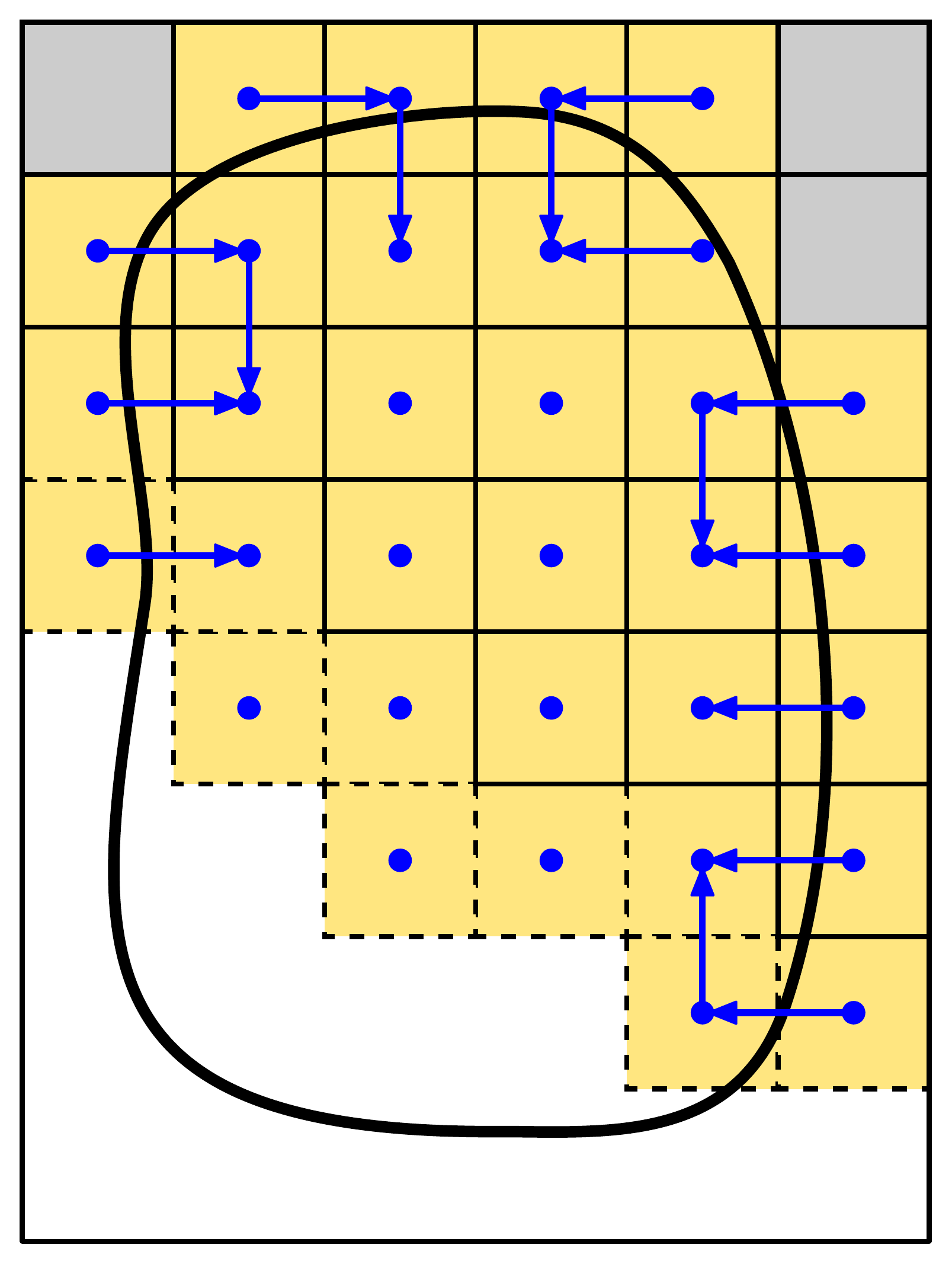}

    \includegraphics[width=0.9\textwidth]{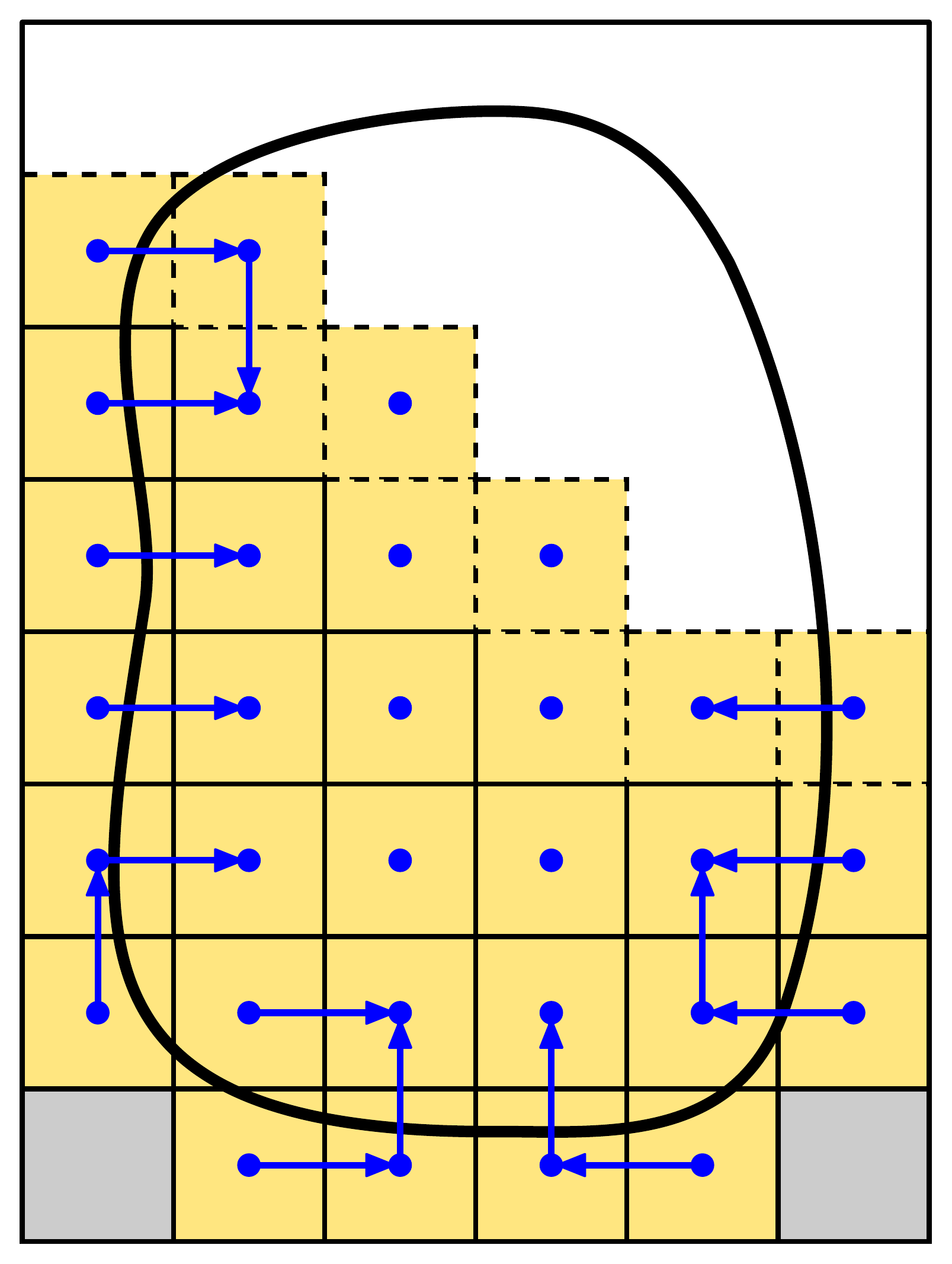}
    \caption{Step 2' (end).} 
    \label{fig:par-aggr-steps-d}
  \end{subfigure}
  \caption{Illustration of the cell aggregation scheme defined in Algorithm \ref{alg:par_agg_sch} for the two sub-domains presented in Fig.~\ref{fig:partition-setup}.  A blue circle indicates cells for which the map $R_s$ is already defined with a valid value. Blue arrows indicate the face neighbor that is selected in line~\ref{lin:par_agg_sch_4} of the algorithm.}
  \label{fig:par-aggr-steps}
\end{figure}

\subsection{Importing data from root cells}\label{sec:data-import}

As we have discussed in Sect.~\ref{sec:par-fe-spaces}, we need to import data associated with remote root cells in order to build the \ac{agg} constraints in Eq.~\eqref{eq:agg_space_s}.  Formally, the set containing the global cell ids of those root cells is defined as $\mathcal{K}^{\rm R}_s\doteq \{k\in \mathcal{K}\setminus\mathcal{K}_s:  k=R_s(\ell) \text{ for some } \ell\in\mathcal{L}_s\}$. For assembling the system matrix and right hand side vector, we need to import the cell-wise nodal coordinates and global \ac{dof} ids for all cells with global ids in $\mathcal{K}^{\rm R}_s$ as it will be justified in Sect. \ref{sec:par-assembly}.
This information is fetched from the corresponding remote sub-domains (as we discuss later) and stored in the local scope of the processor in a separate buffer allocation. In order to address data in this buffer, each cell with global id $k$ in the set  $\mathcal{K}^{\rm R}_s$ is labeled with a new sub-domain local id $z$, defined as the position of $k$ in the list that represents the index set $\mathcal{K}^{\rm R}_s$. This transformation is represented here by a map $Z_s : \mathcal{K}^{\rm R}_s \longrightarrow \mathcal{L}^{\rm R}_s\doteq\{1,\ldots,|\mathcal{K}^{\rm R}_s|\}$ such that, given a cell global id $k\in\mathcal{K}^{\rm R}_s$, it returns the corresponding $z$ id. We denote by $\tilde T^z_s$ the cell associated with label $z\in\mathcal{L}^{\rm R}_s$, namely $\tilde T^z_s\doteq T^k$ with $z=Z_s(k)$. In addition, we denote as $\tilde x^a_{s,z}\doteq x^{a,k}$, $z=Z_s(k)$, the coordinates of the $a$-th node of cell $\tilde T^z_s$ for each $a\in\Lambda(\tilde T^z_s)$. We use an analogous notation for $\tilde g^a_{s,z}\doteq g^{a,k}$, $\tilde{\phi}^a_{s,z}\doteq \phi^{a,k}$, and $\tilde v^a_{s,z}\doteq v^{a,k}$ (i.e., cell-wise global \ac{dof} ids, cell shape functions, and cell  \ac{dof} values of a given \ac{fe} function $v\in\mathcal{V}^{\rm agg}$).

Note that the cells with global id in $\mathcal{K}^{\rm R}_s$ are local cells in the sub-domains with id in the set $\mathcal{Q}^{\rm rcv}_s\doteq\{s'\in\mathcal{S}:\ s'=S^{\rm root}_s(\ell) \text{ for some }\ell\in\mathcal{L}^{\rm cut}_s \}$. That is, we import data from all sub-domains that contain the root cells of the cut cells in the current sub-domain. We note that the root cells of {\em ghost} cut cells also need to be imported as $L_s^{\rm own}(j)$ may be in $\mathcal{L}_s^{\rm G}$ for some $j \in \mathcal{J}^{\rm out}_s$. We also note that we have constructed the map $S^{\rm root}_s$ during the parallel cell aggregation {\rev algorithm} in order to be able to construct the set $\mathcal{Q}^{\rm rcv}_s$. It is also easy to see that the set $\mathcal{R}^{\rm rcv}_{s,s'}\doteq\{k\in\mathcal{K}_{s'}:\ k=R_s(\ell) \text{ for some }\ell\in\mathcal{L}^{\rm cut}_s\ \}$ contains, for a given $s'\in\mathcal{Q}^{\rm rcv}_s$, the global ids of root cells in sub-domain $D_{s'}$ from which we (potentially) need to recover data.

In order to grasp the sets $\mathcal{Q}^{\rm rcv}_s$ and $\mathcal{R}^{\rm rcv}_{s,s'}$, we can consider the cell path that goes from each cut cell to its corresponding root cell, jumping along the face neighbors that are selected during the cell aggregation algorithm (see the blue arrows in Fig~\ref{fig:par-aggr-steps-d}). The set $\mathcal{Q}^{\rm rcv}_s$ contains the ids of the sub-domains, where the cell paths starting in the current sub-domain $D_s$ end. On the other hand, the set $\mathcal{R}^{\rm rcv}_{s,s'}$ contains the ids of the last cell in each path that starts at $D_s$ and ends at $D_{s'}$. The practical construction of the sets $\mathcal{Q}^{\rm rcv}_s$ and $\mathcal{R}^{\rm rcv}_{s,s'}$ from the maps $R_s$ and $S^{\rm root}_s$ is indicated in Alg.~\ref{alg:par_agg_rec_direct}. Note that this procedure is done at each sub-domain independently and does not require communication at all.

\begin{algorithm}[t!]
\caption{Parallel direct path reconstruction}
\label{alg:par_agg_rec_direct}
\vspace*{0.5em}
$\mathcal{Q}^{\rm rcv}_s\leftarrow\emptyset$\;
$\mathcal{R}^{\rm rcv}_{s,s'}\leftarrow\emptyset$ for all $s'\in\mathcal{S}$\;
\For{$\ell\in\mathcal{L}^{\rm cut}_s$}{
$s'\leftarrow S^{\rm root}_s(\ell)$\;
\If{ $s'\neq s$  }{
Add $s'$ to $\mathcal{Q}^{\rm rcv}_s$\;
Add $R_s(\ell)$ to $\mathcal{R}^{\rm rcv}_{s,s'}$\;
}

}
\end{algorithm}

Once the sets $\mathcal{Q}^{\rm rcv}_s$ and $\mathcal{R}^{\rm rcv}_{s,s'}$ are built, we know from which sub-domains, and from which cells in those sub-domains, we need to recover data. However, we still need to know to which sub-domains we need to send data from the current sub-domain, and which is the data to be sent to each of those sub-domains. The "send" counterparts of the sets $\mathcal{Q}^{\rm rcv}_s$ and $\mathcal{R}^{\rm rcv}_{s,s'}$ are formally defined as $\mathcal{Q}^{\rm snd}_s\doteq\{s'\in\mathcal{S}:\ s=S^{\rm root}_{s'}(\ell') \text{ for some }\ell' \in\mathcal{L}^{\rm cut}_{s'}\}$ and, for a given  $s'\in\mathcal{Q}^{\rm snd}_s$, $\mathcal{R}^{\rm snd}_{s,s'}\doteq\{k\in\mathcal{K}_{s}:\ k=R_{s'}(\ell') \text{ for some }\ell'\in\mathcal{L}^{\rm cut}_{s'} \}$. We can interpret the set $\mathcal{Q}^{\rm snd}_s$ as the ids of the sub-domains, where all the cell paths ending at $D_s$ start. On the other hand, the set $\mathcal{R}^{\rm snd}_{s,s'}$ contains the ids of the last cell in each path that ends at $D_s$ and starts at $D_{s'}$. In contrast to the sets $\mathcal{Q}^{\rm rcv}_s$ and $\mathcal{R}^{\rm rcv}_{s,s'}$, it is clear that we cannot construct the sets $\mathcal{Q}^{\rm snd}_s$ and $\mathcal{R}^{\rm snd}_{s,s'}$ without communications since their definitions are in terms of the maps $R_{s'}$ and $S^{\rm root}_{s'}$ seen from remote sub-domains $D_{s'}$ and thus not available at the current one $D_s$. However, thanks to the construction of the auxiliary map $N_s:\mathcal{L}_s\rightarrow\mathcal{K}$ during the parallel cell aggregation {\rev algorithm}, we can build the sets $\mathcal{Q}^{\rm snd}_s$ and $\mathcal{R}^{\rm snd}_{s,s'}$ by using only standard communications between nearest neighbor sub-domains.

The detailed construction of $\mathcal{Q}^{\rm snd}_s$ and $\mathcal{R}^{\rm snd}_{s,s'}$ is given in Alg.~\ref{alg:par_agg_rec_inv}. In the algorithm, we create a tuple for each cell path that starts at the current sub-domain (i.e., one tuple for each local cut cell). These tuples have three items each (see line \ref{lin:par_agg_rec_tuple_init}): the global id of the first cell in the path, the global id of the next cell in the path, and the current sub-domain id. Conceptually, we gradually move forward these tuples along their cell paths until they arrive to the final cell in the path, which can be possibly in another sub-domain. At the end of the process each sub-domain collects all the tuples corresponding to the cell paths that end in this sub-domain. With this information each sub-domain is able to reconstruct the sets  $\mathcal{Q}^{\rm snd}_s$ and $\mathcal{R}^{\rm snd}_{s,s'}$. Note that Alg.~\ref{alg:par_agg_rec_inv} only involves communications with nearest neighbor sub-domains.

Finally, we detail in Alg.~\ref{alg:par_import_data} how we import the required data from remote sub-domains by using the sets $\mathcal{Q}^{\rm snd}_s$ and $\mathcal{R}^{\rm snd}_{s,s'}$ and $\mathcal{Q}^{\rm rcv}_s$ and $\mathcal{R}^{\rm rcv}_{s,s'}$. 

\begin{algorithm}[ht!]
\caption{Parallel inverse path reconstruction}
\label{alg:par_agg_rec_inv}
\vspace*{0.5em}

$\mathcal{A}_s\leftarrow\emptyset$\;
\For{$\ell\in\mathcal{L}^{\rm L}_s\cap\mathcal{L}^{\rm cut}_s$}{
Add the tuple $(K_s(\ell),\ N_s(\ell),\ s)$ to the set $\mathcal{A}_s$\nllabel{lin:par_agg_rec_tuple_init}\;
}
\While(\nllabel{lin:par_agg_rec:goto}){there is a tuple $(k,n,z)\in\mathcal{A}_s$ such that $n\in\mathcal{K}^{\rm L}_s$ and $L_s(n)\in\mathcal{L}^{\rm cut}_s$ }{
\For{$(k,n,z)\in\mathcal{A}_s$ }{
\If{$n\in\mathcal{K}^{\rm L}_s$ and $L_s(n)\in\mathcal{L}^{\rm cut}_s$}{
$n'\leftarrow N_s(L_s(n))$\;
Replace the tuple $(k,n,z)$ with $(k,n',z)$ in the set $\mathcal{A}_s$\;
    }
  }
}

$\mathcal{Q}^{\rm snd}_s\leftarrow\emptyset$\;
$\mathcal{R}^{\rm snd}_{s,s'}\leftarrow\emptyset$ for all $s'\in \mathcal{Q}^{\rm snd}_s$  \;
\For{ $(k,n,z)\in\mathcal{A}_s$  }{
\If{$n\in\mathcal{K}^{\rm L}_s$ and $L_s(n)\in\mathcal{L}^{\rm in}_s$}{
Add $z$ to the set $\mathcal{Q}^{\rm snd}_s$\;
Add $k$ to the set $\mathcal{R}^{\rm snd}_{s,z}$\;
Delete the tuple $(k,n,z)$ from $\mathcal{A}_s$\;
}
}
\For{ $s'\in\mathcal{S}^{\rm nei}_s$}{
$\mathcal{A}^{\rm snd}_{s,s'}\leftarrow\emptyset$\;
\For{ $(k,n,z)\in\mathcal{A}_s$  }{
\If{$s'=S_s(L_s(n))$}{
Add the tuple $(k,n,z)$ to $\mathcal{A}^{\rm snd}_{s,s'}$\;
}
}
Send the set $\mathcal{A}^{\rm snd}_{s,s'}$ to processor $s'$ \nllabel{lin:par_agg_rec:send}\;

}
Delete the set $\mathcal{A}_s$\;
$\mathcal{A}_s\leftarrow\emptyset$\;
\For{ $s'\in\mathcal{S}^{\rm nei}_s$}{
Receive the set $\mathcal{A}^{\rm rcv}_{s,s'}$ from processor $s'$ as a counterpart of the send operation in line \ref{lin:par_agg_rec:send}\;
Add all tuples in $\mathcal{A}^{\rm rcv}_{s,s'}$ to the set $\mathcal{A}_s$\;
}
\If{$\mathcal{A}_{s'}\neq\emptyset$ still in some sub-domain $s'\in\mathcal{S}$}{
Go to line \ref{lin:par_agg_rec:goto} in all processors
}
\end{algorithm}

%
%
%
%

\begin{algorithm}[ht!]
\caption{Import required data from remote root cells}
\label{alg:par_import_data}
\vspace*{0.5em}
\For{$s'\in\mathcal{Q}^{\rm snd}_s$}
{
  $\mathcal{B}_{s,s'}\leftarrow\emptyset$\;
  \For{$k\in\mathcal{R}^{\rm snd}_{s,s'}$}{
  $\ell\leftarrow L_s(k)$\;
  Add to $\mathcal{B}^{\rm snd}_{s,s'}$ the pair $(x^a_{s,\ell},\ g^a_{s,\ell} )$ for all $a\in\Lambda(T^\ell_s)$\;
  }
  Send the buffer $\mathcal{B}^{\rm snd}_{s,s'}$ to sub-domain $s'$ \nllabel{lin:par_import_data:snd}\;
}

\For{$s'\in\mathcal{Q}^{\rm rcv}_s$}
{
  Receive the buffer $\mathcal{B}^{\rm rcv}_{s,s'}$ from sub-domain $s'$ as counterpart of the send operations in line \ref{lin:par_import_data:snd}\;

  \For{$k\in\mathcal{R}^{\rm rcv}_{s,s'}$}{
   $z\leftarrow Z_s(k)$\;
   \For{$a\in\Lambda(\tilde T^z_s)$}{
   Get the next pair  $(x,\ g)$ in $\mathcal{B}^{\rm rcv}_{s,s'}$\;
   $\tilde x^{a}_{s,z}\leftarrow x$;  $\tilde g^{a}_{s,z}\leftarrow g$\;
   }
   }
}
\end{algorithm}

\subsection{Linear algebra data distribution layout}\label{sec:data-layout}

Before detailing the \ac{fe} assembly, we need to specify the data model used for representing the matrices and vectors in a distributed-memory context. Here, we consider row-wise partitioned globally addressable matrices and vectors (as we detail in the next paragraph). This choice is motivated by the fact that the vast majority of parallel linear algebra packages (e.g., \petsc~\cite{petsc-web-page,petsc-user-ref}, or \trilinos~\cite{Heroux2005}) use this data distribution model. Other data layouts are also possible, such as sub-assembled distributed matrices for non-overlapping domain decomposition
solvers \cite{Badia.etal.ACME.2013,badia_multilevel_2016}, but they are out of the scope of the current work.

In a row-wise partitioning model, the global matrix and the global vector are partitioned by rows, i.e., each processor owns a consecutive, non-overlapping, range of global rows of $\mathbf{A}$ and $\mathbf{b}$, respectively. Since each row of $\mathbf{A}$ and $\mathbf{b}$ is associated with a global \ac{dof} id in $\mathcal{I}^{\rm in}$, the consecutive row partition implies that each sub-domain $D_s$ owns a consecutive subset of the global \ac{dof} ids in $\mathcal{I}^{\rm in}$, denoted as $\mathcal{I}^{\rm in, O}_s$, where  $\{\mathcal{I}^{\rm in, O}_s\}_{s\in\mathcal{S}}$ is a partition of $\mathcal{I}^{\rm in}$. Note that we do not need to partition the set of constrained \ac{dof} ids $\mathcal{I}^{\rm out}$ since, in our approach, these \acp{dof} are removed from the global linear system. On the other hand, we consider linear algebra data structures that are {\em globally addressable}. That is, a given processor with sub-domain id $s\in\mathcal{S}$ may contribute to a global matrix row $[\mathbf{A}]_{gi}$, $i\in\mathcal{I}^{\rm in}$ (or a global vector entry $[\mathbf{b}]_g$), which is not necessarily owned by it, namely $g\in\mathcal{I}^{\rm in}\setminus\mathcal{I}^{\rm in, O}_s$. This is possible since, once the corresponding entries of the matrix and vector are set in the assembly process, a final communication stage transfers all (touched) values $[\mathbf{A}]_{gi}$, $[\mathbf{b}_g]$ such that $g\in\mathcal{I}^{\rm in}\setminus\mathcal{I}^{\rm in, O}_s$ to the corresponding remote processor that owns the global id  $g$. This communication phase is typically resolved by the linear algebra library and, thus, we do not cover it here for simplicity.

Even though  globally addressable linear algebra objects allow one to assign entries that are not owned by the current processor, it is desirable to minimize the number of such assignments in order to reduce inter processor communication. To this end, we define the set of owned global \ac{dof} ids $\mathcal{I}^{\rm in, O}_s$ following the usual approach in \ac{fe} applications. Let $\mathcal{I}^{\rm in}_s\doteq\{i\in\mathcal{I}^{\rm in}:\ i=g^a_{s,\ell} \text{ for some } a\in\Lambda(T^\ell_s) \text{ and } \ell\in\mathcal{L}^{\rm L}_s\}$ be the set of interior global \ac{dof} ids touched only by local cells of sub-domain $D_s$ and let $\mathcal{B}_s(i)\doteq\{\ell\in\mathcal{L}_s:\ i=g^a_{s,\ell} \text{ for some } a\in\Lambda(T^\ell_s) \}$ be the ids of the cells around a given id $i\in\mathcal{I}^{\rm in}_s$. With this notation, we define the set of owned global \ac{dof} ids as $\mathcal{I}^{\rm in, O}_s\doteq\{i\in\mathcal{I}^{\rm in}_s:\ s\leq S_s(\ell) \text{ for all } \ell \in \mathcal{B}_s(i)\}$. Note that each processor owns all \ac{dof} ids which are surrounded solely by cells in $\mathcal{T}_{s}^{\rm L}$ plus some of the \ac{dof} ids that are on the interface between cells in $\mathcal{T}_{s}^{\rm L}$ and $\mathcal{T}_{s}^{\rm G}$. With this model, most of the cells in $\mathcal{T}^{\rm L}_s$ contribute to locally owned \ac{dof} ids in $\mathcal{I}^{\rm in, O}_s$ during the assembly process (i.e., most of the values set from sub-domain $D_s$ will correspond to values owned by $D_s$ if we assume that the number of \acp{dof} in the interior of the sub-domain $D_s$ is larger than number of \acp{dof} on its boundary). On the other hand, for \ac{dof} ids $i\in\mathcal{I}^{\rm in}_s$ that are located at the interface between cells in $\mathcal{T}_{s}^{\rm L}$ and $\mathcal{T}_{s}^{\rm G}$ (i.e., at sub-domain boundaries), we define the owner sub-domain as the sub-domain with smaller id among all sub-domain that share the \ac{dof} id $i$.  We have used this tie-breaker strategy in order to decide which sub-domain is the owner for interface \acp{dof}, but other criteria can be also considered as long as they are consistent across all processors and require no communication.

Even though it is not strictly required, it is convenient, in order to facilitate the interaction with the parallel linear algebra package, that the sets $\{\mathcal{I}^{\rm in, O}_s\}_{s\in\mathcal{S}}$ contain contiguous, sequentially increasing \ac{dof} ids, namely $\mathcal{I}^{\rm in, O}_1=\{1,\ldots,m^{\rm in}_1\}$, $\mathcal{I}^{\rm in, O}_2=\{m^{\rm in}_1+1,\ldots,m^{\rm in}_1+m^{\rm in}_2\}$, etc., where $m^{\rm in}_s\doteq|\mathcal{I}^{\rm in, O}_s|$ for $s\in\mathcal{S}$.
Methods for the generation of cell-wise global \ac{dof} ids $g^a_{s,\ell}$ that result in sets $\{\mathcal{I}^{\rm in, O}_s\}_{s\in\mathcal{S}}$ fulfilling this property are well known in distributed-memory \ac{fe} analysis (see, e.g., \cite{bangerth_algorithms_2012} for the algorithm implemented in \dealii~and \cite{Badia2019b} for the algorithm implemented in \FEMPAR). These techniques can readily be used to generate the cell-wise global \ac{dof} ids $g^a_{s,\ell}$ in the context of the  \ac{agg} method by considering the following steps. First, the values $g^a_{s,\ell}$ are distributed over interior cell ids $\ell\in\mathcal{L}^{\rm in}_s$ by applying the conventional algorithms \cite{bangerth_algorithms_2012,Badia2019b} restricted to the interior portion of the mesh (i.e., as though cut cells were removed from the mesh). In a second phase, global \ac{dof} ids $g^a_{s,\ell}$ are assigned to the vertices, edges and faces of cut cells that are in contact with an interior cell (if any) by taking the corresponding values of the neighbor interior cell. Note that, at the end of  this process, the cell-wise global \ac{dof} ids $g^a_{s,\ell}$ are not defined in the interior of cut cells nor on the parts of their boundaries that are not in contact with interior cells. This is not a problem since this information is not required (see Sect.~\ref{sec:par-assembly}).

\subsection{Distributed-memory finite element assembly} \label{sec:par-assembly}

Finally, we discuss the assembly of the global matrix $\mathbf{A}$ and vector $\mathbf{b}$ in a distributed-memory context. The assembly at the current sub-domain $D_s$ is performed with a loop in local cells in $\mathcal{T}^{\rm L}_s$. At each cell $T^\ell_s\in\mathcal{T}^{\rm L}_s$, we compute the corresponding elemental matrix and vector, namely $\mathbf{A}^\ell_s$ and $\mathbf{b}^\ell_s$, whose entries are defined as $[\mathbf{A}^\ell_s]_{ab}\doteq \mathrm{a}^\ell_s(\phi^a_{s,\ell},\phi^b_{s,\ell})$ and $[\mathbf{b}^\ell_s]_a\doteq \mathrm{b}^\ell_s(\phi^a_{s,\ell})$, where $\mathrm{a}^\ell_s(\cdot,\cdot)$ and $\mathrm{b}^\ell_s(\cdot)$ represent the restriction of the bi-linear form $\mathrm{a}(\cdot,\cdot)$ and the linear form $\mathrm{b}(\cdot)$ to cell $T^\ell_s$, resp.  Then, the contributions of the elemental matrix and vector $\mathbf{A}^\ell_s$ and $\mathbf{b}^\ell_s$ are added to the corresponding entries of the global matrix and vector $\mathbf{A}$ and $\mathbf{b}$.

The assembly operations to be performed are closely related to the ones detailed in Sect. \ref{sec:serial-assembly} for the serial case. In order to grasp the main differences between the serial and the distributed memory implementations, we rewrite the definition of $\mathcal{V}^{\rm agg}_s$ in Eq. \eqref{eq:agg_space_s} using the notation introduced in Sect. \ref{sec:data-import}, namely
\begin{equation}\label{eq:agg_space_s_bis}
\mathcal{V}^{\rm agg}_s \doteq \{ v_s \in \mathcal{V}^{\rm std}_s  :   v^j_s = \sum_{g\in \mathcal{M}_s(j)} C^{jg}_s\ [v^g]_s \text{ for any } j\in \mathcal{J}^{\rm out}_s \},
\end{equation}
 where $\mathcal{M}_s(j)\subset\mathcal{I}^{\rm in}$ denotes the set containing the global \ac{dof} ids that constrain the value associated with a given sub-domain local id $j\in\mathcal{J}^{\rm out}_s$. The values $C^{jg}_s$, for $g\in\mathcal{M}_s(j)$ and $j\in\mathcal{J}^{\rm out}_s$, are the coefficients of the cell aggregation constraints that are needed to be computed in sub-domain $D_s$, and the value $[v^g]_s$  represents a variable that stores (in the local scope of sub-domain $D_s$) the value of $v^g$ for each master global \ac{dof} id $g\in\mathcal{M}_s(j)$ and each $j\in\mathcal{J}^{\rm out}_s$. These quantities are formally defined as
\begin{equation}    
\left.\begin{array}{l}
  \mathcal{M}_s(j) \doteq\{g^a_{s,\ell}\}_{a\in\Lambda(T^\ell_s)},\quad C^{jg}_s \doteq \phi^b_{s,\ell} (x^j_s),\quad [v^g]_s \doteq v^b_{s,\ell} \quad\text{ if } k\in\mathcal{K}_s\\
\mathcal{M}_s(j) \doteq\{\tilde{g}^a_{s,z}\}_{a\in\Lambda(\tilde{T}^z_s)},\quad C^{jg}_s \doteq \tilde{\phi}^b_{s,z} (x^j_s),\quad [v^g]_s \doteq \tilde{v}^b_{s,z} \quad\text{ otherwise}
 \end{array} \right.
\end{equation}
with $\ell=L_s(k)$,  $z=Z_s(k)$, $k = R_s(L^{\rm own}_s(j))$, and
being  $b$ such that $g=g^b_{s,\ell}$ or $g=\tilde{g}^b_{s,z}$.  Note that the definition of these quantities takes one form or another depending on whether the global cell id $k = R_s(L^{\rm own}_s(j))$ belongs  to $\mathcal{K}_s$ or not (i.e., if the root cell $T^k$ associated with $j$ is a locally relevant cell in $D_s$ or not). If $k\in\mathcal{K}_s$, the quantities can be computed from data stored in the local scope of the current processor. Otherwise, if $k\notin\mathcal{K}_s$ (which implies $k\in\mathcal{K}^{\rm R}_s$), the quantities need to be computed from the specially allocated variables $\tilde{g}^b_{s,z}$, $\tilde{x}^b_{s,z}$, and $\tilde{v}^b_{s,z}$ presented in Sect.~\ref{sec:data-import}.

The main differences between the serial and distributed-memory assembly are readily grasped by comparing formulas \eqref{eq:aggr-FES} and \eqref{eq:agg_space_s_bis}.
The first main difference is the following. Note that, in both the serial and distributed implementations, the master \acp{dof} are represented with global \ac{dof} ids. However, the constrained \acp{dof} are represented with global ids in the serial implementation and with sub-domain local ids in the distributed-memory version. For this reason, it is not needed to generate a global numbering for the constrained \acp{dof}. The second major difference is that, in the distributed memory implementation, one has to import data from remote sub-domains in order to compute the quantities $\mathcal{M}_s(j)$, $C^{jg}_s$, and $[v^g]_s$, for $g\in\mathcal{M}_s(j)$, if the root cell associated with the constrained \ac{dof} id $j\in\mathcal{J}^{\rm out}_s$ is not a locally relevant cell.  The assembly operations to be performed in the distributed-memory case are obtained by modifying the operations of the serial case (Table \ref{tab:assembly-cases}) in accordance with these two main differences and taking into account that cells are represented by sub-domain local ids in the parallel implementation. The result is detailed in Table~\ref{tab:par-assembly-cases}.

\begin{table}[ht!]

\setlength{\lrow}{0.6em}

\begin{tabular}{p{0.25\textwidth}p{0.7\textwidth}}
\toprule
Case & Assembly operation\\
\midrule
Right hand side vector\\
$j^a_{s,\ell}\in\mathcal{J}^{\rm in}_s$ & $[\mathbf{b}]_i \text{ += }   [\mathbf{b}^\ell_s]_a$, $i=g^a_{s,\ell}$ \quad (standard case)\\[\lrow]
$j^a_{s,\ell}\notin\mathcal{J}^{\rm in}_s$ & $[\mathbf{b}]_i \text{ += }   C^{ni}_s\ [\mathbf{b}^\ell_s]_a$,  $n=j^a_{s,\ell}$,  $\forall i\in\mathcal{M}_s(n)$\\[2\lrow]
System matrix\\
$j^a_{s,\ell}\in\mathcal{J}^{\rm in}_s$ and  $j^b_{s,\ell}\in\mathcal{J}^{\rm in}_s$ & $[\mathbf{A}]_{ij}\text{ += }  [\mathbf{A}^\ell_s]_{ab}$, $i=g^a_{s,\ell}$, $j=g^b_{s,\ell}$ \quad (standard case)\\[\lrow]
$j^a_{s,\ell}\notin\mathcal{J}^{\rm in}_s$ and  $j^b_{s,\ell}\in\mathcal{J}^{\rm in}_s$ & $[\mathbf{A}]_{ij}\text{ += }  C^{ni}_s\ [\mathbf{A}^\ell_s]_{ab}$, $n=j^a_{s,\ell}$, $j=g^b_{s,\ell}$, $\forall i\in\mathcal{M}_s(n)$\\[\lrow]
$j^a_{s,\ell}\in\mathcal{J}^{\rm in}_s$ and  $j^b_{s,\ell}\notin\mathcal{J}^{\rm in}_s$  & $[\mathbf{A}]_{ij}\text{ += }  C^{mj}_s\ [\mathbf{A}^\ell_s]_{ab}$, $i=g^a_{s,\ell}$ , $m=j^b_{s,\ell}$, $\forall j\in\mathcal{M}_s(m)$\\[\lrow]
$j^a_{s,\ell}\notin\mathcal{J}^{\rm in}_s$ and  $j^b_{s,\ell}\notin\mathcal{J}^{\rm in}_s$ & $[\mathbf{A}]_{ij}\text{ += }  C^{ni}_s C^{mj}_s\ [\mathbf{A}^\ell_s]_{ab}$, $n=j^a_{s,\ell}$, $m=j^b_{s,\ell}$, $\forall i\in\mathcal{M}_s(n)$, $\forall j\in\mathcal{M}_s(m)$\\
\bottomrule
\end{tabular}
\vspace{0.5em}

\caption{Summary of the assembly operations related with the \ac{agfe} space $\mathcal{V}^{\rm agg}$ in the distributed-memory implementation. 
}
\label{tab:par-assembly-cases}
\end{table}

\section{Numerical examples}\label{sec:numericals}

{\rev
\subsection{Experimental environment} \label{sec:experimental-environment}

The numerical experiments are run at the Marenostrum-IV supercomputer~\cite{mn4-ug}, hosted by the Barcelona Supercomputing Center. The Marenostrum-IV is a petascale machine  equipped with 3,456 compute nodes interconnected with the Intel OPA HPC network. Each node has 2x Intel Xeon Platinum 8160 multi-core CPUs, with 24 cores each (i.e. 48 cores per node) and 96 GBytes of RAM.

With respect to the software, we used the MPI-parallel implementation of the AgFEM method available at \FEMPAR{}~\cite{badia_fempar:_2017} linked against \p4est{} v2.0~\cite{burstedde_p4est_2011} as the Cartesian grid manipulation engine, and \petsc{} v3.9.0~\cite{petsc-user-ref} for distributed-memory linear algebra data structures and solvers. These software were compiled with Intel v18.0.5 compilers using system recommended optimization flags and linked against the Intel MPI Library (v2018.4.057) for message-passing and the BLAS/LAPACK available on the Intel MKL library for optimized dense linear algebra kernels. All floating-point operations were performed in IEEE double precision.
}

\subsection{Model problem}\label{sec:model-problem}

As a model problem for the numerical examples, let us consider the Poisson equation with Dirichlet boundary conditions. After scaling with the diffusion term, the equation reads: find $u \in H^1(\Omega)$ such that
\begin{equation}
-\Delta u = f  \quad \text{in } \ \Omega, \qquad
u=g  \quad\text{on } \ \Gamma\doteq\partial\Omega,
\label{eq:PoissonEq}
\end{equation}
where  $f\in H^{-1}(\Omega)$ is the source term and $g\in H^{1/2}(\Gamma)$ is the prescribed value on the  Dirichlet boundary. We denote by $\mathcal{V}^{\rm x}$ any of the \ac{fe} spaces $\mathcal{V}^\mathrm{std}$, or $\mathcal{V}^\mathrm{agg}$, when it is not necessary to distinguish between them. We approximate problem ~\eqref{eq:PoissonEq}  in $\mathcal{V}^{\rm x}$ with the following variational equation: find $u_h\in \mathcal{V}^{\rm x}$ such that $\mathrm{a}(v_h,u_h)=\mathrm{b}(v_h)$ for all $v_h\in \mathcal{V}^{\rm x}$, where
\begin{equation}
\label{eq:weak-PoissonEq}
\begin{array}{l}
\displaystyle \mathrm{a}(v,u)\doteq\int_{\Omega} \nabla v\cdot \nabla u \mathrm{\ d}\Omega + \int_{\Gamma} \left( \tau v u  - v \left(n \cdot \nabla u\right) -  u \left(n \cdot \nabla v\right) \right)  \mathrm{\ d}{\Gamma}, \text{ and}
\\
\displaystyle \mathrm{b}(v)\doteq \int_{\Omega} v f \mathrm{\ d}\Omega + \int_{\Gamma}  \left( \tau  v g  -  \left(n \cdot \nabla v\right)g^\mathrm{D} \right)  \mathrm{\ d}\Gamma,
\end{array}
\end{equation}
with $n$ being the outward unit normal on $\Gamma$.
Note that the forms $\mathrm{a}(\cdot,\cdot)$ and $\mathrm{b}(\cdot)$ include the usual terms resulting from the integration by parts of \eqref{eq:PoissonEq} plus additional terms associated with the weak imposition of Dirichlet boundary conditions with Nitsche's method \cite{Becker2002,nitsche_uber_2013}. The coefficient $\tau>0$ is a mesh-dependent parameter that has to be large enough to ensure the coercivity of $\mathrm{a}(\cdot,\cdot)$.
 We consider a weak imposition of boundary conditions since it is not straightforward to  include prescribed values in the approximation space in a strong manner, when dealing with an unfitted grid. Nitsche's method is commonly used in the embedded boundary community for circumventing this problem (see, e.g.,
\cite{burman_cutfem:_2015,Schillinger2015}) since it provides a consistent numerical scheme with optimal converge rates (also for high-order \acp{fe}).

When using the \ac{agfe} space (i.e., $\mathcal{V}^{\rm x}=\mathcal{V}^{\rm agg}$), we compute the mesh-dependent parameter $\tau$ as
\begin{equation}
\label{eq:tau-est}
\tau=\dfrac{\beta}{h},
\end{equation}
where $h$ is the cell size in the computational grid and $\beta$ is a user defined constant parameter. In the numerical examples, we take two different values of $\beta$, namely $\beta=10.0$ and $\beta=100.0$. The choice of $\beta=10.0$ is common in the literature of unfitted methods (see, e.g., reference \cite{Massing2014}). We also have considered $\beta=100.0$ in order to study the sensibility of the method with respect to this parameter. We note that using the space $\mathcal{V}^{\rm agg}$ and formula \eqref{eq:tau-est} for the computation of $\tau$, the bilinear form in \eqref{eq:weak-PoissonEq} is coercive independently of the position of the cuts. We refer to \cite{Badia2018} for the mathematical proof of this statement. However, when using the standard \ac{fe} space $\mathcal{V}^{\rm std}$, formula \eqref{eq:tau-est} does not necessarily lead to a coercive bi-linear form. For the simulations with the standard space $\mathcal{V}^{\rm std}$, we compute $\tau$ following an alternative approach in order recover a coercive problem independently of the position of the cuts (see \cite{DePrenter2017} for details). For the sake of completeness, we briefly present this methodology here. It consists in building the penalty parameter $\tau$ as an element-wise constant function, where the value in a generic cut cell, namely $\tau_k\doteq\tau|_{T^k}$, $k\in\mathcal{K}^{\rm cut}$, is computed as  $\tau_k = \beta \lambda^{\rm max}_k$. In this context, $\beta$ is also a user-defined constant (with values $\beta=10.0$ and $\beta=100.0$ in the numerical examples). On the other hand, $\lambda^{\rm max}_k$ is the maximum eigenvalue of the following generalized eigenvalue problem: find $\eta_k\in\mathcal{V}(T^k)$ and $\lambda_k\in\mathbb{R}$ such that
\begin{equation}
\int_{T^k\cap\Omega} \nabla\eta_k\cdot\nabla\xi_k \mathrm{\ d}\Omega = \lambda_k  \int_{T^k\cap\Gamma} (\nabla\eta_k\cdot n) (\nabla\xi_k\cdot n) \mathrm{\ d}\Gamma\text{ for all } \xi_k\in\mathcal{V}(T^k)\text{ for each cut cell id } k\in\mathcal{K}^{\rm cut} .
\end{equation}
Note that this problem is local at each cut cell, and its definition only depends on how the physical domain intersects the cut cell and the functional space used within the cell.  We refer to \cite{DePrenter2017} for the specific details on the solution of this problem. As indicated in \cite{DePrenter2017}, one of the main drawbacks of using the standard \ac{fe} space $\mathcal{V}^{\rm std}$ is that the optimal value of $\tau$ needed for ensuring coercivity of the bilinear form $\mathrm{a}(\cdot,\cdot)$ (i.e.,  $\lambda^{\rm max}_k$) might be arbitrarily large if the active part of a cell $T\in\mathcal{T}$, namely $\Omega\cap T$, has a volume that tends to zero. This fact leads to strongly ill-conditioned systems of linear equations when considering the standard space $\mathcal{V}^{\rm std}$. This issue is solved when considering the \ac{agfe} space $\mathcal{V}^{\rm agg}$ as the optimal value of $\tau$, which leads to a coercive bi-linear form, is independent of the position of the cuts (see \cite{Badia2018} for further details).

\subsection{Linear solver setup}\label{sec:ls-setup}

The \ac{fe} discretization of the weak problem represented by the bi-linear and linear forms in \eqref{eq:weak-PoissonEq} leads to a system of linear algebraic equations, namely $\mathbf{A}\mathbf{x}=\mathbf{b}$, which, at large scales, requires efficient parallel linear solvers.
Our goal is to show that the  usage of the \ac{agfe} spaces allows one to efficiently solve the linear systems using solvers that are well established for conventional \ac{fe} analysis based on body-fitted grids. To this end, we consider the widely available linear solvers available in the \petsc~library \cite{petsc-user-ref}. In particular, we use a conjugate gradient method from the \texttt{KSP} module of \petsc, preconditioned with a smoothed-aggregation \ac{amg} preconditioner called  \texttt{GAMG} (see \cite{May2016} for specific details).

We set up the linear solver by relying as much as possible on the default configuration given by \texttt{GAMG} to effectively show that the \ac{agfe} spaces lead to linear systems that are efficiently solved using standard multigrid tools. The modifications that we have introduced in the default parameters are, basically, in order to build a preconditioner tailored for symmetric definite positive matrices (since the default configuration of \texttt{GAMG} is designed for the non-symmetric case). We report the configuration file used to set up the \petsc~solver and preconditioner in Listing~\ref{lst:petscrc}. The file is written using a domain specific language provided by the \petsc~library (see the \petsc~users' manual \cite{petsc-user-ref} for further details on the file syntax).

\begin{lstlisting}[float,caption=Contents of the used \petsc{} configuration file.,label=lst:petscrc]
-ksp_type cg(*@\label{lst:petscrc:1}@*)
-ksp_rtol 1.0e-6(*@\label{lst:petscrc:2}@*)
-ksp_converged_reason
-ksp_max_it 500
-ksp_norm_type unpreconditioned (*@\label{lst:petscrc:5}@*)
-pc_type gamg (*@\label{lst:petscrc:6}@*)
-pc_gamg_type agg (*@\label{lst:petscrc:7}@*)
-mg_coarse_sub_pc_type cholesky (*@\label{lst:petscrc:8}@*)
-mg_levels_esteig_ksp_type cg (*@\label{lst:petscrc:9}@*)
-pc_gamg_square_graph 0 (*@\label{lst:petscrc:10}@*)
-build_twosided redscatter (*@\label{lst:petscrc:11}@*)
\end{lstlisting}

In line~\ref{lst:petscrc:1} of Listing~\ref{lst:petscrc}, we select a conjugate gradient method to compute an approximation $\mathbf{x}^{\rm cg}\approx\mathbf{x}$ of the solution of the symmetric definite positive system $\mathbf{A}\mathbf{x}=\mathbf{b}$. We declare convergence when $\| \mathbf{r} \|_2/ \| \mathbf{b} \|_2 < 10^{-6}$ within the first $500$ iterations (see lines ~\ref{lst:petscrc:2}-\ref{lst:petscrc:5}), where $\mathbf{r}\doteq \mathbf{b}-\mathbf{A}\mathbf{x}^{\rm cg}$ is the un-preconditioned residual and $\|\cdot\|_2$ denotes the Euclidean norm. In lines \ref{lst:petscrc:6}-\ref{lst:petscrc:7}, we select the smoothed aggregation \ac{amg} preconditioner called \texttt{GAMG} and in line \ref{lst:petscrc:8} we select a Cholesky factorization in order to solve the problem on the coarsest level of the \ac{amg} hierarchy. The \texttt{GAMG} preconditioner uses Chebyshev smoothers  in order to approximate the solution at the finest and  intermediate levels of the \ac{amg} hierarchy. In order to set up a Chebyshev smoother one needs to estimate the maximum eigenvalue of the corresponding system matrix (see, e.g., \cite{Hu2003} for details). To compute this eigenvalue estimate, \texttt{GAMG} uses internally an iterative Krylov sub-space method. In line \ref{lst:petscrc:9}, we instruct \texttt{GAMG} to use a conjugate gradient iteration for this purpose. {\rev In particular, we are using the default number of iterations provided by {\petsc} (i.e., 10 iterations) in the conjugate gradient method used to estimate the eigenvalue.} In the process of generating the prolongation and restriction operators associated with the \ac{amg} hierarchy, \texttt{GAMG} allows one to use the graph of the matrix $\mathbf{A}^T\mathbf{A}$ instead of the graph of the original matrix $\mathbf{A}$. By default, \texttt{GAMG} uses the graph of  $\mathbf{A}^T\mathbf{A}$ in the finest level and the graph of (the coarse representation of) $\mathbf{A}$ at the other levels. In line \ref{lst:petscrc:10}, we instruct \texttt{GAMG} to use the graph of $\mathbf{A}$ also in the finest level. We have observed that this setting leads to a faster solver (since the product $\mathbf{A}^T\mathbf{A}$ is not computed) even though the complexity of the \ac{amg} preconditioner is slightly increased. Finally, we have included the option in line \ref{lst:petscrc:11} in order to bypass a software BUG in the Intel MPI implementation {\rev used in the experiments}.
Please refer to the \texttt{GAMG} package documentation \cite{GAMGweb} for the default methods and parameters that are not covered in Listing~\ref{lst:petscrc}.

\subsection{Problem setup}\label{sec:problem-setup}

 In this section, we detail the rest of the parameter values and methods used in the numerical examples. In order to show that the methods are robust enough to deal with different geometrical data, we have performed all the experiments for {\rev three} different complex geometries. On the one hand, we have considered a bulky 3D body whose shape reminds the one of a popcorn flake (see Fig.~\ref{fig:geom-popconr}). This "popcorn-flake" geometry is often used in the literature to study the performance of unfitted FE methods (see, e.g., \cite{burman_cutfem:_2015}). The popcorn flake geometry considered here is obtained by taking the one defined in \cite{burman_cutfem:_2015}, scaling it by a factor of 0.5 and translating it {\rev by} a value of 0.5 in each direction such that the body fits in the unit cube $[0, 1]^3$. The second geometry is a massive spiral pipe (see Fig.~\ref{fig:geom-spiral}). The radius of its tubular cross section is 0.1, whereas the radius of the spiral central axis is 0.875. {\rev The third geometry, called "Swiss cheese", is a complex body often used is the literature of unfitted \ac{fe} methods. The precise definition of this geometry is found in \cite{burman_cutfem:_2015}.  In all cases}, the bounding box used to define the background mesh is the unit cube $[0, 1]^3$, which is also depicted in Fig.~\ref{fig:geom}.

\begin{figure}[ht!]
  \centering
  
  \begin{subfigure}{0.99\textwidth}
    \centering
    \includegraphics[width=0.45\textwidth]{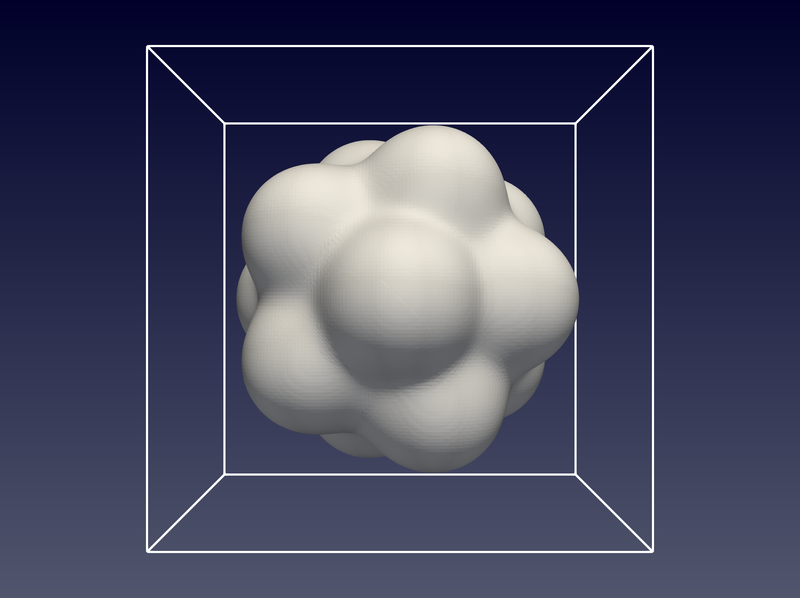}
    \hspace{1em}
    \includegraphics[width=0.45\textwidth]{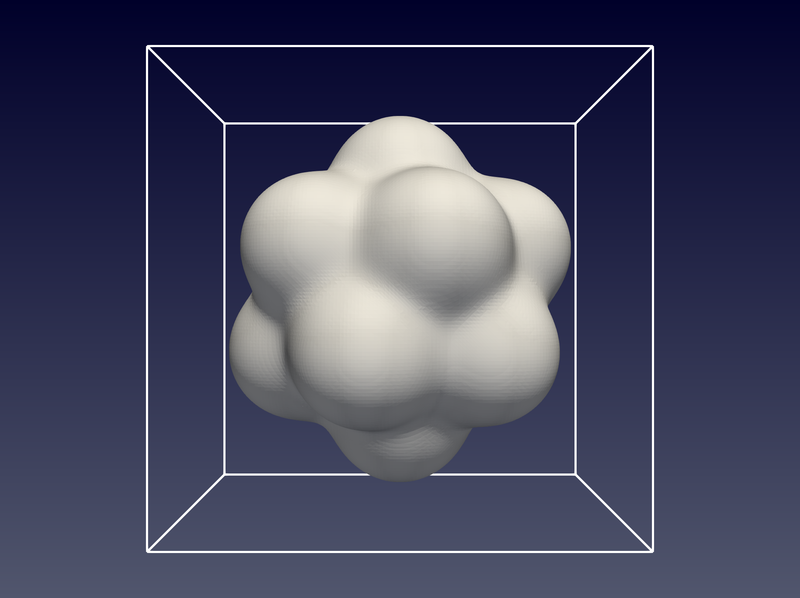}
    \caption{Popcorn flake.}
    \label{fig:geom-popconr}
  \end{subfigure}
  
  \vspace{1em}
  \begin{subfigure}{0.99\textwidth}
    \centering
    \includegraphics[width=0.45\textwidth]{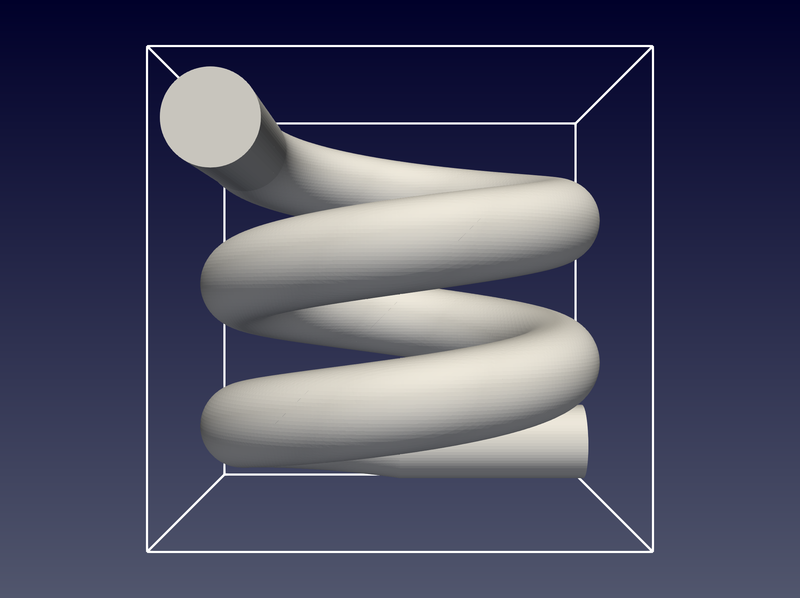}
    \hspace{1em}
    \includegraphics[width=0.45\textwidth]{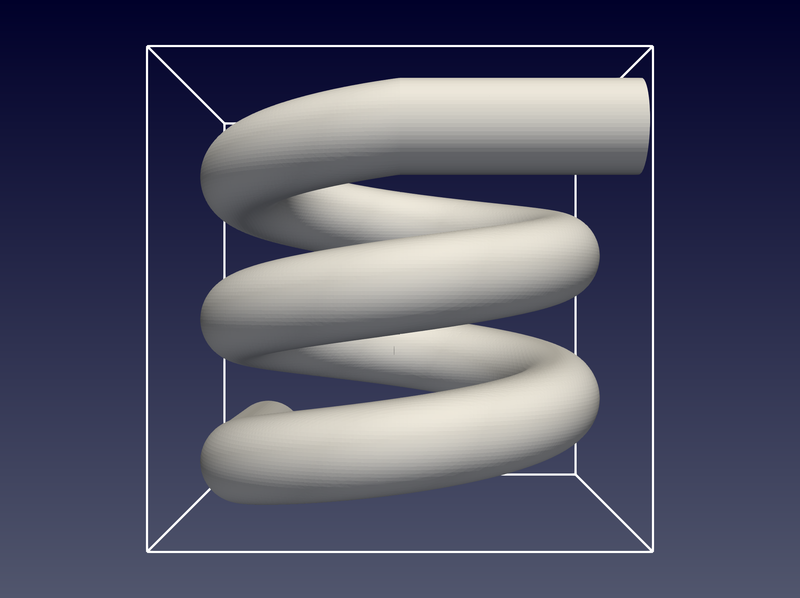}
    \caption{Spiral.}
    \label{fig:geom-spiral}
  \end{subfigure}
  
  \vspace{1em}
  \begin{subfigure}{0.99\textwidth}
    \centering
    \includegraphics[width=0.45\textwidth]{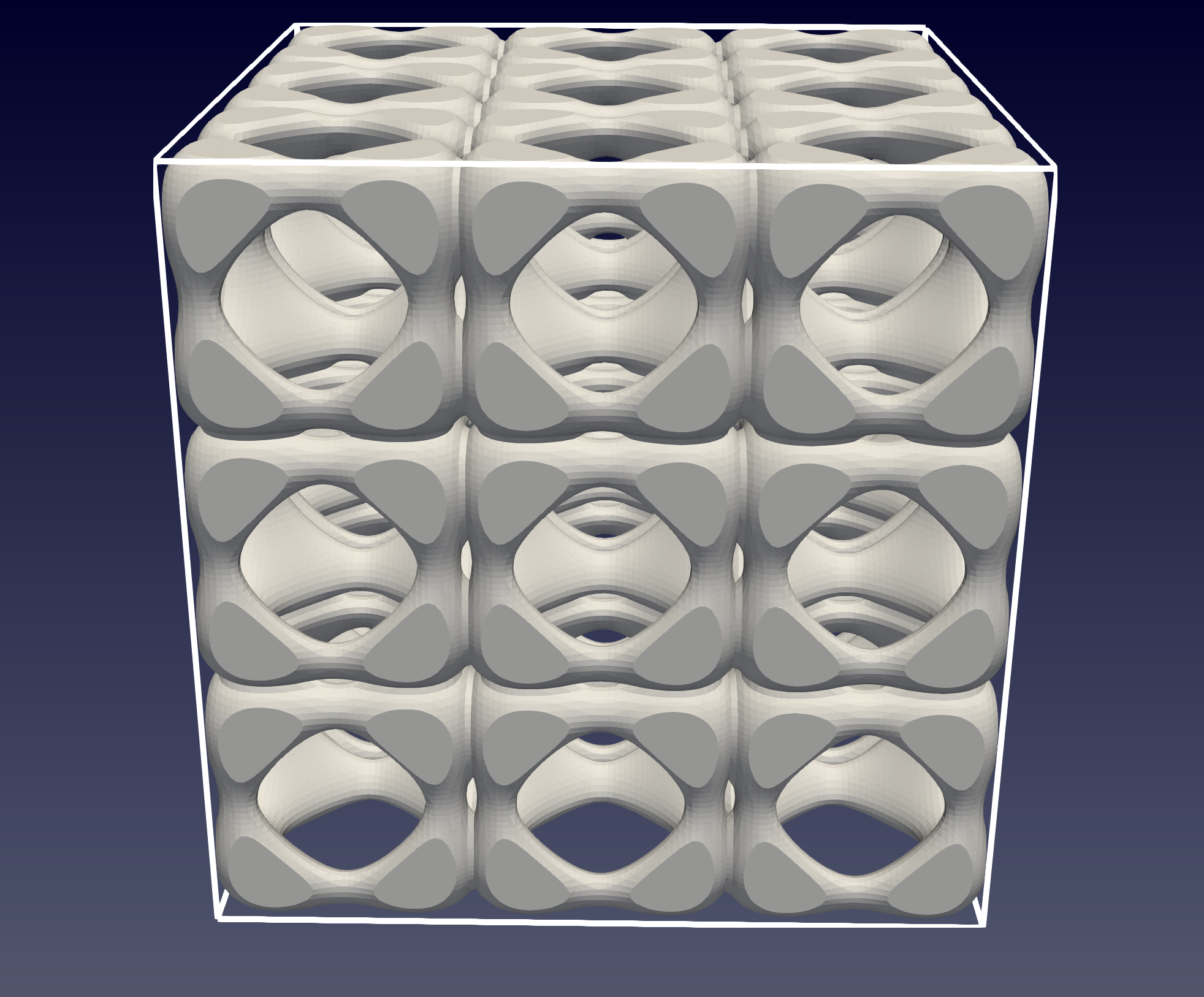}
    \hspace{1em}
    \includegraphics[width=0.45\textwidth]{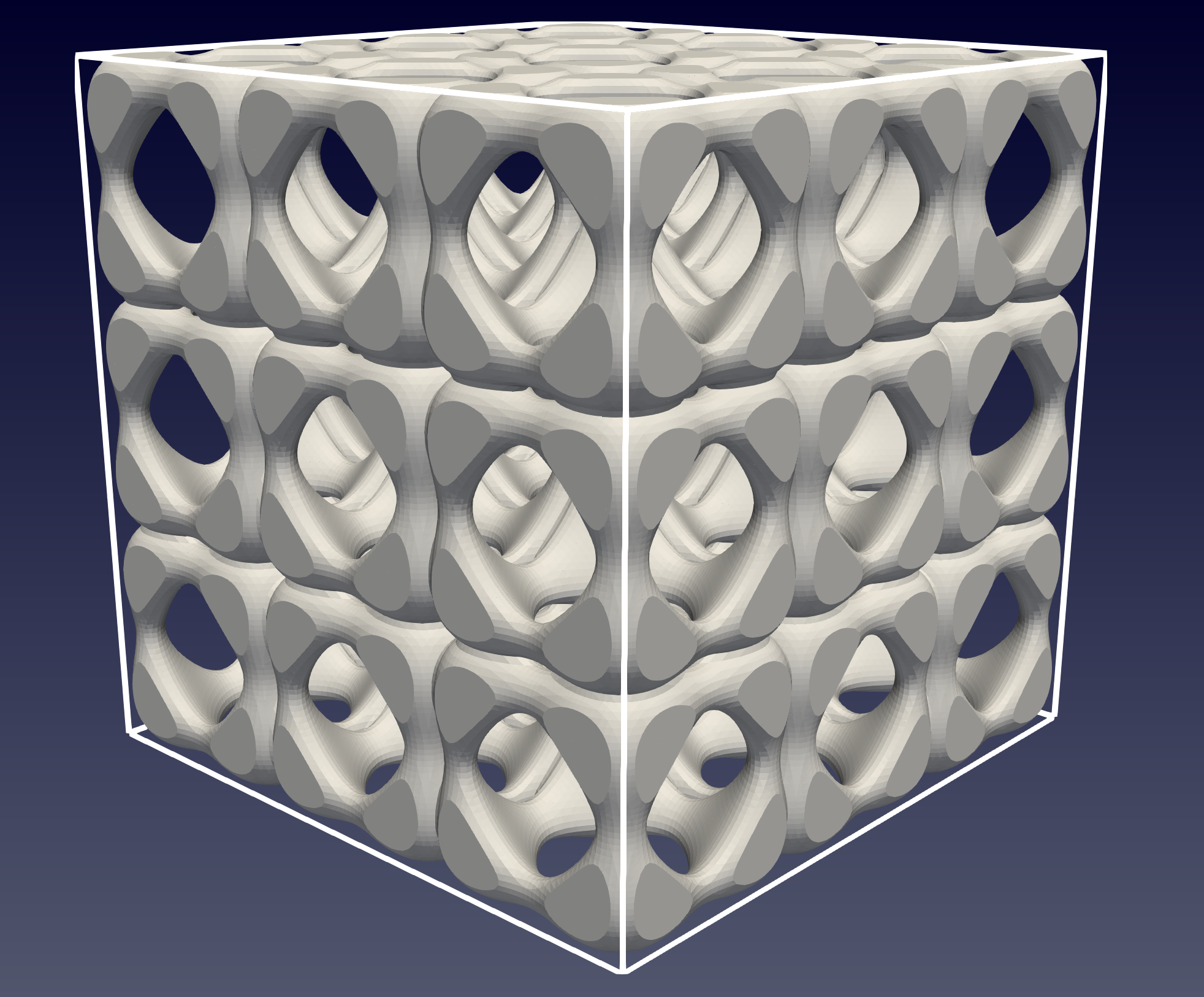}
    \caption{\rev Swiss Cheese.}
    \label{fig:geom-cheese-block}
  \end{subfigure}

\caption{Views of the {\rev three} geometries and the bounding boxes considered in the numerical examples.}
\label{fig:geom}
\end{figure}

The source term function $f$ and the Dirichlet condition function $g$ of the Poisson problem \eqref{eq:PoissonEq} are defined using the method of the manufactured analytic solutions so that problem \eqref{eq:PoissonEq} has the following exact solution: $u(x,y,z)=x+y+z$. This exact solution is used to check the quality of the computed \ac{fe} approximations. Since the exact solution $u$ belongs to the considered \ac{fe} spaces, the computed \ac{fe} approximations $u_h$ have to coincide with $u$ up to the tolerance used to declare convergence in the linear solver. In order to check that this is true in our experiments, we have computed the (relative) error norms $|u-u_h|_{H^1}/ | u |_{H^1}$ and $\| u-u_h \| _{L^2}/ \ \| u \|_{L^2}$, where
\begin{equation}
| u |_{H^1} \doteq \left( \int_{\Omega} \nabla u\cdot \nabla u \mathrm{\ d}\Omega \right)^{1/2}\text{ and } \| u \|_{L^2}\doteq\left( \int_{\Omega} u^2 \mathrm{\ d}\Omega  \right)^{1/2}
\end{equation}
are the $H^1$ semi-norm and $L^2$ norm, resp.

As indicated in Sect.~\ref{sec:imb_setup}, we consider background meshes of hexahedral cells generated and partitioned in parallel with the \p4est library. We generate the meshes by recursive sub-division of the unit cube $[0,1]^3$. That is, in the first refinement level,  the cube is decomposed into $(2^1)^3=8$ cells, in the second refinement level into $(2^2)^3=64$ cells, etc. For the weak scaling test shown below, we generate three families of partitions such that the average number of cells per sub-domain is constant as the mesh is refined (see Table~\ref{tab:partitions}). The first one (referred as load 1) has approximately $1,000$ cells per sub-domain, and the second and the third ones (load 2 and load 3) have approximately $8,000$ and $64,000$ cells per sub-domain, resp. {\rev The Swiss cheese geometry is not simulated with the coarsest mesh since it is unable to capture the main geometrical features.} The finest mesh considered in the examples has $1,073,741,824$ cells (corresponding to refinement level 10), and it is partitioned into $16,777$ sub-domains (which are mapped to the same {\rev number} of MPI tasks). {\rev We used the default MPI task placement policy of Intel MPI (v2018.4.057) with partially filled nodes. 
The number of nodes $N$ is selected as $N=\texttt{ceil}(|\mathcal{S}|/48)$, where $\texttt{ceil}(x)$ rounds $x$ to the nearest integer that is greater or equal to $x$. 
If $|\mathcal{S}|$ is not multiple of $48$, the placement policy fully populates the first $N-1$ nodes with $48$ MPI tasks per node.
The remaining $\texttt{mod}(|\mathcal{S}|,48)$ MPI tasks are placed on the last node, with $\texttt{mod}$
denoting integer division remainder.}

\begin{table}[ht!]
\centering
\begin{small}
\begin{tabular}{rrrrr}
\toprule
 & \multicolumn{3}{c}{Number of cells}\\
\cmidrule{2-4}
 {\rev $|\mathcal{S}|$} & load 1 & load 2& load 3 \\
\midrule
4 &  4,096 {\rev (34\% 40\% \hspace{0.8em}-\hspace{0.8em})}  & 32,768 {\rev (29\% 29\% 56\%)} & 262,144 {\rev (26\% 24\% 47\%)}\\
32 & 32,768 {\rev (29\% 29\% 56\%)}   & 262,144 {\rev (26\% 24\% 47\%)} & 2,097,152 {\rev (25\% 22\% 41\%)}\\
262 & 262,144 {\rev (26\% 24\% 47\%)} & 2,097,152 {\rev (25\% 22\% 41\%)} & 16,777,216 {\rev (24\% 20\% 38\%)}\\
2,097 & 2,097,152 {\rev (25\% 22\% 41\%)} & 16,777,216 {\rev (24\% 20\% 38\%)} & 134,217,728 {\rev (24\% 20\% 36\%)}\\
16,777 & 16,777,216 {\rev (24\% 20\% 38\%)} & 134,217,728 {\rev (24\% 20\% 36\%)} & 1,073,741,824 {\rev (24\% 20\% 35\%)}\\
\bottomrule
\end{tabular}
\end{small}
\caption{Number of sub-domains {\rev ($|\mathcal{S}|$)} and number of total cells in the background mesh for the three families of partitioned meshes considered in the examples. {\rev In parentheses, the percentage of active cells (over the total number of background cells) for each of the three geometries is given. Within a parenthesis, the three given quantities correspond to the "pop-corn flake", the spiral, and the "Swiss cheese" geometries, respectively.}}
\label{tab:partitions}
\end{table}

Each mesh is decomposed using the load-balancing mechanism of the \p4est library. By default, the method splits the mesh such that the number of cells in each of the resulting sub-domains is approximately equally balanced. \p4est also includes the option to provide a user-defined weight to each of the cells. In this case, the algorithm leads to partitions, where the sum of the weights in each sub-domain is equally balanced across the entire partition. We use this functionality in order to give different weights to active and exterior cells. In this context, one has to find a trade-off between two conflicting goals. On the one hand, it is desirable to load balance the storage of the background mesh and the cost of the geometrical computations that run on top of it. In this case, one has to assign the same weight to active and exterior cells. On the other hand, it is also desirable to load balance the assembly, storage, and solution of the underlying systems of linear algebraic equations. In that case, exterior cells do not play any role, and one has to balance the number of active cells per sub-domain. This is achieved in practice by assigning a much smaller weight to exterior cells than to active ones. In the experiments,
we have used a ratio {\rev $w=10/1$ (i.e., weight  $10$ for active cells and $1$ for exterior ones)}. We observed experimentally {\rev (see Figs. \ref{fig:time-w-phases} and \ref{fig:time-w-s-phases})} that this value
of $w$ leads to the best trade-off among the two aforementioned factors {\rev for the considered geometries. The wall clock time of the main phases of the AgFEM method (Fig. \ref{fig:time-w-phases}) and the wall clock time of the linear solver (Fig. \ref{fig:time-w-s-phases}) decreases significantly when the weight factor $w$ moves from $w=1$ to $w=10$. Beyond this point, the timings tend to stagnate or slightly increase for $w\geq 10$. An important exception is the assembly of the algebraic linear system, whose wall clock time significantly increases again for $w\geq 10$ (see the orange curve with pentagonal markers in Fig. \ref{fig:time-w-phases}). For this reason, we have chosen a value of $w=10$ for the subsequent numerical experiments.}

\begin{figure}[ht!]
  \centering
   
  \begin{subfigure}{0.45\textwidth}
    \includegraphics[scale=1]{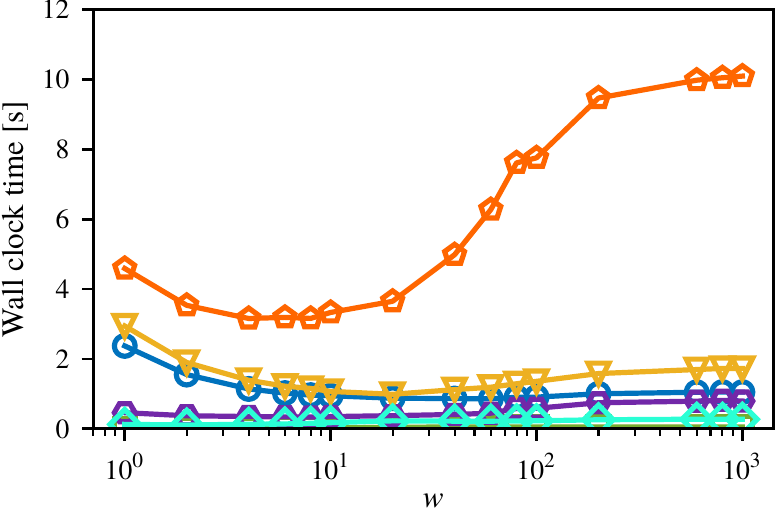}
        \caption{Popcorn flake.}
    \label{fig:time-phases-w-a}
  \end{subfigure}
  \hspace{2em}
  \begin{subfigure}{0.45\textwidth}
    \includegraphics[scale=1]{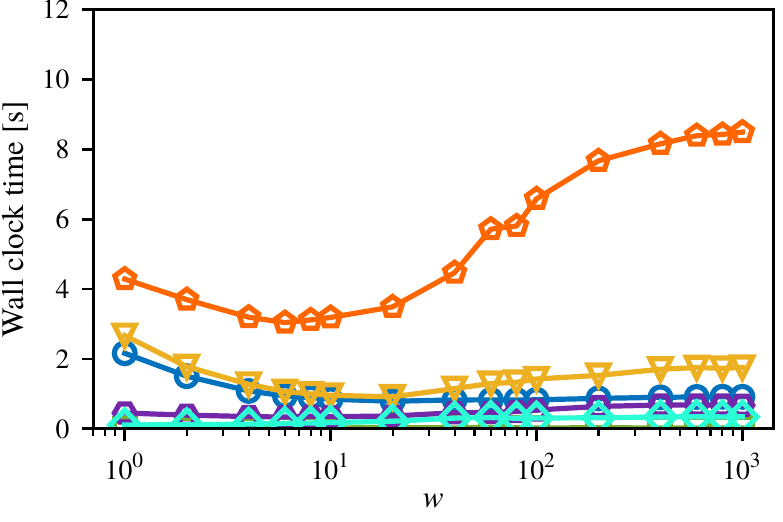}
    \caption{Spiral.}
    \label{fig:time-phases-w-b}
  \end{subfigure}
  
    \begin{subfigure}{0.45\textwidth}
    \includegraphics[scale=1]{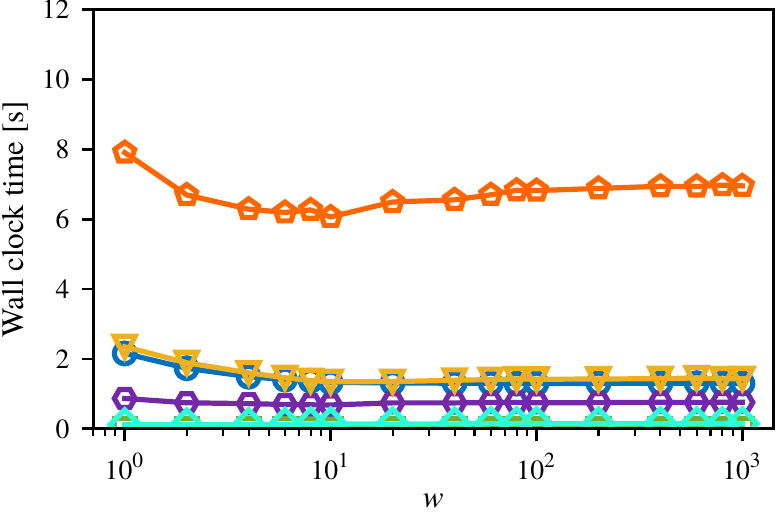}
    \caption{{Swiss cheese.}}
    \label{fig:time-phases-w-c}
  \end{subfigure}
      \hspace{2em}
      \begin{subfigure}{0.45\textwidth}
    \includegraphics[scale=1]{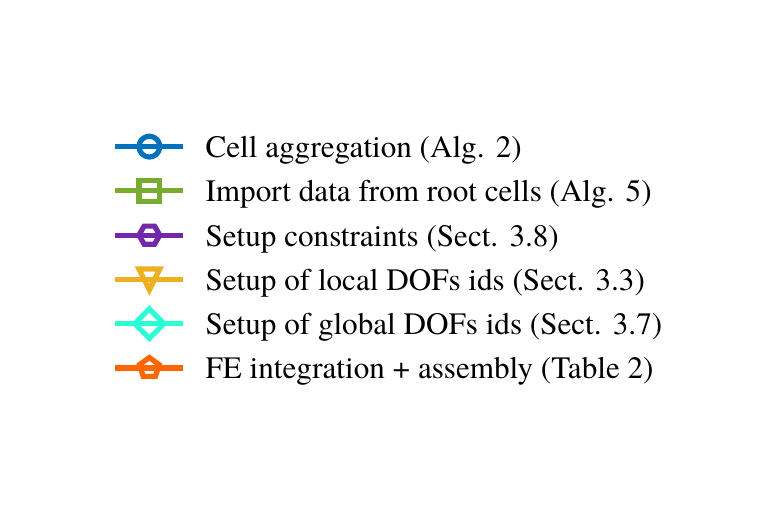}
  \end{subfigure}  
\caption{ \rev Wall clock time for the major phases of the \ac{agg} method versus the weighting factor $w$. Results computed with \ac{agfe} spaces, $\beta=10$, and load~3 on 262 CPU cores. }
\label{fig:time-w-phases}
\end{figure}

\begin{figure}[ht!]
  \centering
   
  \begin{subfigure}{0.45\textwidth}
    \includegraphics[scale=1]{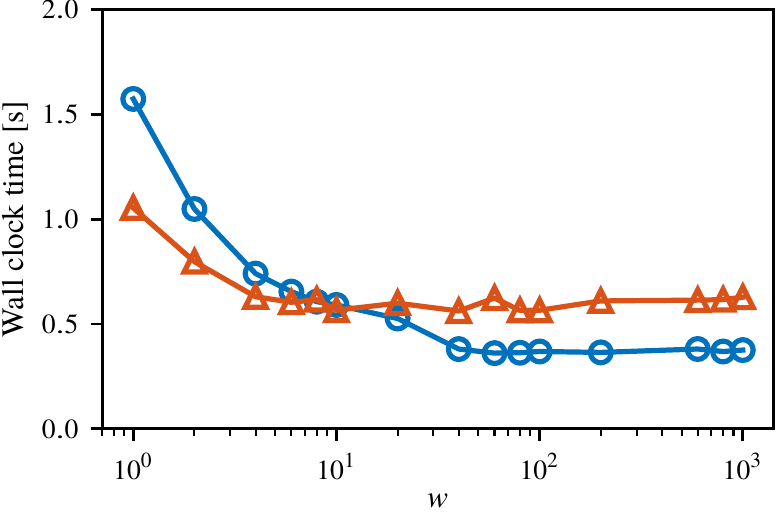}
        \caption{Popcorn flake.}
    \label{fig:time-phases-w-s-a}
  \end{subfigure}
  \hspace{2em}
  \begin{subfigure}{0.45\textwidth}
    \includegraphics[scale=1]{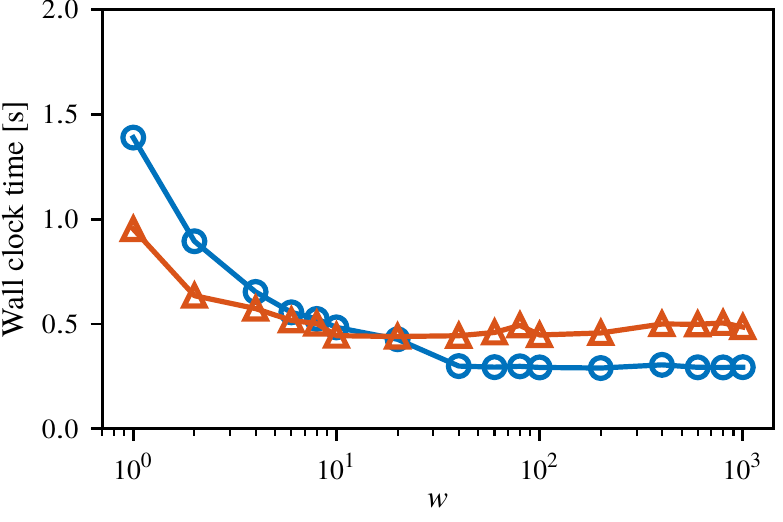}
    \caption{Spiral.}
    \label{fig:time-phases-w-s-b}
  \end{subfigure}
  
    \begin{subfigure}{0.45\textwidth}
    \includegraphics[scale=1]{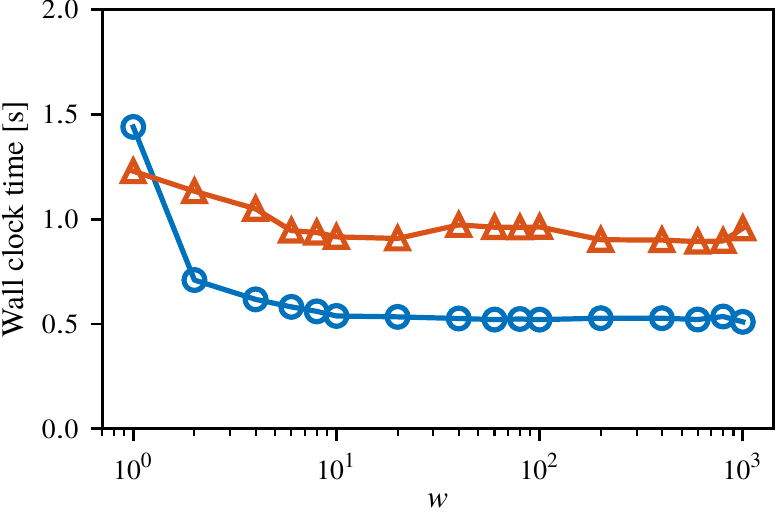}
    \caption{Swiss cheese.}
    \label{fig:time-phases-w-s-c}
  \end{subfigure}
      \hspace{2em}
      \begin{subfigure}{0.45\textwidth}
    \includegraphics[scale=1]{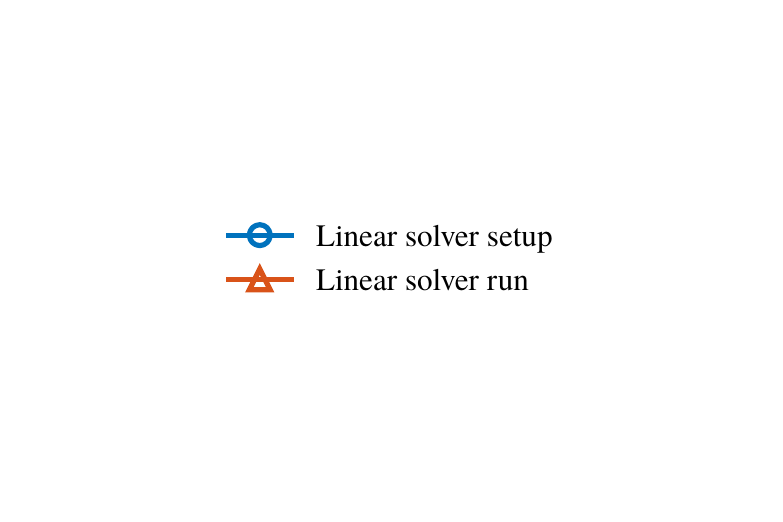}
  \end{subfigure}  
\caption{\rev Wall clock time for the main linear solver phases versus the weighting factor $w$. Results computed with \ac{agfe} spaces, $\beta=10$,  and load~3 on 262 CPU cores. }
\label{fig:time-w-s-phases}
\end{figure}

{\rev We  build the \ac{fe} spaces $\mathcal{V}^{\rm std}$ and $\mathcal{V}^{\rm agg}$ using} hexahedral elements with continuous piece-wise trilinear shape functions.  We consider both \ac{agfe} and conventional spaces in order to evaluate the benefits of using the \ac{agg} method. A summary of the main parameters and methods considered in the numerical examples is given in Table~\ref{tab:params}.

\begin{table}[ht!]
\begin{small}
\begin{tabular}{lll}
\toprule
Description &  Considered methods/values & \\
\midrule
Model problem & Poisson equation\\
Problem geometry & "Popcorn flake", massive spiral{\rev , and "Swiss cheese"}\\
Parallel mesh generation and partitioning & \p4est library\\
\ac{fe} spaces & Aggregated $\mathcal{V}^{\rm agg}$ and standard $\mathcal{V}^{\rm std}$\\
{\rev Cell type} & {\rev Hexahedral cells} \\
{\rev Interpolation} & {\rev Piece-wise trilinear shape functions} \\
Linear solver & Preconditioned conjugate gradients (\petsc)\\
Preconditioner & Smoothed-aggregation \ac{amg} (\texttt{GAMG} preconditioner from \petsc)\\
Weighting factor for space filling-curve & {\rev $w=10$} \\
Coefficient used in Nitsche's penalty terms  & $\beta=10.0$ and $\beta=100.0$ & \\

\bottomrule
\end{tabular}
\end{small}
\caption{Summary of the main parameters and computational strategies used in the numerical examples.}
\label{tab:params}
\end{table}

\newcommand{\secttitle}{Weak scaling test}
\subsection{\secttitle}\label{sec:weak-scaling-test}

We start this section devoted to the weak scaling test by checking that the problems were accurately solved in the experiments. To this end, Fig.~\ref{fig:errors-l2} shows the relative $L^2$ norm of the error associated with the computed numerical solutions versus problem size for the three loads per processor detailed in Table~\ref{tab:partitions}. {\rev The results in Figs.~\ref{fig:errors-l2-a},~\ref{fig:errors-l2-b}, and~\ref{fig:errors-l2-c} are for the "popcorn flake", the spiral and the "Swiss cheese" geometries.  Note that the results are qualitatively equivalent for the three geometries}.  {\rev In Fig.~\ref{fig:errors-l2}},  some points of the red lines (the ones obtained with standard \ac{fe} spaces) are missing since the linear solver was not able to provide a solution. Either the setup of the preconditioner failed or the computed preconditioner was not positive definite (which makes it unusable in a conjugate gradients iteration). These problems are related with the strong ill-conditioning that results from the usage of standard \ac{fe} spaces. Note, however,  that the linear solvers properly worked in all cases, when considering the \ac{agfe} spaces. Note that (as expected) the relative errors of the computed numerical solution are of the same order of magnitude as the tolerance imposed to declare convergence of the linear solver. {\rev Moreover, the computed errors in the $L^2$ norm are reduced about three orders of magnitude, if the problem is solved by setting a linear solver tolerance three orders of magnitude smaller, i.e., $10^{-9}$ instead of $10^{-6}$ (see Fig. \ref{fig:errors-l2-rtol}). This clearly demonstrates that the solution quality is dictated by the linear solver tolerance as one would expect.} Results for the $H^1$ norm of the error are qualitatively equivalent as the ones showed in Fig.~\ref{fig:errors-l2} and they are not included for the sake of brevity.


\begin{figure}[ht!]
  \centering
   
  \begin{subfigure}{0.45\textwidth}
    \includegraphics[scale=1]{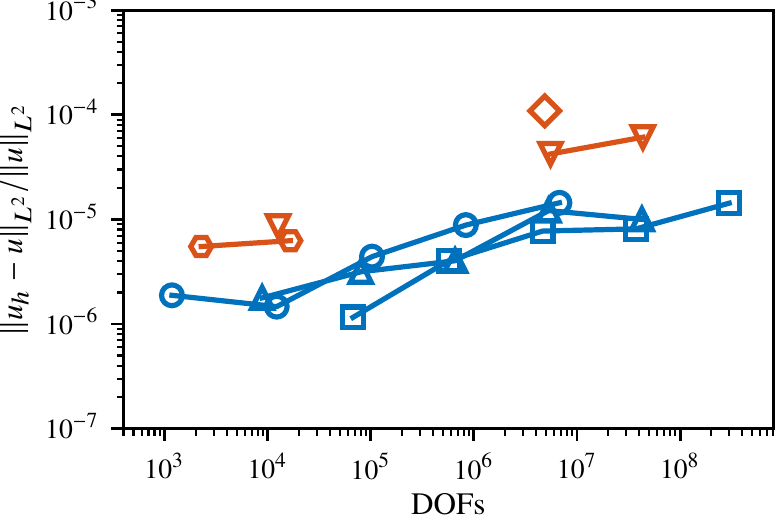}
        \caption{Popcorn flake.}
    \label{fig:errors-l2-a}
  \end{subfigure}
  \hspace{2em}
  \begin{subfigure}{0.45\textwidth}
    \includegraphics[scale=1]{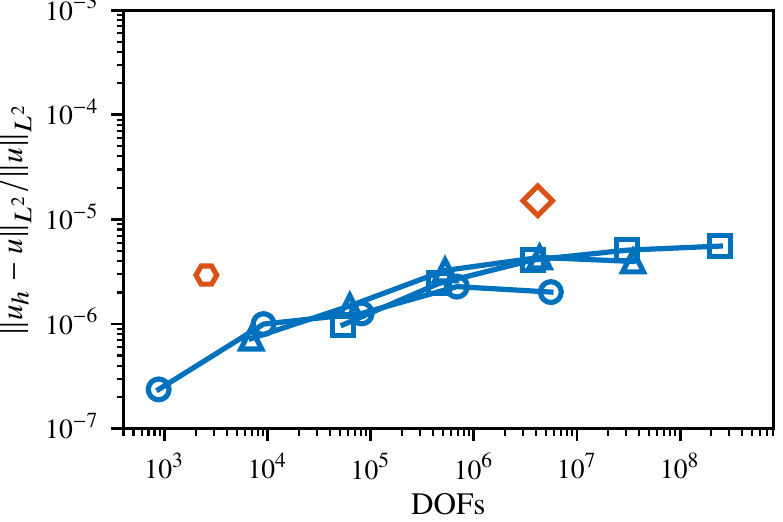}
    \caption{Spiral.}
    \label{fig:errors-l2-b}
  \end{subfigure}
  
  \begin{subfigure}{0.45\textwidth}
    \includegraphics[scale=1]{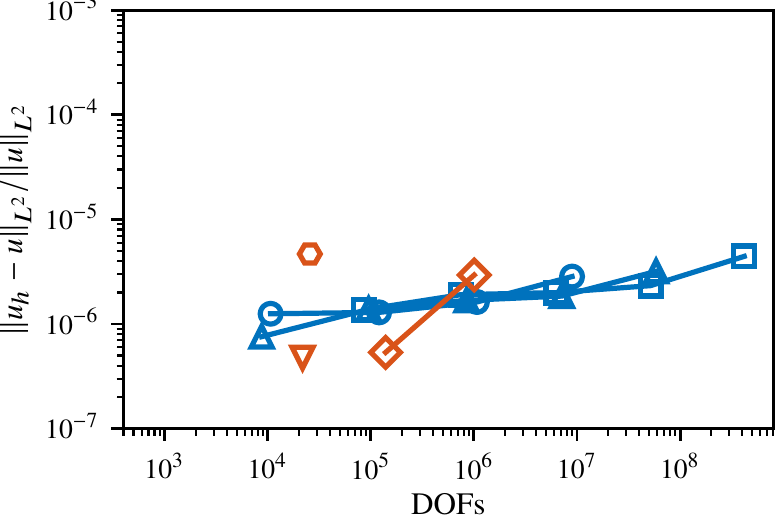}
    \caption{{\rev Swiss cheese.}}
    \label{fig:errors-l2-c}
  \end{subfigure}
    \hspace{2em}
      \begin{subfigure}{0.45\textwidth}
    \includegraphics[scale=1]{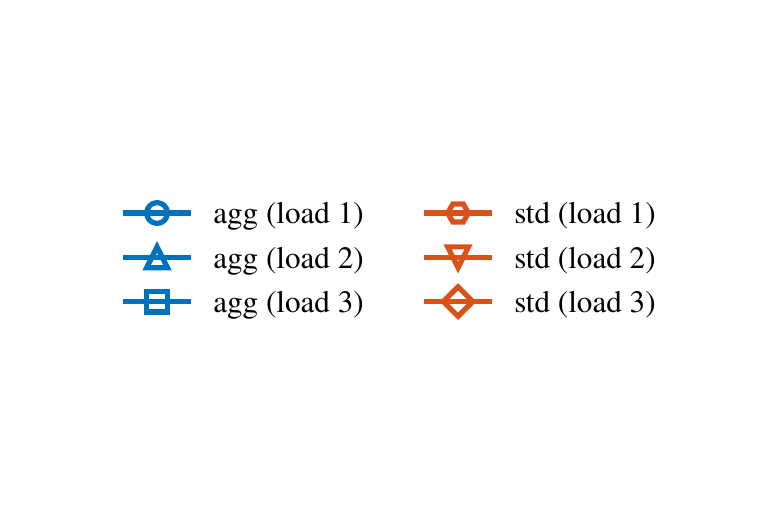}
  \end{subfigure}
\caption{\secttitle: Quality check of the computed numerical solution. The figure shows the (relative) $L^2$ norm of the error versus problem size for three different loads per processor. Results computed with parameters  $\beta=10$ and {\rev $w=10$}. }
\label{fig:errors-l2}
\end{figure}

\begin{figure}[ht!]
  \centering
   
  \begin{subfigure}{0.45\textwidth}
    \includegraphics[scale=1]{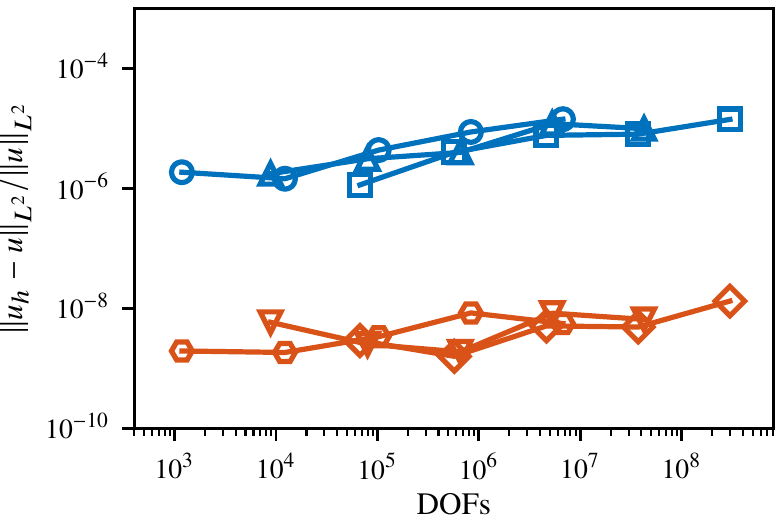}
        \caption{Popcorn flake.}
    \label{fig:errors-l2-rtol-a}
  \end{subfigure}
  \hspace{2em}
  \begin{subfigure}{0.45\textwidth}
    \includegraphics[scale=1]{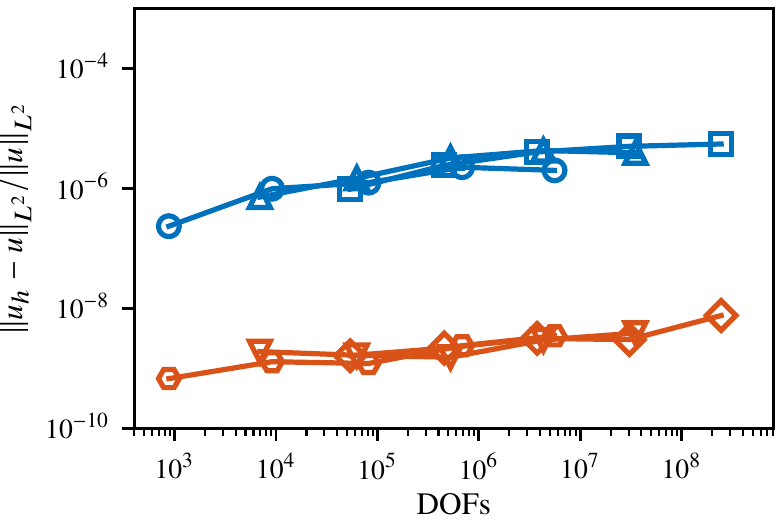}
    \caption{Spiral.}
    \label{fig:errors-l2-rtol-b}
  \end{subfigure}
  
  \begin{subfigure}{0.45\textwidth}
    \includegraphics[scale=1]{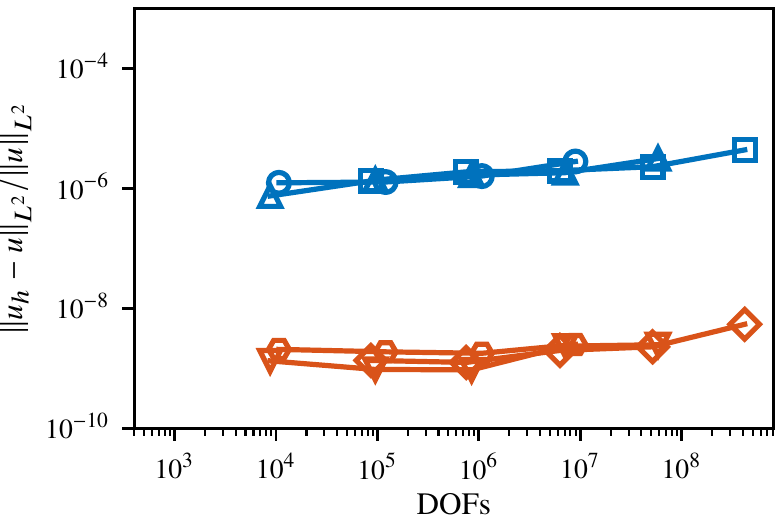}
    \caption{{\rev Swiss cheese.}}
    \label{fig:errors-l2-rtol-c}
  \end{subfigure}
    \hspace{2em}
      \begin{subfigure}{0.45\textwidth}
    \includegraphics[scale=1]{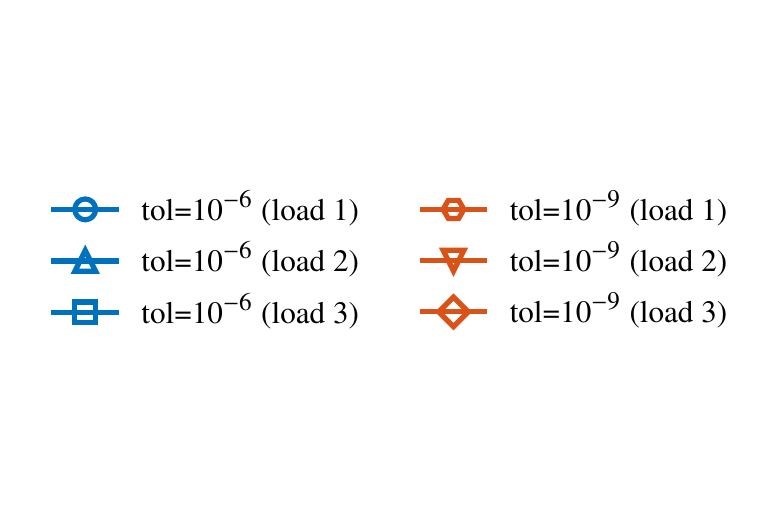}
  \end{subfigure}
\caption{ {\rev \secttitle: Quality check of the computed numerical solution. The figure shows the (relative) $L^2$ norm of the error versus problem size for three different loads per processor and for two different values of the linear solver tolerance. Results computed with aggregated \ac{fe} spaces and parameters  $\beta=10$ and {\rev $w=10$}.} }
\label{fig:errors-l2-rtol}
\end{figure}

We start the reporting of the weak scaling test, by discussing the impact of the different parameters on the algorithmic performance of the considered \ac{amg} preconditioners.
First, we consider the impact of using \ac{agfe} or standard \ac{fe} spaces in the number of solver iterations. As the results in Fig. \ref{fig:cg-iter-agg-vs-std} show, it is clear that the usage of \ac{agg}  is beneficial (and essential) for the good performance of the linear solver. As it was already discussed for Fig.~\ref{fig:errors-l2}, some results in Fig.~\ref{fig:cg-iter-agg-vs-std} (the ones associated with the standard \ac{fe} spaces) are missing since the linear solver was not able to provide a converged solution. This shows that standard \ac{fe} spaces are unusable in practice since they do not lead to a robust discretization method. In contrast, the usage of \ac{agfe} spaces allows one to effectively solve the underlying linear systems. In particular, note that the results obtained with \ac{agfe} spaces (blue lines) are very close to the expected optimal performance of multigrid methods (i.e., number of linear solver iterations is asymptotically independent of the problem size). {\rev This optimal behavior is observed  for  all} considered geometries, which further demonstrates the usability and robustness on the computational strategy based on \ac{agfe} spaces.

\begin{figure}[ht!]
  \centering
  \begin{subfigure}{0.45\textwidth}
    \includegraphics[scale=1]{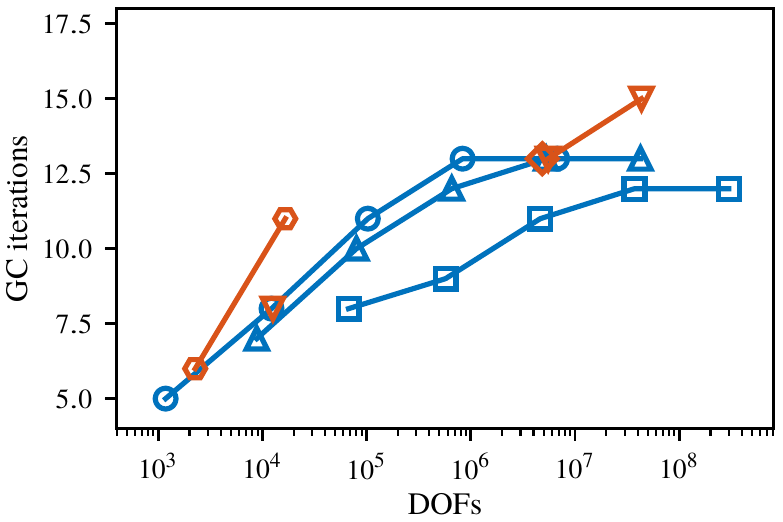}
        \caption{Popcorn flake.}
    \label{fig:cg-iter-agg-vs-std-a}
  \end{subfigure}
  \hspace{2em}
  \begin{subfigure}{0.45\textwidth}
    \includegraphics[scale=1]{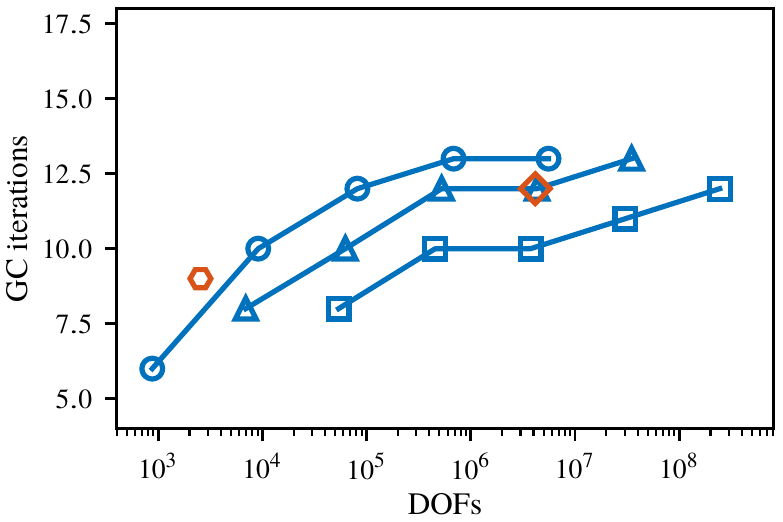}
    \caption{Spiral.}
    \label{fig:cg-iter-agg-vs-std-b}
  \end{subfigure}
  
  \begin{subfigure}{0.45\textwidth}
    \includegraphics[scale=1]{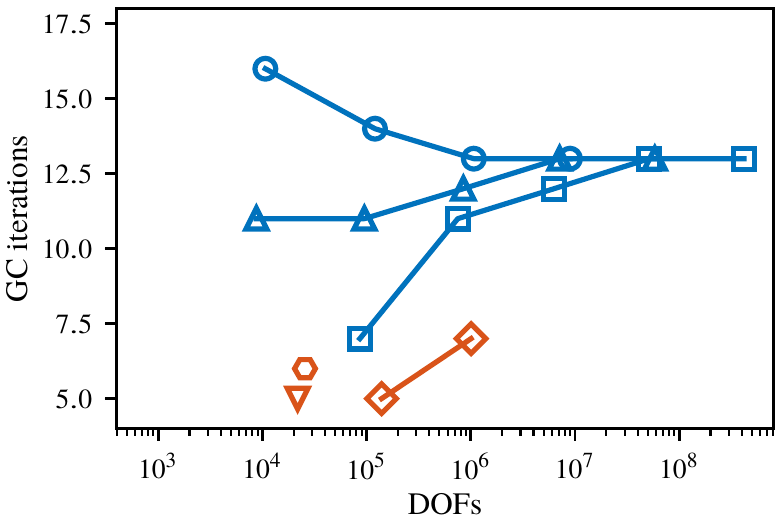}
    \caption{{\rev Swiss cheese.}}
    \label{fig:cg-iter-agg-vs-std-c}
  \end{subfigure}
    \hspace{2em}
      \begin{subfigure}{0.45\textwidth}
    \includegraphics[scale=1]{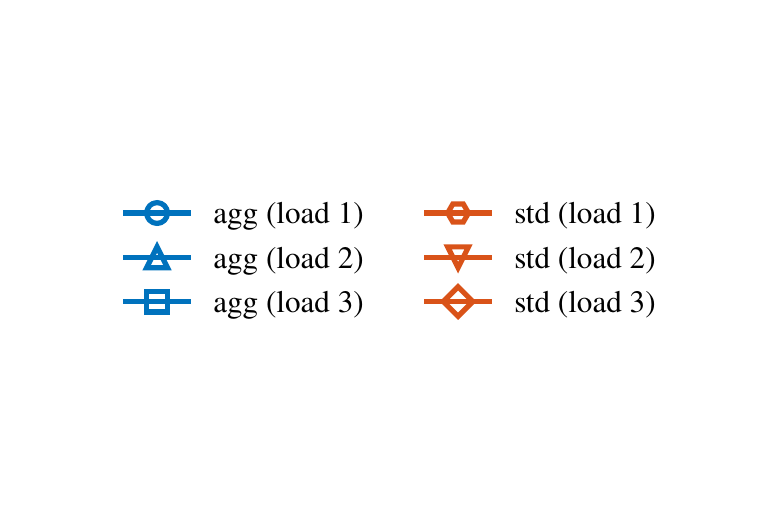}
  \end{subfigure}
\caption{\secttitle: Influence of standard or \ac{agfe} spaces on the linear solver behavior. The figure shows the conjugate gradient (CG) iterations versus problem size for three different loads per processor. Results computed parameters  $\beta=10$ and {\rev $w=10$}. }
\label{fig:cg-iter-agg-vs-std}
\end{figure}


Next, we study the influence the user-defined coefficient $\beta$ associated with Nitsche's method (see Sect.~\ref{sec:model-problem}) has in the performance of the solver. Since it is clear from Figs.~\ref{fig:errors-l2} and \ref{fig:cg-iter-agg-vs-std} that the method based on standard \ac{fe} spaces is unusable in practice, from now on,  we only provide results associated with the \ac{agfe} spaces. The results in Fig.~\ref{fig:beta-influence} show the influence of $\beta$ in the solver performance for the two values of $\beta$ previously detailed in Table~\ref{tab:params} and for the three loads per processor detailed in Table~\ref{tab:partitions}. 
 As it is seen in the figure, the number of linear solver iterations are weakly scalable for the two values of $\beta$ (i.e., iterations asymptotically independent of the problem size). However, the absolute number of iterations is significantly larger for $\beta=100$ than for $\beta=10$, which suggests that $\beta$ has to be chosen large enough to ensure coercivity of the underlying bi-linear form, but cannot be chosen arbitrary large since it has a negative impact on the solver performance. This requirement is well known in the literature of Nitsche-based \ac{fe} methods (see, e.g., \cite{burman_cutfem:_2015}).

\begin{figure}[ht!]
  \centering
  \begin{subfigure}{0.45\textwidth}
    \includegraphics[scale=1]{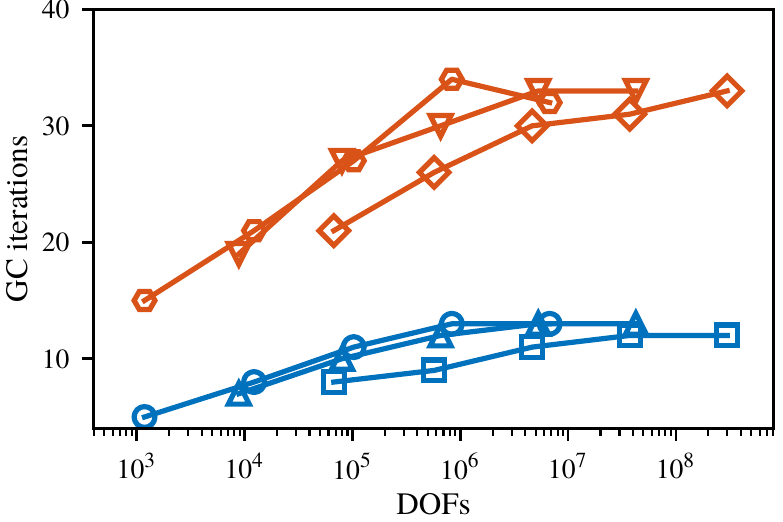}
        \caption{Popcorn flake.}
    \label{fig:beta-influence-a}
  \end{subfigure}
  \hspace{2em}
  \begin{subfigure}{0.45\textwidth}
    \includegraphics[scale=1]{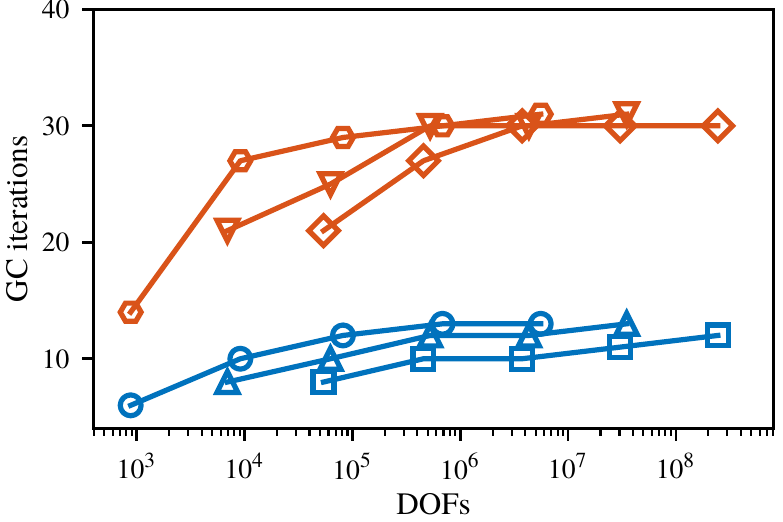}
    \caption{Spiral.}
    \label{fig:beta-influence-b}
  \end{subfigure}
  
  \begin{subfigure}{0.45\textwidth}
    \includegraphics[scale=1]{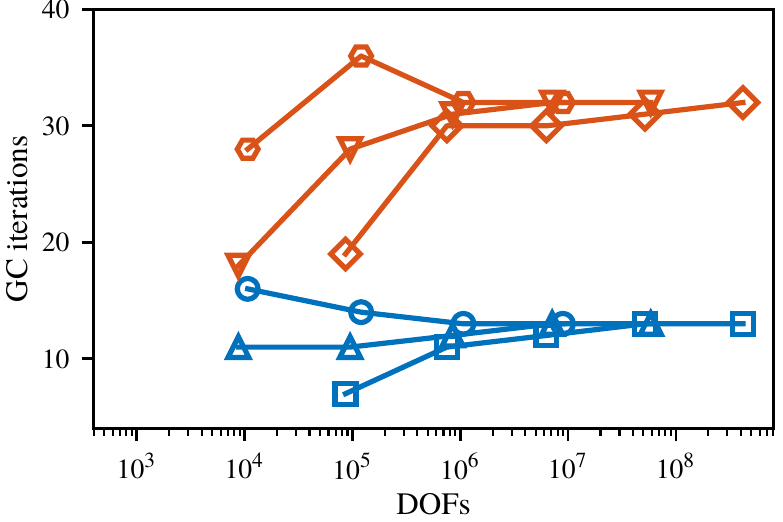}
    \caption{{\rev Swiss cheese.}}
    \label{fig:beta-influence-c}
  \end{subfigure}
    \hspace{2em}
      \begin{subfigure}{0.45\textwidth}
    \includegraphics[scale=1]{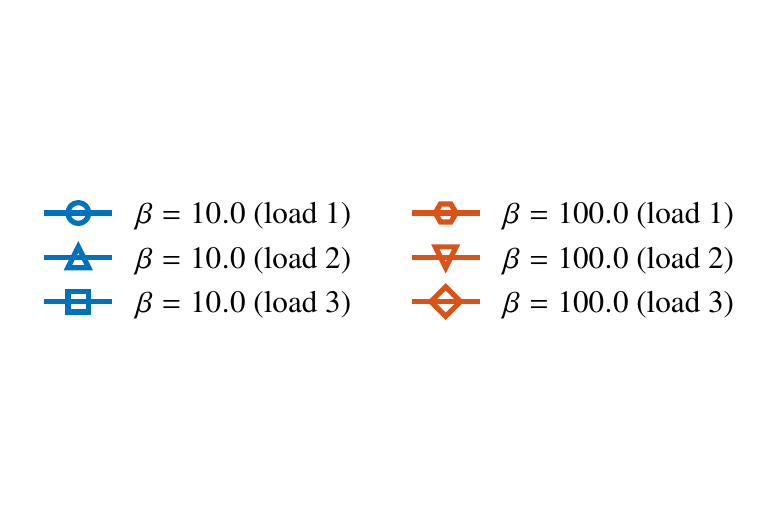}
  \end{subfigure}  
  
\caption{\secttitle: Influence of the coefficient $\beta$ on the linear solver behavior. The figure shows the conjugate gradient (CG) iterations versus problem size for three different loads per processor. Results computed with \ac{agfe} spaces and {$w=10$}.}
\label{fig:beta-influence}
\end{figure}

\begin{figure}[ht!]
  \centering
   
  \begin{subfigure}{0.45\textwidth}
    \includegraphics[scale=1]{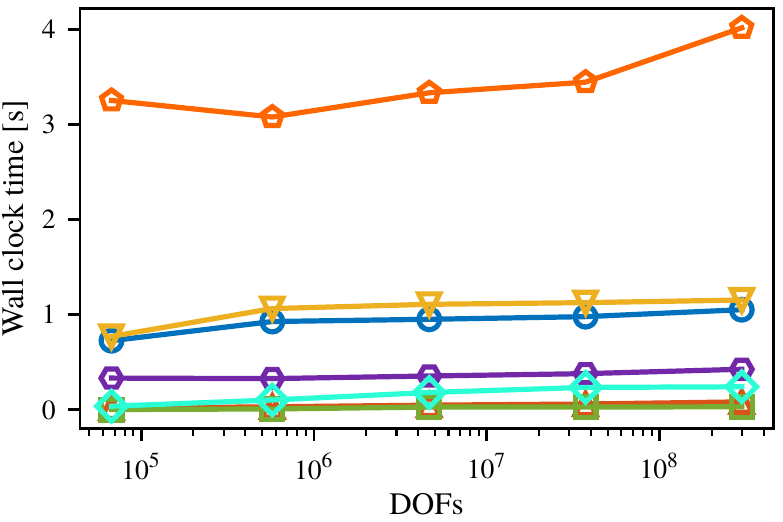}
        \caption{Popcorn flake.}
    \label{fig:time-phases-a}
  \end{subfigure}
  \hspace{2em}
  \begin{subfigure}{0.45\textwidth}
    \includegraphics[scale=1]{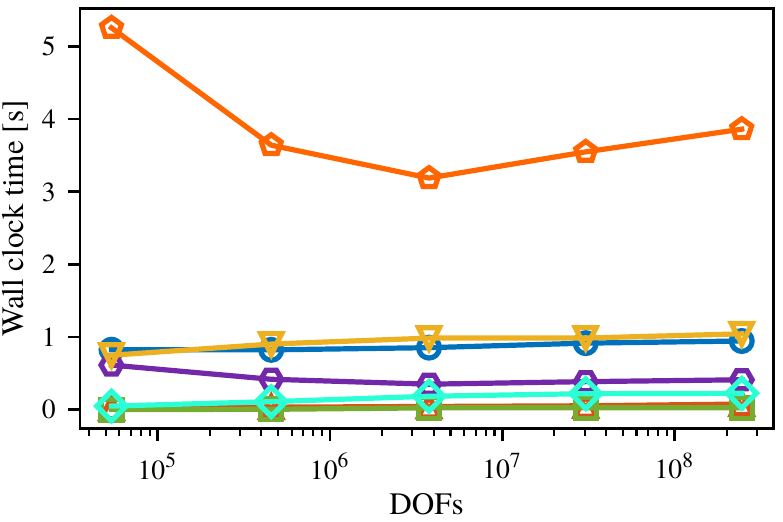}
    \caption{Spiral.}
    \label{fig:time-phases-b}
  \end{subfigure}
  
    \begin{subfigure}{0.45\textwidth}
    \includegraphics[scale=1]{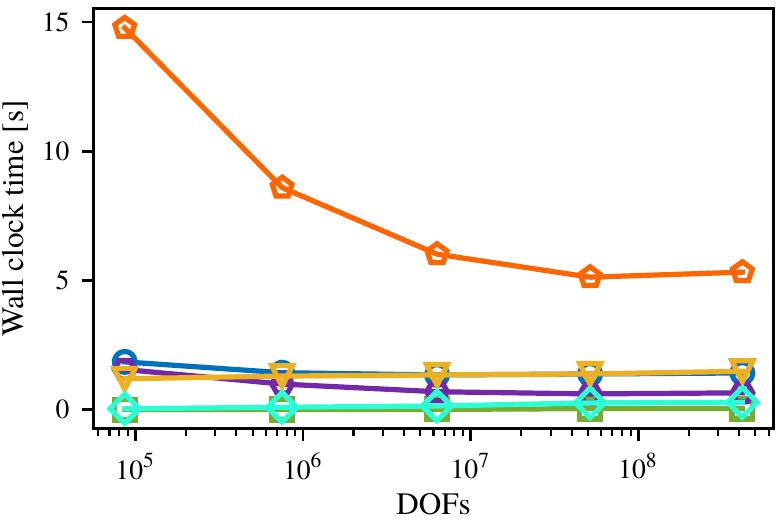}
    \caption{{\rev Swiss cheese.}}
    \label{fig:time-phases-c}
  \end{subfigure}
      \hspace{2em}
      \begin{subfigure}{0.45\textwidth}
    \includegraphics[scale=1]{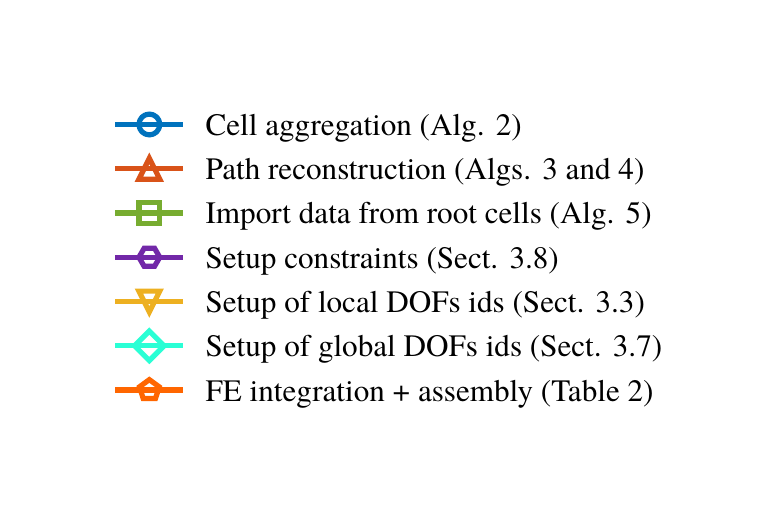}
  \end{subfigure}  
\caption{\secttitle: Wall clock time for the major phases of the \ac{agg} method versus problem size. Results computed with \ac{agfe} spaces, $\beta=10$,  {\rev $w=10$}, and load~3. }
\label{fig:time-phases}
\end{figure}

We conclude the numerical examples by reporting the wall clock time consumed in the main phases of the \ac{agg} method (see Fig.~\ref{fig:time-phases}) and the one in the linear solver step (see Fig.~\ref{fig:time-solver}). As it is observed in Fig.~\ref{fig:time-phases}, all phases  have remarkable weak scaling.
For the linear solver (see Fig.~\ref{fig:time-solver}), we report the wall clock times of the solver setup (i.e., the setup of the \ac{amg} preconditioner) and the wall clock time of the solver run (i.e., the preconditioned conjugate gradient iterations). Note that the scaling of the linear solver is not as optimal as for the different phases reported in Fig.~\ref{fig:time-phases}. However, it can still be  considered excellent since the solver wall clock times increase only moderately as the problem is scaled, and also if one takes into account that we have customized the linear solver as {\rev little} as possible. As an illustrative example, for the popcorn flake body and for load 3, the total solver wall clock time (setup plus iterations) scales from $0.36$ to $4.33$ seconds, while the problem size scales from   $67,192$ to $298,911,096$ \acp{dof}. That is to say, the total solver time increases by a factor of $12x$, while the problem size increases by a factor of $4,449x$, which are excellent results.


\begin{figure}[ht!]
  \centering

  \begin{subfigure}{0.45\textwidth}
    \includegraphics[scale=1]{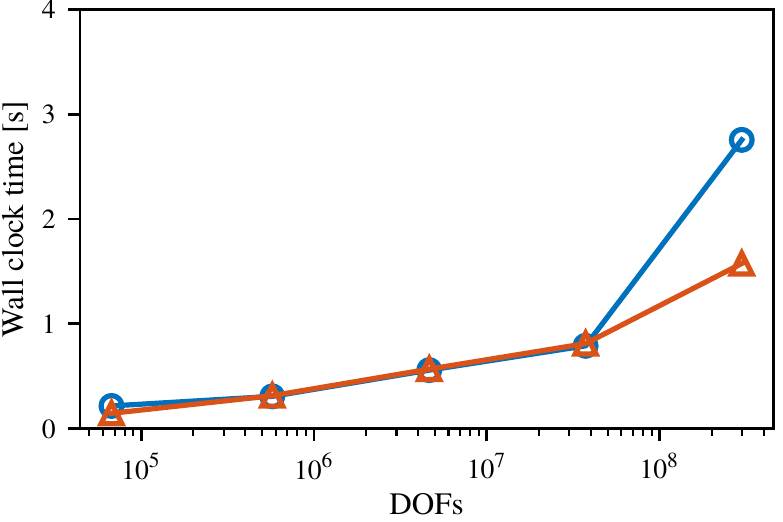}
        \caption{Popcorn flake.}
    \label{fig:time-solver-a}
  \end{subfigure}
  \hspace{2em}
  \begin{subfigure}{0.45\textwidth}
    \includegraphics[scale=1]{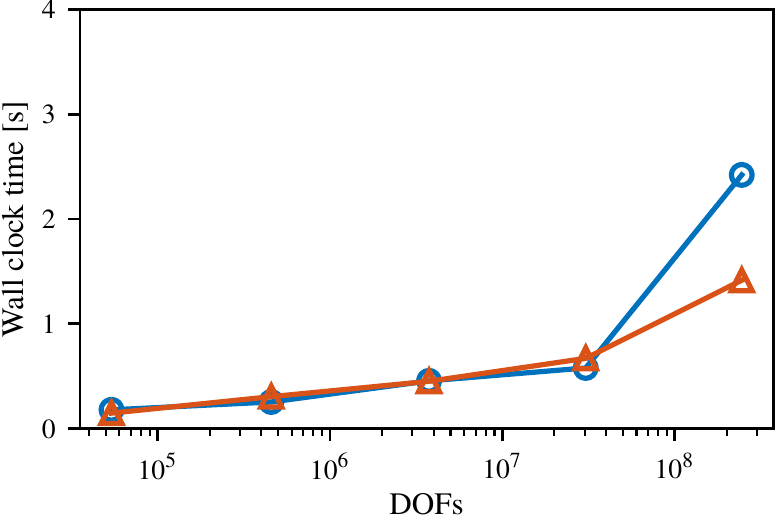}
    \caption{Spiral.}
    \label{fig:time-solver-b}
  \end{subfigure}

    \begin{subfigure}{0.45\textwidth}
    \includegraphics[scale=1]{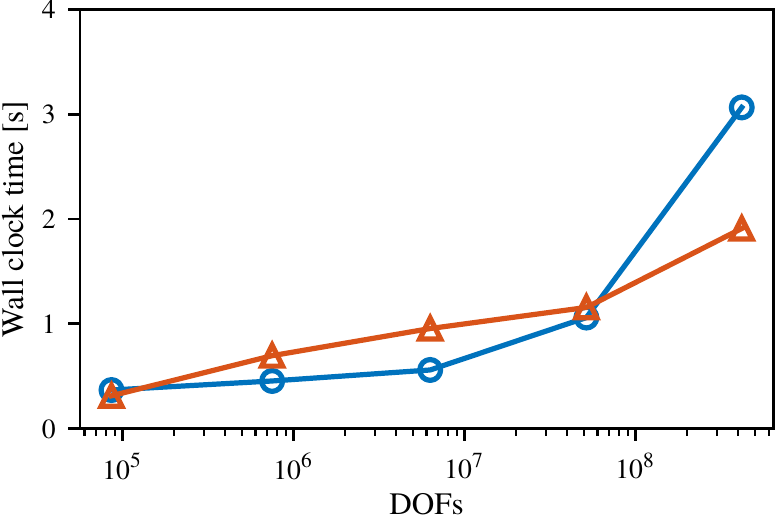}
    \caption{{\rev Swiss cheese.}}
    \label{fig:time-solver-c}
  \end{subfigure}
      \hspace{2em}
      \begin{subfigure}{0.45\textwidth}
    \includegraphics[scale=1]{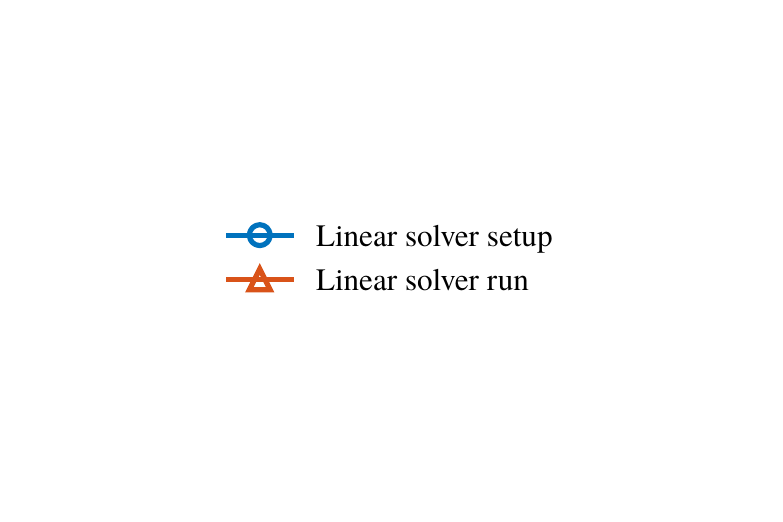}
  \end{subfigure}    
  
\caption{\secttitle: Wall clock time for the two main linear solver phases (setup and run) versus problem size. Results computed with \ac{agfe} spaces, $\beta=10$,  {\rev $w=10$}, and load~3.}
\label{fig:time-solver}
\end{figure}

\renewcommand{\secttitle}{Weak scaling test for the Poisson equation on the unit cube}
\subsection{\rev \secttitle}\label{sec:weak-scaling-test-full}

{\rev 

In order to demonstrate that the (moderate) degradation of the linear solver time observed in Fig. \ref{fig:time-solver} is not caused by the \ac{agg} method, we have performed (as a reference) a weak scaling test with a Poisson problem defined on the unit cube, with a body-fitted mesh optimally partitioned into sub-domains by means of a Cartesian partition. The test is done for three different local loads (referred to as "load 1", "load 2" and "load 3") with $20^3$, $30^3$, and $40^3$ cells per sub-domain respectively. The largest mesh considered has $1,053,696,000$ cells which are distributed over $16,464$ sub-domains (see Table \ref{tab:partitions-fullbox}). The linear solver considered here is the same as the one used before in Section~\ref{sec:weak-scaling-test} (i.e., we have used the configuration file displayed in Listing~\ref{lst:petscrc}). The results of this weak scaling test are reported in Fig.~\ref{fig:full-box}. Note that the number of iterations of the conjugate gradient solver exhibits the optimal scaling behavior expected for this problem type (see Fig.~\ref{fig:full-box-a}). However, even though this is the most favorable setup one can consider for an \ac{amg} preconditioner (i.e., Poisson problem, Cartesian body fitted meshes, optimal load balance, no \ac{agg} constraints), the obtained weak scaling of the solver time is not the expected optimal one (constant time as the problem size increases in the same proportion as the number of processors). This can be clearly seen in Figs.~\ref{fig:full-box-b} and \ref{fig:full-box-c}. These results show that the moderate loss of scalability previously observed also for the computations based on the \ac{agg} method in Fig.~\ref{fig:time-solver} cannot be directly attributed to the \ac{agg} method itself. The degradation of the linear solver time could be caused by other factors such as the underlying parallel system (i.e., the HPC cluster architecture plus the software stack that runs on top of it), but confirming this is out of the scope of the current work. In any case, by performing the weak scaling test on the unit cube with a body-fitted mesh, we have discarded that the loss of efficiency is caused by the \ac{agg} method.}


{\rev
\begin{table}[ht!]
\centering
\begin{small}
\begin{tabular}{rrrrr}
\toprule
 & \multicolumn{3}{c}{\rev Number of cells}\\
\cmidrule{2-4}
 {\rev $|\mathcal{S}|$} &\rev load 1 &\rev load 2&\rev load 3 \\
\midrule
\rev 48 &   \rev   384,000 & \rev 1,296,000 & \rev 3,072,000 \\
\rev 1,296 & \rev 10,368,000 & \rev 34,992,000 & \rev 82,944,000\\
\rev 3,072 & \rev 24,576,000 & \rev 82,944,000 & \rev 196,608,000\\
\rev 6,000 & \rev 48,000,000 & \rev 162,000,000 &\rev  384,000,000\\
\rev 10,368 & \rev 82,944,000 &\rev  279,936,000 & \rev 663,552,000\\
\rev 16,464 & \rev 131,712,000 & \rev  444,528,000 &\rev  1,053,696,000\\
\bottomrule
\end{tabular}
\end{small}
\caption{\rev \secttitle: Number of sub-domains ($|\mathcal{S}|$) and number of cells of the considered meshes.}
\label{tab:partitions-fullbox}
\end{table}
}

\begin{figure}[ht!]
  \centering

  \begin{subfigure}{0.45\textwidth}
    \includegraphics[scale=1]{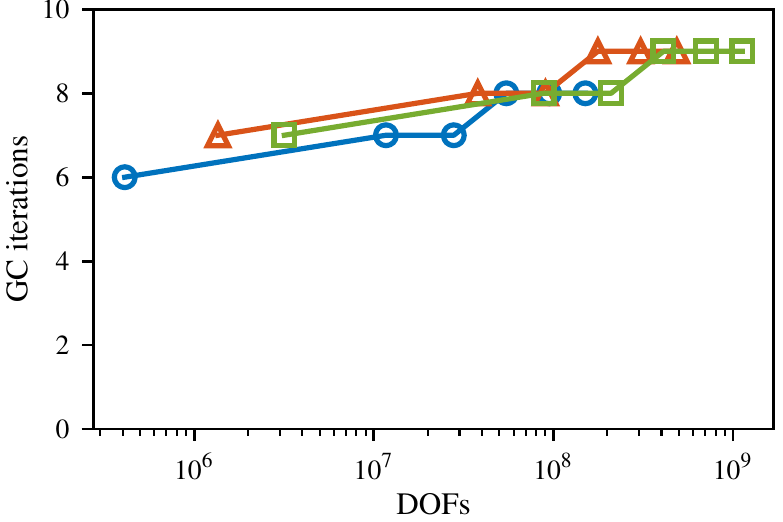}
        \caption{\rev Weak scaling of CG iterations.}
    \label{fig:full-box-a}
  \end{subfigure}
  \hspace{2em}
  \begin{subfigure}{0.45\textwidth}
    \includegraphics[scale=1]{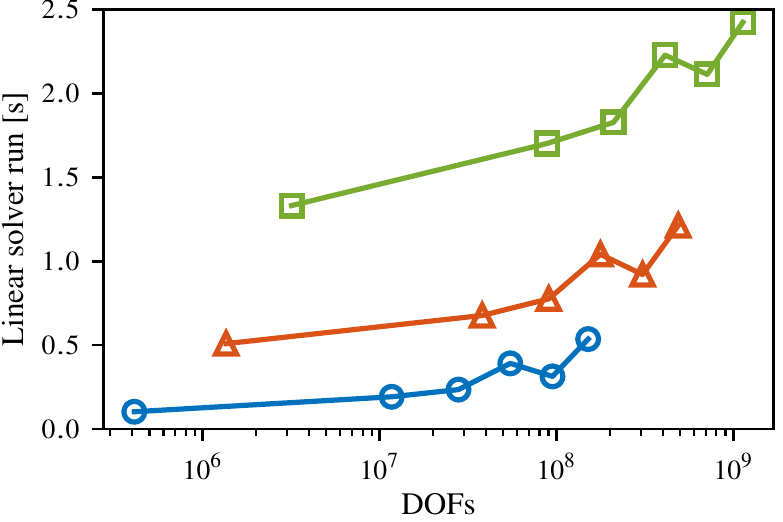}
    \caption{\rev Weak scaling of linear solver (iteration phase).}
    \label{fig:full-box-b}
  \end{subfigure}

    \begin{subfigure}{0.45\textwidth}
    \includegraphics[scale=1]{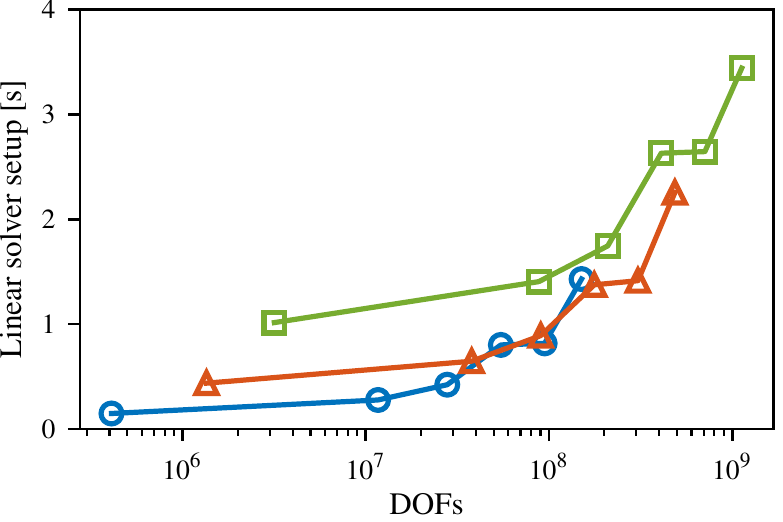}
    \caption{\rev Weak scaling of linear solver (setup phase).}
    \label{fig:full-box-c}
  \end{subfigure}
      \hspace{2em}
      \begin{subfigure}{0.45\textwidth}
    \includegraphics[scale=1]{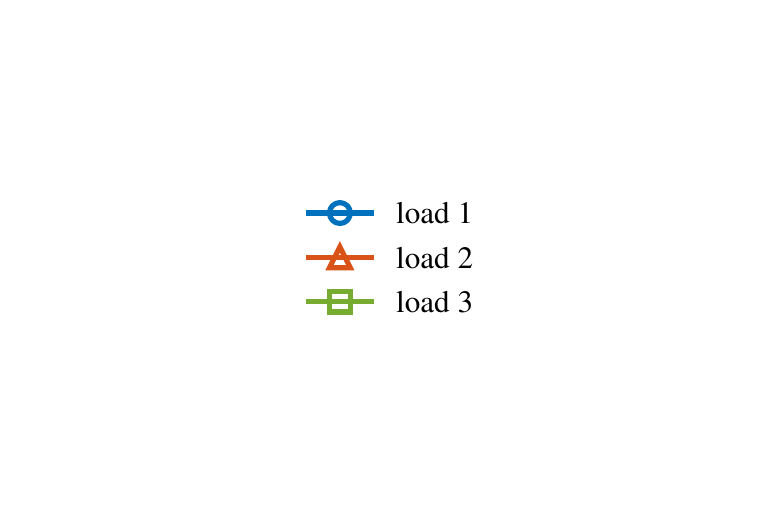}
  \end{subfigure}    
  
\caption{\rev \secttitle: Performance of the linear solver.}
\label{fig:full-box}
\end{figure}

\section{Conclusions} \label{sec:conclussions}

In this work, we have presented the distributed-memory implementation and performance of the so-called \ac{agg} method for the solution of large-scale problems with unfitted \ac{fe} techniques. The method is based on removal of basis functions associated with badly cut cells by introducing carefully designed constraints, which result in well-conditioned systems of linear algebraic equations independently of the position of the cuts. As the main contributions of the paper, (a) we have proposed a parallel implementation of \ac{agg}, (b) shown that the proposed parallel setup is scalable, and (c) shown that, by using \ac{agg}, the underlying system of linear equations can be efficiently solved using preconditioned iterative linear solvers designed for conventional \ac{fe} analysis. In this case, we have considered well-known \ac{amg} preconditioners available in \petsc. This is in contrast to previous works, which considered highly customized preconditioners in order to handle the severely ill-conditioned operators resulting from the discretization based on unfitted grids. Also in contrast to previous works, which mainly consider non-scalable serial algorithms, we have studied the performance and scalability of the method with weak scaling tests up to 16K processors and up to 300M \acp{dof}. The proposed parallel \ac{agg} method can be easily incorporated in existing large-scale codes since (a) it allows one to use standard parallel linear solvers, and (b) the algorithmic phases of the method can be implemented using standard functionality available in distributed-memory \ac{fe} codes. In this case, the parallel implementation of the AgFEM method has been done in the large-scale \ac{fe} package \FEMPAR.

\section*{Acknowledgments}

\thethanks

\begin{small}

\setlength{\bibsep}{0.0ex plus 0.00ex}
\bibliographystyle{myabbrvnat}
\bibliography{art034}

\end{small}

\end{document}